\input amssym.def
\input amssym.tex
\magnification=\magstep1
\font\twelvebf=cmbx12
\thinmuskip = 2mu
\medmuskip = 2.5mu plus 1.5mu minus 2.1mu  
\thickmuskip = 4mu plus 6mu
\font\teneusm=eusm10
\font\seveneusm=eusm7
\font\fiveeusm=eusm5
\newfam\eusmfam
\textfont\eusmfam=\teneusm
\scriptfont\eusmfam=\seveneusm
\scriptscriptfont\eusmfam=\fiveeusm

\font\tenmib=cmmib10
\font\sevenmib=cmmib7
\font\fivemib=cmmib5
\newfam\mibfam
\textfont\mibfam=\tenmib
\scriptfont\mibfam=\sevenmib
\scriptscriptfont\mibfam=\fivemib

\font\tenss=cmss10
\font\sevenss=cmss8 scaled 833
\font\fivess=cmr5
\newfam\ssfam
\textfont\ssfam=\tenss
\scriptfont\ssfam=\sevenss
\scriptscriptfont\ssfam=\fivess
\def\ss{\fam\ssfam}
\font\mathnine=cmmi9
\font\rmnine=cmr9
\font\cmsynine=cmsy9
\font\cmexnine=cmex10 scaled 905
\def\msmall#1{\hbox{$\displaystyle \textfont0=\rmnine
\textfont1=\mathnine \textfont2=\cmsynine \textfont3=\cmexnine
{#1}$}}
\hyphenation{Lip-schit-zian Lip-schitz}
\def\cc{{\Bbb C}}
\def\jj{{\Bbb J}}
\def\rr{{\Bbb R}}
\def\nn{{\Bbb N}}
\def\pp{{\Bbb P}}
\def\zz{{\Bbb Z}}
\def\cinf{C^\infty}
\def\loc{{\rm loc}}
\def\st{_{\rm st}}
\def\coker{{\ss Coker}\,}

\def\exp{{\ss exp}}
\def\ssh{{\ss H}}
\def\hom{{\ss Hom}\,}
\def\homr{{\ss Hom}\vph_\rr}
\def\id{{\ss Id}}
\def\im{{\ss Im}\,}
\def\ker{{\ss Ker}\,}
\def\lim{\mathop{\ss lim}}

\def\ord{{\ss ord}}
\def\supp{{\ss supp}\,}

\def\pr{{\ss pr}}

\def\epsi{\varepsilon}
\def\bss{\backslash}
\def\ogran{{\hskip0.7pt\vrule height8.5pt depth4.5pt\hskip0.7pt}}
\def\comp{\Subset}
\def\d{\partial}
\def\dbar{\overline\partial}
\def\dbarj{\dbar_J}
\def\ddef{\mathrel{{=}\raise0.23pt\hbox{\rm:}}}
\def\deff{\mathrel{\raise0.23pt\hbox{\rm:}{=}}}
\def\ge{\geqslant}
\def\inv{^{-1}}
\def\fraction#1/#2{\mathchoice{{\msmall{ #1\over#2}}}%
{{ #1\over #2 }}{{#1/#2}}{{#1/#2}}}

\def\le{\leqslant}
\def\vph{^{\mathstrut}}
\def\lrar{\longrightarrow}
\def\emptyset{\varnothing}
\def\scirc{\mathchoice{\msmall{\circ}}{\msmall{\circ}}%
{{\scriptscriptstyle\circ}}{{\scriptscriptstyle\circ}}}
\def\state#1. {\medskip\noindent{\bf#1. }}
\def\qed{\par\nobreak\hfill\ \hbox{q.e.d.}}
\def\Chi{\raise 2pt\hbox{$\chi$}}
\def\eg{\hskip1pt plus1pt{\sl{e.g.\/\ \hskip1pt plus1pt}}}
\def\ie{\hskip1pt plus1pt{\sl i.e.\/\ \hskip1pt plus1pt}}
\def\iff{\hskip1pt plus1pt{\sl iff\/\hskip1pt plus1pt }}

\def\sli{{\sl i)} }
\def\slii{{\sl ii)} }
\def\sliii{{\sl iii)} }
\def\sliv{{\sl iv)} }
\def\slv{{\sl v)} }
\def\slvi{{\sl vi)} }

\def\diff{{\cal D}\mskip -.6mu i\mskip -5.1mu f\mskip-6mu f}
\def\star{\mathop{\msmall{*}}}
\def\calu{{\cal U}}
\def\calm{{\cal M}}
\def\calc{{\cal C}}
\def\calh{{\cal H}}
\def\caln{{\cal N}}
\def\calv{{\cal V}}

\def\epsi{\varepsilon}
\def\bss{\backslash}
\def\ogran{{\hskip0.7pt\vrule height8.5pt depth4.5pt\hskip0.7pt}}
\def\comp{\Subset}
\def\d{\partial}
\def\dbar{\overline\partial}
\def\dbarj{\dbar_J}
\def\ddef{\mathrel{{=}\raise0.23pt\hbox{\rm:}}}
\def\deff{\mathrel{\raise0.23pt\hbox{\rm:}{=}}}
\def\ge{\geqslant}
\def\inv{^{-1}}
\def\fraction#1/#2{\mathchoice{{\msmall{ #1\over#2}}}%
{{ #1\over #2 }}{{#1/#2}}{{#1/#2}}}

\def\le{\leqslant}
\def\vph{^{\mathstrut}}
\def\lrar{\longrightarrow}
\def\emptyset{\varnothing}
\def\scirc{\mathchoice{\msmall{\circ}}{\msmall{\circ}}%
{{\scriptscriptstyle\circ}}{{\scriptscriptstyle\circ}}}
\def\state#1. {\medskip\noindent{\bf#1. }}
\def\qed{\par\nobreak\hfill\ \hbox{q.e.d.}}

\ifx \twelvebf\undefined \font\twelvebf=cmbx12\fi
\centerline{\twelvebf Pseudo-holomorphic curves and
envelopes of meromorphy}
\smallskip\centerline{\twelvebf of two-spheres in $\cc\pp^2$.}

\smallskip \smallskip

\centerline{\rm S.~Ivashkovich\footnote{$^*$}{Supported in part by
Schwerpunkt "Komplexe Mannigfaltigkeiten".},
V.~Shevchishin\footnote{$^{**}$}{Supported in part by SFB-237.}
}

\smallskip
\centerline{January  1995.}

\bigskip
\def\longpoints{\leaders\hbox to 0.5em{\hss.\hss}\hfill \hskip0pt}
\noindent \bf Table of  contents.

\smallskip\smallskip\noindent\sl
\line{0. Introduction. \longpoints pp.2--7.}

0.1. Envelopes of meromorphy.

0.2. Sketch of the proof.

0.3. Smoothness of the~moduli space in a neighbourhood of a cusp-curve.

0.4. The Genus formula.

\smallskip\noindent
\line{1. Deformations. \longpoints pp.7--9.}

1.1. Deformation of structure.

1.2. The Genus formula for immersed symplectic surfaces.

\smallskip\noindent
\line{2. The $D_{u, J}$-operator and its properties. \longpoints pp.10--22.}

2.1. The definition of $D_{u, J}$-operator.

2.2. Operator $\dbar_{u, J}$ and a holomorphic structure on the~induced
bundle.

2.3. Surjectivity of $D_{u, J}$.

2.4. Deformation of cusp-curves.

\smallskip\noindent\setbox1=\hbox{3. }
3. Singular points of pseudo-holomorphic curves and positivity of geometric \par
\line{\hskip\wd1%
self-intersections. \longpoints pp.22--32.}

3.1. Unique continuation lemma.

3.2. Inversion of $\dbar_J +R$.

3.3. Perturbing cusps of pseudo-holomorphic curves.

3.4. Positivity of geometric self-intersection.

\smallskip\noindent
4. The Bennequin index and the Genus formula for pseudo-holomorphic \par
\line{\hskip\wd1%
curves. \longpoints pp.32--39.}

4.1. Transversality and the Bennequin index of the cusp.

4.2. The Genus formula and the lower estimate of the Bennequin index.

\smallskip\noindent
\line{5. Continuity principles and envelopes. \longpoints pp.39--47.}

5.1. Continuity principles relative to the K\"ahler spaces.

5.2. Construction of the family of curves and the~proof of the~Theorem 1.

5.3. Examples.

\smallskip\noindent
\line{6. Appendix. \longpoints pp.47--52.}

6.1. A complex structure lemma.

6.2. A matching structures lemma.

6.3. Gromov topology and versal deformations of noncompact

curves.

\smallskip\noindent
\line{References. \longpoints pp.52--54.}

\break
\noindent\bf 0. Introduction.

\smallskip\noindent
\sl 0.1. Envelopes of meromorphy.

\rm Consider a complex surface $X$, {\sl i.e.} a~two-dimensional complex
manifold, which is connected, Hausdorff and countable at infinity. Suppose
that $X$ admits a K\"ahler form $\omega $. Further let $u:S^2\to
X$ be an immersion  of a standard two-dimensional sphere into $X$. We shell
call this immersion  $\omega $-positive (or $\omega $-symplectic) if
$u^*\omega$ is strictly positive 2-form on $S^2$ with respect to the standard
orientation of $S^2$. We suppose moreover, that
self-intersections of $M\deff u(S^2)$ are transversal and positive. We denote
by $\delta$ the~number of pairs $s_1\not= s_2\in S^2$ with $u(s_1)=u(s_2)$.

Denote by $c_1(X)$ a 2-form on $X$ representing the first Chern class
of $X$. For a~closed immersed two-surface $M$ in $X$ denote by $c_1(X)\cdot
[M]$ the~value of the~first Chern class on $[M]$, {\sl i.e.} the~number
$\int_Mc_1(X)$. Let $U$ be an open neighbourhood of $u(S^2)=M$. Recall that the
envelope of meromorphy of open set $U$ is a maximal domain $(\hat U,\pi )$
over $X$ with an inclusion $i:U\longrightarrow \hat U$ such that every
meromorphic function $f$, which is defined on $U$, extends to a meromorphic
function $\hat f$ on $\hat U$. Such envelope always exists, see $\S 5$. We
shall say that the~domain $(\hat U,\pi )$ over $X$ is unbounded if $\pi (\hat
U)$ is not relatively compact in $X$.

In this paper we call a K\"ahler surface $(X,\omega)$ positive if for any
almost complex structure $J$ on $X$ tamed by $\omega$ (see definition below)
 any $J$-holomorphic sphere $C$ satisfies $[C]^2 \ge 0$.

\smallskip
Our first result here is the~following

\smallskip
\noindent
\bf Theorem 1. \it Let $M$ be a $\omega $-positive immersed  two-sphere in a
positive K\"ahler surface $(X,\omega )$, having only positive
self-intersections. Suppose that $c_1(X)\cdot [M] - \delta =p\ge 1$. Then the
envelope of meromorphy $(\hat U,\pi )$ of any neighbourhood $U$ of $M$ is
either unbounded or contains a one-dimensional compact analytic set $C$ such
that:

\smallskip
$(1)$ $\pi (C) = \bigcup _{k=1}^NC_k$ is a union of rational curves.

\smallskip
$(2)$ $\displaystyle \sum _{k=1}^N \left(c_1(X)\cdot [C_k]
-\delta_k -\varkappa_k \right) \ge p$.

\smallskip
\rm Here $\delta_k$ and $\varkappa_k$ are nonnegative integers which we define
below. $\delta_k$ is a number of nodes of $C_k$ taken with appropriate
multiplicities, and $\varkappa_k$ is a sum of conductors of the~cusps of $C_k$.

\bigskip
\rm Consider for example the case $X = \cc\pp^2$ with the usual Fubini-Studi
form $\omega _{FS}$. Then $\ssh_2(\cc\pp^2,\zz)=\zz$ and for any surface $M
\sim n\cdot \cc\pp^1$ one has $[M]^2=n^2 \ge0$. So $\cc\pp^2$ is positive in
our sense. Also, as it is well known $c_1(\cc\pp^2) = 3\cdot \omega_{FS}$.
Remark further, that any \sl imbedded \rm two-sphere in $\cc\pp^2$ is
homologous to zero, $\pm \cc\pp^1$ or $\pm 2\cc\pp^1$, in the third case it
means to the quadric. After reversing the orientation, if necessary, we can
rest on the case with +. Further, since $\omega$-positivity implies $\int
_S\omega>0$, the~condition on the~Chern class is now satisfied automatically.

If $M$ is symplectically immersed in $\cc\pp^2$ with positive
self-intersections, then the~condition $c_1(\cc\pp^2)[M] > \delta$ means that
$$
c_1(\cc\pp^2)[M] > \delta = { [M]^2 - c_1(\cc\pp^2)\cdot [M]
\over2} + 1,\eqno(0.1)
$$
(see the~Genus formula for symplectic surfaces in {\sl Lemma 1.3.1}). Thus
$(0.1)$ can be rewritten as $3m > (m^2-3m)/2 +1$ where $[M]=m[\cc\pp^1]$.

\smallskip
So one has the~following

\smallskip\noindent
\bf Corollary 1. \it The~envelope of meromorphy of any immersed symplectic
sphere in $\cc\pp^2$ of degree not more than 8 with only positive
self-intersections contains a rational curve and thus, in fact coincides with
$\cc\pp^2$ itself. In particular, this is true for every symplectically
embedded sphere in $\cc\pp^2$.

\smallskip\noindent\bf
Remark. \rm  1. As examples of symplectic surfaces  in K\"ahler surfaces one
can consider ``$C^1$-small perturbations" of holomorphic  curves, because the
condition  to be symplectic ($\omega $-positive) is obviously open.

\noindent
2. Another example, which is discussed in $\S 5$, shows that the envelope of
meromorphy of every imbedded symplectic sphere in $\cc\pp^1\times \cc\pp^1$
with natural K\"ahler structure contains a rational curve, which is a graph of
a rational function.

\noindent
3. In {\sl 5.3} the surfaces of higher genus are also discussed. In fact the
envelope of meromorphy of any symplectic surface in $\cc\pp^2$ of degree
$\le $8 (i.e. of genus $\le $21) coincides with $\cc\pp^2$, see {\sl Corollary }
5.3.2.

\noindent
4. In fact we prove more general statement than that is given in
{\sl Theorem 1}. Namely, one can replace meromorphic functions ({\sl i.e.}
meromorphic mappings into $\cc\pp^1$) by meromorphic mappings with values
into arbitrary  disk-convex K\"ahler spaces, see {\sl 5.2}.

\noindent
5. As a corollary of the {\sl Theorem }1 we obtain a result on rigidity of
symplectic imbeddings of two-spheres in $\cc\pp^2$ under biholomorphisms,
see {\sl Corollary }5.2.2.

We want to point out the difference of the situation studied here with that
of 2-spheres in $\cc^2$, see [B-K], [Sh]. Our spheres are \sl symplectic, \rm
in particular, they are not homologous to zero.

\medskip
\noindent \sl 0.2. Sketch of the proof.

\rm We give now a short sketch  of the proof of {\sl Theorem 1}.

\smallskip
\noindent {\sl Step 1: Deformation of a structure.} Using $\omega$-positivity
of $M$ we construct a smooth curve $\{ J_t\}_{t\in [0, 1]}$ of almost complex
structures on $X$ satisfying the following properties:

a) $J_0$ is the~given integrable structure on $X$.

b) $\{x\in X:  J_t(x)\not= J_0(x)\}\subset U_1\subset
\subset U$ for each $t\in [0, 1]$, {\sl i.e.} each $J_t$ is different from
$J_0$ only
inside some relatively compact open subset $U_1$ (the~same for all $J_t$)
of the  given neighbourhood $U$ of $M$.

c) $M$ is $J_1$-holomorphic.

d) All $\{ J_t\} $ are "tamed" by our K\"ahler form $\omega$, {\sl i.e.}
$\omega(v, J_tv)>0 $ for every nonzero $v\in TX$.

\smallskip \noindent
{\sl Step 2: Deformation of the~sphere $M$.} We construct a "piece-wise
continuous" family of (possibly reducible) surfaces $M_t$ such that:

a) $M_t$ is $J_t$-holomorphic.

b) $M_1 = M$.

c) For all $t\in [0, 1]$   $M_t$ is a finite union of $J_t$-holomorphic
 spheres $\{ M_t^j \} $ with
$$
\sum c_1(X)[M_t^j]-\delta_j -\varkappa_j \ge p.
$$

d) The family is "piece-wise continuous" in the following sense. There is a
partition of unit interval $0=t_0<t_1<\ldots<t_r=1$, natural numbers $N_1\ge
\ldots \ge N_r= 1$, and $J_t$-holomorphic maps $u_t:\bigsqcup_{N_i} S^2 \to X$
for $t\in (t_{i-1}, t_i]$, such that $u_t$ is continuous in $t\in (t_{i-1},
t_i]$, $M_t\deff u_t(\bigsqcup _{N_i}S^2)$ converges
in the~Hausdorff metric when $t\searrow t_i$, and ${\cal H}$-$lim_ {t\searrow
t_{i-1}}M_t$ contains $M_{t_{i-1}}$.

Note that by ``tameness" of all $J_t$ and closeness of $\omega$ the ${\ss area}
(M_t)\sim \int _{M_t}\omega$ is bounded uniformly for all $t$. So by Gromov's
compactness theorem a~Hausdorff limit of $M_t$ exists and is a union of
$J_{t_{i-1}}$-holomorphic spheres.

\smallskip
\noindent {\sl Step 3: Kontinuit\"atssatz}. Let $(\hat U,\pi )$ be the
envelope of
$U\supset M$. Note that $\Sigma _t = M_t \setminus U_1$ is a holomorphic curve
in $X$ with boundary on $\partial U_1$, and $\Sigma _1 = \emptyset $. So
by the ``continuity principle", see $\S 5$ for details, $\Sigma _t$ and thus
$M_t$
could be lifted onto the envelope $\hat U$. But $M_0$ is holomorphic and
satisfies (c) from {\sl Step 2}. Thus the proof is complete.

For the construction of a  family of curves in the {\sl Step 2} we need two
results about pseudo-holomorphic curves.

\bigskip \noindent
\sl 0.3. Smoothness of the~moduli space in a neighbourhood of a cusp-curve.
\rm

Let $(X, J)$ be an almost complex manifold of a complex dimension $n$. An
almost complex structure is always supposed to have smoothness of class $C^1$.
Recall that an (irreducible) pseudo-holomorphic (or more precisely
$J$-holomorphic) curve in $X$ is an image of a (connected) Riemannian surface
$(S, J_S)$ in $X$ under a $C^1$-mapping $u:S\to X$,
which satisfies the equation

\smallskip
$$
du + J\circ du\circ J_S = 0\eqno(0.2)
$$
\smallskip
In the case of integrable $J$ it is just a Cauchy-Riemann equation.
If $J$ is $C^{k,\alpha}$-smooth with $k\ge1$ and $0<\alpha<1$, then
$u$ is $C^{k+1,\alpha}$-smooth due to elliptic regularity.

Denote by $[\gamma ]$ some homology class in $\ssh_2(X,\zz)$. Fix a compact
Riemannian surface $S$, {\sl i.e.} a compact, connected, oriented, smooth
manifold of real dimensions 2. Fix $p$ with $2<p<\infty$ and consider
a~Banach manifold
$$
{\cal S} =\bigl \{ u \in L^{1,p}(S, X):u(S)\in [\gamma ]\bigr\}
$$
of all $L^{1,p}$-maps from $S$ to $X$, representing the~class $[\gamma]$,
(see \S2 for details).

Denote by ${\cal J}$ the Banach manifold of $C^1$-smooth almost
complex structures on $X$, and by ${\cal J}_S$ the Banach manifold of
$C^1$-smooth almost complex structures on $S$. Consider also a subset
${\cal P}\subset {\cal S}\times {\cal J}_S \times {\cal J}$ consisting of an
all triples $(u, J_S, J)$ with $u$ being $(J_S, J)$-holomorphic, {\sl i.e.}

\smallskip
$$
{\cal P} = \bigl\{(u, J_S, J)\in {\cal S}\times {\cal J}_S\times {\cal J}: du +
J\circ du\circ  J_S = 0 \bigl\} \eqno(0.3)
$$

\smallskip
Denote by $\pr_{\cal J} : {\cal S}\times {\cal J}_S\times {\cal J}\to
{\cal J}$ the natural projection onto the third factor.

We are interested when $\pr_{\cal J}:\cal P \to {\cal J}$ is a surjective
map, or more precisely, when for a given pseudo-holomorphic curve $u:(S,J_S)
\to (X,J)$ one can find a neighbourhood $V\subset {\cal J}$ of $J$ and
a~continuous section of a fibration $\pr_{\cal J} :{\cal P} \to {\cal J}$ over
$V$, which passes through $(u,J_S,J)\in {\cal P}$. In other words, when for
every small perturbation $\tilde J$ of $J$ one can find a small perturbation
$\tilde u$ of $u$, which is $\tilde J$-holomorphic.

We show that for the~pseudo-holomorphic curve $u:(S,J_S) \to (X,J)$
the~pulled-back bundle $E\deff u^*TX$ possesses the~natural {\sl holomorphic}
structure, such that the~differential $du:TS\to E$ is a {\sl holomorphic}
homomorphism. This allows us to define the order of vanishing of
the~differential
$du$ at a~point $s\in S$. We denote this number by ${\ss ord}_sdu$.

Let ${\cal O}(N_0)$ be a free part of the~quotient
${\cal O}(E)/du({\cal O}(TS))$, {\sl i.e.}
$$
0\longrightarrow {\cal O}(TS) \buildrel du \over\longrightarrow
{\cal O}(E) \buildrel \pr \over\longrightarrow
{\cal O}(N_0)\oplus {\cal N}_1 \longrightarrow 0,
$$
where ${\cal N}_1$ is supported on the~finite set of cusps of $u$ ({\sl i.e.}
points of vanishing of $du$). On the~space $L^{1,p}(S,N_0)$ of sections
of the~bundle $N_0$ the~natural Gromov operator $D_N: L^{1,p}(S,N_0) \to
L^p(S,\Lambda^{0,1}S \otimes N_0)$ is defined. We prove
the~following

\medskip\nobreak\noindent
\bf Theorem 2. \it Suppose that $D_N$ is surjective. Then

\nobreak
\sli for some neighbourhood $U\subset {\cal S} \times {\cal J}_S
\times {\cal J}$ of $(u, J_S, J)$ the set ${\cal P} \cap U$ is a Banach
submanifold of $U$;

\slii there exists a $C^1$-map $f$ from some neighbourhood $V$ of
$J\in {\cal J}$ into $U$ with $f(V)\subset {\cal P}$ and $\pr_{\cal J}
\scirc f = \id_V$, such that $f(J)=(u, J_S, J)$.
\rm

\medskip\noindent
\bf Remark. \rm When $X$ is a surface and $N_0$ is a~line bundle
the~sufficient condition for surjectivity of $D_N$ is $c_1(N_0) > 2g-2$,
where $g$ is the~genus of $S$, see [Gro], [Hf-L-Sk] and {\sl Corollary 2.3.3.}
So, if $c_1(E) > \sum_{s\in S}{\ss ord}_sdu $,
one has the surjectivity
of projection $\pr_{\cal J}$.

It is important for us that such a map $f$ still exists if even $u$ has cusps,
but the sum of their orders is strictly bounded by the first Chern class of
the~induced bundle. This result extends the result of [Hf-L-Sk] from the
case of immersed curves to the case of pseudoholomorphic curves with cusps.

In fact one has also the next result about Moduli space of nonparametrized
pseudo holomorphic curves, see 2.4.

\state Corollary 2. \it  Let $u: (S, J_S) \to (X, J)$ be a nonconstant
irreducible pseudo-holomorphic map, such that $\ssh^1_D(S, N_0)=0$. Than in
a neighbourhood of $M\deff u(S)$ the Moduli space of {\sl nonparametrized $J$-
holomorphic curves}
${\cal M}_{[\gamma], g, J}$ is a~manifold with the tangent space $T_M
{\cal M}_{[\gamma], g, J} = \ssh^0_D(S, N_0) \oplus \ssh^0(S, {\cal N}_1)$.

\bigskip
\noindent \sl 0.4. The Genus formula.

\rm We want now to state our third result. It concerns with the version of so
called genus formula for pseudo-holomorphic curves in almost complex
surfaces.

By the geometric self-intersection number of pseudo-holomorphic curve $M =
u(S)$ we understand the number of pairs $(s_1, s_2)$ of distinct points of $S$
such that $u(s_1)=u(s_2)$, taken with appropriate multiplicities. We call this
points {\sl nodes} of $M$. We prove that there only a finite number of such
pairs and multiplicities of intersections are always positive, see $\S 3$.
Defined in such way the~geometric self-intersection number we denote by
$\delta$. This definition in an obvious way extends to the~case of reducible
curves, {\sl i.e.} curves $M$ which are a union of a finite number of
irreducible ones.

We show that for the $J$-holomorphic curve $u : S\to X$ with $J\in C^1$
the order of zero of differential ${\ss ord}_adu $ at point $a\in S$
is well defined. The points $a$ with ${\ss ord}_a du \ge1$ are called {\sl
cusp-points} of $M$. A compact pseudo-holomorphic curve have only finitely
many cusp-points. We define in $\S 4$, using the Bennequin index, the number
$\varkappa _i$ for a cusp point $a_i$ of a pseudo-holomorphic curve, and show
that the number $\varkappa_i$ is bounded from below by ${\ss ord}_{a_i}du$,
see {\sl Corollary 4.2.2}. Following classical terminology we call the~number
$\varkappa_i$ the~{\sl conductor} of a cusp-point $a_i$. In the~holomorphic
case it is also called the~Milnor number of $a_i$. The sum of $\varkappa_i$
over all cusp points is denoted by $\varkappa$. In the case of reducible curve
one should take a sum over all irreducible components.

As usually, we denote by $[M]^2$ the homological self-intersection of a surface
in four-dimensional manifold. Having $J$, we can define a first Chern class
$c_1(X, J)$ which, in fact, does not change when $J$ varies continuously.
Therefore we shall usually omit the dependence $c_1(X)$ on $J$. Now we are
ready to write down the Genus formula.

\smallskip
\noindent
\bf Theorem 3. \rm(Genus formula). \it Let $M = \bigcup _{i=1}^dM_i$ be
a~pseudo-holomorphic curve in an almost complex surface $(X, J)$ and $g_i$
the~genera of its irreducible components $M_i$. Suppose that $J$ is of class
$C^1$ and that $M_i$ and $M_j$ are distinct for $i\not=j$. Then

\smallskip
$$
\sum_{i=1}^d g_i= {[M]^2 - c_1(X)[M]\over2} + d - \delta - \varkappa.
\eqno(0.4)
$$

\smallskip\noindent
\bf Remark. \rm
1. In the case of a holomorphic curve in a~complex surface (0.4) becomes the~
usual genus formula, or the~adjunction formula, or the~Pl\"ucker formula, well
known in classical algebraic geometry, see for example [Gr-Hr] or [Wa].

\noindent
2. When the structure $J$ is $C^\infty$-smooth the following important result
is due to McDuff: \sl if $M$ is an irreducible $J$-holomorphic curve of genus
$g$ in almost complex surface $(X, J)$, then ${1\over2}([M]^2-c_1(X)[M]) + 1
\ge g$, and the equality take place \iff $M$ is smoothly imbedded,\rm see
[McD-1].

\noindent
3. Recently, the result similar to our {\sl Theorem 3} was proved in [Mi-Wh],
assuming $C^2$-differentiability of $J$ using different methods and another
description of the numbers $\varkappa_i$.

We develop here the methods, which enable us to study the local and global
properties of pseudo-holomorphic curves under the assumption of $C^1$-
differentiability of an almost-complex structure $J$.

Now we want to sketch the construction of the family of curves, needed for
the~{\sl Step 2} of the~proof of {\sl Theorem 1}. We suppose for the~sake of
simplicity that $M$ is imbedded. Start with $M_1$ and use the condition
$c_1[M]\ge 1$. This allows us to apply  Theorem 2 and construct our family
of curves for $t$ in a neighbourhood of $1$ in $[0,1]$. Now let $t$ fall down
to $t_{r-1}$. The Genus formula goes out in the construction in the following
way. Because $M_t$ is an embedded smooth pseudo-holomorphic sphere for $t\in
(t_{r-1}, 1]$, the Genus formula tells us that

\smallskip
$$
0 = {[M_t]^2 - c_1(X)[M_t]\over2} + 1 \eqno(0.5)
$$

\smallskip
As we have already explained, $M_t$ have uniformly bounded areas. Thus the
Gromov compactness theorem provides that ${\cal H}$-$lim_{t\searrow t_{r-1}}
M_t = {\overline M}_{t_{r-1}}$ breaks into say $d$ spheres $M_1,\ldots, M_d$.
This system of spheres is connected as a~Hausdorff limit of connected sets, so
it has at least $d-1$ points of geometric self-intersection. Suppose in this
sketch of proof that multiplicities does not occur, {\sl i.e.} all $M_1,
\ldots, M_d$ are distinct. Now the Genus formula tells us once more that

\smallskip
$$
0 = {[{\overline M}_{t_{r-1}}]^2 - c_1(X)[\overline M_{t_{r-1}}]\over2} + d -
\delta - \varkappa \eqno(0.6)
$$

\smallskip
But the fractions in (0.5) and (0.6) are equal, because they depend only
on homology classes. So we have that $1 = d -\delta - \varkappa $, where
$\delta$ is at least $d-1$ as we explained. So $\varkappa = 0$ and $\delta
=d-1$. So all irreducible components $M_i$ of ${\overline M}_{t_{r-1}}$ are
smooth imbedded pseudo-holomorphic spheres again. Now since $\sum_{i=1}^d
c_1(X)[M_i] = p\ge 1$ we can choose from $\{ M_i\}$ a subset $\{M_{i_1},
\ldots, M_{i_k}\} $ such that $c_1(X)[M_{i_l}]\ge 1$ for $l=1,\ldots, k$
and $\sum_{l=1}^kc_1(X)[M_{i_l}]\ge p$. Repeating this procedure gives us
the construction of the family. At the end we want to point out that
the genus formula allows us to prove the closeness of the set of $t\in
[0, 1]$ for which such family exists. For its openness one needs {\sl
Theorem 2} and thus the surjectivity of the $D^N_{u,J}$-operator, which is
associated with a $J$-holomorphic sphere $u$. If $M$ is immersed and
multiplicities can occur one should be more careful and employ
the~positivity condition on the~K\"ahler surface $(X,\omega)$.

\smallskip
The authors would like to express their gratitude to the A.Vitushkin, who
asked us the question about the meromorphic envelopes of a "small
perturbation" of a line at infinity in $\cc\pp^2$. We also want to give
our thanks to J.-C.~Sikorav, for the usefull conversations and remarks.

\bigskip\bigskip\noindent
\bf  1. Deformations. \rm

\smallskip
In this paragraph we shall describe the {\sl Step 1} of the proof of the {\sl
Theorem 1}. Also, we shell make some preparatory steps for the proof of the
Genus formula.

\smallskip\noindent
\sl 1.1. Deformation of structure.

\rm First  we recall some elementary
facts about orthogonal complex structures in $\rr^4$.

In $\rr^4$ with coordinates $x_1, y_1, x_2, y_2$ consider the~standard
symplectic form $\omega = dx_1\wedge dy_1 + dx_2\wedge dy_2 $ and the~standard
complex structure $J\st$ defined by the operator

\medskip
$$
J\st =
\left( \matrix{0&-1&0&0\cr
1&0&0&0\cr
0&0&0&-1\cr
0&0&1&0\cr}\right).
$$
\medskip

Let  $\jj$ denote the set of all orthogonal complex structures in $\rr^4$
giving $\rr^4$ the same orientation as $J\st$. This means that for any pair
$x_1, x_2$ of $J$-independent vectors the basis $x_1, Jx_1, x_2, Jx_2$ gives
the same orientation of $\rr^4$ as $J\st$. Note that this orientation does
not depend on the particular choice of $x_1, x_2$. Orthogonality here means
just that $J$ is an orthogonal matrix. The following lemma summarizes the
elementary facts, which we need for the sequel.

\smallskip
\noindent
\bf Lemma 1.1.1. \it The elements of $\jj$ have the form

\smallskip
$$
 J_{s} =
\left( \matrix{0&-s &c_1&c_2\cr
s &0&c_2&-c_1\cr
-c_1&-c_2&0&-s \cr
-c_2&c_1&s &0\cr}\right),\eqno(1.1.1)
$$

\smallskip
\noindent
with $c_1^2 + c_2^2 + s^2 = 1$. One also has for $x\in \rr^4$

\smallskip
$$
 \omega (x, J_s x) = s \Vert x\Vert^2 . \eqno(1.1.2)
$$

\smallskip
\rm
The proof will be omitted. We remark that the set $\jj$ is a unit
two-dimensional sphere $S^2$ in $\rr^3$ with coordinates $c_1, c_2,s $. We
note also that the number $\omega (x, J_s x)$ does not depend on the
choice of a \sl unit \rm vector $x$. One also remarks that the standard
structure corresponds to the north pole of $S^2$ and structures tamed by
$\omega $ constitute the upper half-sphere.

Let $M_1$ and $M_2$ be two smooth oriented surfaces in the unit ball $B\subset
\rr^4$ with zero as a common point. Let $v_1, w_1$ and $v_2, w_2$ are oriented
bases of $T_0M_1$ and $T_0M_2$ respectively. Suppose that $M_1$ and $M_2$
intersect transversally at zero, {\sl i.e.} $v_1, w_1, v_2, w_2$ is the basis
of $\rr^4$. We say that their intersection is positively if this basis gives
the same orientation of $\rr ^4$ as standard one.

An immersed surface $M$ in an almost complex manifold $(X, J)$ is called
(non\-pa\-ra\-met\-rized)  $J$-holomorphic curve if $TM$ is $J$-invariant.

\smallskip
\noindent
\bf Lemma 1.1.2. \it Let $M$ be a $\omega $-positive compact surface immersed
into a K\"ahler surface $(X,\omega )$ with only double positive
self-intersections,
and let $U_1\subset \subset U$ be  neighbourhoods of $M$. Then there exists a
smooth curve $\{ J_t\} _{t\in [0, 1]} $ of almost complex structures on $X$
such that:

\smallskip
a) $J_0$ is the  given integrable complex structure on $X$;

b) for each $t\in [0, 1]$ the set $\{ x\in X: J_{t}(x) \not= J_0(x) \} $ is
contained in $U_1$;

c) $M$ is $J_1$-holomorphic;

d) all ${J_t}$ are tamed by the given form $\omega $, {\sl i.e.}
$\omega(v,J_tv) > 0$ for every nonzero $v\in TX$.

\medskip
\noindent
\bf Proof. \rm Let $N$ be a normal bundle to $M$ in $X$ and $V_1$ a
neighbourhood of the zero section in $N$. Shrinking $V_1$ and $U_1$ we can
assume
that $U_1$ is an image of $V_1$ under an $\exp$-map in $J_0$-Hermitian metric
$h$ associated to $\omega $. More precisely, one should take the  Riemannian
metric associated to $h$, {\sl i.e.} $g = {\ss Re} h$. Shrinking $V_1$ and
$U_1$
once more, if necessary, we can extend the distribution of tangent planes to
$M =$ (zero section of $N$) to the distribution $\{ L_x\} _{x\in V_1}$ of
$\omega $-positive planes on $V_1$. Here we do not distinguish between
$\omega$ and its lift onto $V_1$ by $\exp$. Denote by $N_x$ the subspace
in $T_xV_1$ which is $g$-orthogonal to $L_x$. We can choose the distribution $\{
 L_x\}$
in such a way that if $\exp(x)= \exp(y)$ for some $x\not= y$ from $V_1$ then
$d\exp_x(L_x)=d\exp_y(N_y)$ and $d\exp_y(L_y)=d\exp_x(N_x)$. In particular,
we suppose that $M$ intersects itself $g$-orthogonally. For every $x\in U_1$
choose an orthonormal basis $e_1(x), e_2(x)$ of $L_x$ such that $\omega_x
(e_1(x), e_2(x)) > 0$, and orthonormal basis $e_3(x), e_4(x)$ of $N_x$ such
that the basis $e_1(x), e_2(x), e_3(x), e_4(x)$ gives the same orientation of
$T_xV_1$ as $\omega^2$.

Define an almost complex structure $J$ on $U_1$ by $Je_1 = e_2, Je_3 = e_4$.
$J$ depends smoothly on $x$, even when $e_j(x)$ are not smooth
in $x$. Further $\omega _x(e_1(x), J_xe_1(x)) > 0$ and thus from
{\sl Lemma 1.1.1} we have that $\omega _x(e_3(x), J_xe_3(x)) > 0$. This means
that our $J$ is tamed by $\omega $. Note also that $M$ is $J$-holomorphic.

Denote by $\jj_x$ the sphere of $g$-orthogonal complex structures on $T_x
X$ as
in {\sl Lemma 1.1.1}. Let $\gamma_x$ be the shortest geodesic on $\jj_x$
joining $J_0(x)$, the~given integrable structure, with $ J_x$. Put $J_t(x)
= J_{\gamma _x(t\cdot \Vert \gamma _x\Vert \cdot \phi (x ))}$. Here $\Vert
\gamma_x \Vert $ denotes the length of $\gamma _x$, $\phi$ is smooth on $X$
with support in $U_1$ and identically one in the neighbourhood of $M$. The
curve $\{ J_t \}$ satisfies all the~conditions of our lemma.\qed

\bigskip\noindent
\sl 1.2. The Genus formula for  immersed symplectic surfaces.

\smallskip
\rm
Let us prove now the Genus formula for \sl immersed symplectic \rm surfaces.
Let $u:\bigsqcup_{j=1}^dS_j\to (X,\omega)$ be a reduced compact symplectic
surface immersed into a symplectic four-dimensional manifold. Let $g_j$ denote
the genus of $S_j$ and $M_j=u(S_j)$. Denote by $[M]^2$ the homological
self-intersection number of $M$. Define a geometrical self-intersection number
$\delta $ of $M$ in the following way. Perturb $M$ to obtain (also symplectic)
surface $\widetilde M$ with only transversal double points. Then $\delta$
will be
the sum of indices of intersection over those double points. Those indices
can be equal to $1$ or $-1$.

Surface $X$ carries an almost complex structure $J$, which is tamed by
$\omega$, {\sl i.e.} $\omega(\xi, J\xi)>0$ for any nonzero $\xi \in TX$,
see [Wn]
or [Gro]. Denote by $c_1(X, J)$ the first Chern class of $X$ with respect to
$J$. Since, in fact, $c_1(X, J)$ doesn't depend on continuous changes of $J$
and since  the set of $\omega$-tamed almost complex structures is
contractible, see [Gro], we usually omit the dependence of $c_1(X)$ on $J$.

\smallskip\noindent
\bf Lemma 1.2.1. \it Let $M=\bigcup_{j=1}^d M_j$ be a compact immersed
symplectic surface with transversal self-intersections in four-dimensional
symplectic manifold $X$. Then

\smallskip
$$
\sum_{j=1}^d g_j = {[M]^2 - c_1(X)[M]\over2} + d - \delta .\eqno(1.2.1)
$$

\smallskip\noindent\bf Proof.
\rm Replacing every $M_j$ by its small perturbation, we can suppose that $M_j$
has only transversal double self-intersection points. Let $N_j$ be a normal
bundle to $M_j$ and let $\widetilde M_j$ denote the zero section
of $N_j$. Let also $\exp_j $ be the exponential map from a~neighbourhood $V_j$
of $\widetilde M_j\subset N_j$ onto the neighbourhood $W_j$ of $M_j$.
Lift $\omega$ and $J$ onto $V_j$. Since $\widetilde M_j$ is embedded to $V_j$,
we can apply the {\sl Lemma 1.1.2} to obtain the $\omega$-tame
almost complex structure $J_j$ on $V_j$ in which $\widetilde M_j$ is
pseudo-holomorphic.

For every $j$ we have now the following exact sequence of complex bundles:

\smallskip
$$
0\longrightarrow TS_j {\buildrel du\over \longrightarrow }
E_j \buildrel pr \over \longrightarrow
N_j \longrightarrow 0 \eqno(1.2.2)
$$
\smallskip
Here $E_j=(u^*TX)\ogran_{S_j}$ is endowed with complex structure given by $J$.
Since $du$ is nowhere degenerate {\sl complex} linear morphism, $N_j \deff
E_j/du(TM_j)$ be a complex line bundle over $S_j$. From (1.2.2) we get

\smallskip
$$
c_1(E_j) = c_1(TS_j) + c_1(N_j). \eqno(1.2.3)
$$
\smallskip
Observe now that $c_1(E_j) = c_1(X)[M_j]$ and that $c_1(TS) = \sum_{j=1}^d
c_1(TS_j) = \sum_{j=1}^d(2-2g_j) = 2d - 2\sum_{j=1}^d g_j$. Further,
$c_1(N_j)$ is the algebraic number of zeros of a generic smooth section
of $N_j$. To compare this
number with self-intersection of $M_j$ in $X$, note that if we move $M_j$
generically to obtain $M_j'$, then the~intersection number ${\ss int}(M_j,
M_j')$ is equal to the algebraic number of zeros of generic section
of $N_j$ plus two times the
sum of intersection numbers of $M_j$ in self-intersection points, {\sl i.e.}
$[M_j]^2 = c_1(N_j) + 2\delta_j$. So

\smallskip
$$
c_1(X)[M_j] = 2 - 2 g_j + [M_j]^2 - 2\delta _j . \eqno(1.2.4)
$$

\smallskip\noindent
Now everything that left is to take the sum over $j=1,\ldots, d$ and to
remark that intersection points of $M_i$ with $M_j$ for $i\not=j$ make double
contribution to $[M]^2$.\qed

\bigskip\bigskip\bigskip
\noindent
\bf 2. The $D_{u, J}$-operator and its properties.

\nobreak\smallskip
\noindent \sl
\rm
In this paragraph we want to introduce, following Gromov, $D_{u, J}$-operator
associated to the $J$-holomorphic curve $u : S \to X$ in almost complex
manifold $(X, J)$, where $J$ is supposed to be of class $C^1$.

\medskip\noindent
\sl 2.1. The definition of $D_{u, J}$-operator.


\nobreak\smallskip\rm
Recall that a $C^1$ map $u:S\to X$ from a Riemannian surface $S$ with a complex
structure $J_S$ is called a parameterized pseudo-holomorphic curve with respect
to $J_S$ and $J$ if it satisfies the equation

\smallskip
$$
du + J\scirc du\scirc J_S = 0 \eqno(2.1.1)
$$

This simply means that $du\scirc J_S = J\scirc du$, {\sl i.e.} $du:T_pS \to
T_{u(p)}X$ is {\sl complex} linear map for every $p\in S$. The equation (2.1.1)
is an elliptic quasi-linear PDE of the order one. We are interested in behavior
of solutions of (2.1.1), in particular, when the coefficients $J$ and $J_S$
change. So we need to choose appropriate functional spaces both for solutions
and for coefficient of (2.1.1). The authors' choice is based on the following
facts:

{\sl a)} the minimal reasonable smoothness of an almost complex structure $J$
on $X$, for which the Gromov operator $D_{u, J}$ can be defined, is $C^1$,
see explicit formula (2.1.6);

{\sl b)} it is more convenient to operate with Banach spaces and manifolds and
thus with finite smoothnesses like $C^k$, $C^{k,\alpha}$ or $L^{k,p}$, than
with $C^\infty$-smoothness defining only Frechet-type topology, see
{\sl Lemmas 3.2.4} and {\sl 3.3.1};

{\sl c)} the equation (2.1.1) is defined also for $u$ lying in Sobolev-type
spaces $L^{k,p}(S, X)$ with $k\ge1$, $1\le p\le\infty$, and $kp>2$; such
solutions are $C^1$-smooth and the topology on the space of solution is,
in fact, independent of the particular choice of such Sobolev space, see
{\sl Corollary 3.2.2}.

{\sl d)} for $J\in C^k$ the coefficients of the Gromov operator $D_{u, J}$
are, in general, only $C^{k-1}$-continuous, see (2.1.6), and hence solutions
of the ``tangential equation'' $D_{u, J}v=0$ are only $L^{k,p}$-smooth, $1\le
p<\infty$; thus that for obtaining a {\sl smooth} structure on a space
of (parameterized) pseudo-holomorphic curves one should use Sobolev spaces
$L^{k', p}(S, X)$ with $k'\le k$.

\medskip
Fix a compact Riemannian surface $S$, {\sl i.e.}\ a compact, connected,
oriented, smooth manifold of real dimensions 2. Recall, that the Sobolev space
$L^{k, p}(S, X)$, $kp>2$, consists of those continuous maps $u: S \to X$,
which are represented by $L^{k, p}$-functions in local coordinates on $X$ and
$S$. This is a Banach manifold, and the tangent space $T_uL^{k, p}(S, X)$
to $L^{k, p}(S, X)$ in $u$ is the space $L^{k,p }(S, u^*(TX))$ of all
$L^{k,p}$-sections of the pull-back under $u$ of the tangent bundle $TX$.
Besides, one has the Sobolev imbeddings
$$
\matrix{
L^{k, p}(S, X)& \hookrightarrow & L^{k-1, q}(S, X),&
\quad & \hbox{for $1\le p<2$}   &  \hbox{and} & 1\le q\le {2p\over2-p},
\cr
L^{k, p}(S, X)& \hookrightarrow & C^{k-1, \alpha}(S, X),&
\quad & \hbox{for $2<p\le\infty$}& \hbox{and} & 0\le\alpha\le 1-{2\over p}.
\cr}
\eqno(2.1.2)
$$


\smallskip
Let $[\gamma ]$ be some homology class in  $\ssh_2(X, \zz)$. Fix $p$ with
$2<p<\infty$ and consider the~Banach manifold
$$
{\cal S} = \{u\in L^{1,p}(S, X):u(S)\in [\gamma]\}
$$
of all $L^{1,p }$-smooth mappings from $S$ to $X$, representing the~class
$[\gamma]$.

Denote by ${\cal J}$ the Banach manifold of $C^1$-smooth almost
complex structures on $X$. In other words, ${\cal J} = \{ J\in C^1
(X, {\ss End}(TX)):J^2=-\id  \} $. The tangent space to ${\cal J}$ at $J$
consists of $C^1$-smooth $J$-antilinear endomorphisms of $TX$,
$$
T_J{\cal J} = \{ I\in C^1(X, {\ss End}(TX)) : JI + IJ = 0\}
\equiv C^1(X, \Lambda^{0,1}X \otimes TX),
$$
where $\Lambda^{0,1}X$ denote the complex bundle of $(0,1)$-form on $X$.

Denote by ${\cal J}_S$ the Banach manifold of $C^1$-smooth complex structures
on $S$. Thus ${\cal J}_S = \{ J_S\in C^1(S, {\ss End}(TS)):J_S^2=-\id  \}$
and the tangent space to ${\cal J}_S$ at $J_S$ is
$$
T_{J_S}{\cal J} = \{ I\in C^1(S, {\ss End}(TS)) : J_S I + IJ_S = 0\}
\equiv C^1(S, \Lambda^{0,1}S\otimes TS).
$$

\smallskip
Consider also the subset ${\cal P}\subset {\cal S}\times {\cal J}_S \times
{\cal J}$ consisting of all triples $(u, J_S, J)$ with $u$ being
$(J_S,J)$-holomorphic, {\sl i.e.}

\smallskip
$$
{\cal P} = \{(u, J_S, J)\in {\cal S}\times {\cal J}_S\times {\cal J}: du +
J\scirc du\scirc  J_S = 0 \} \eqno(2.1.3)
$$

\state Lemma 2.1.1.
\it Let $J$ and $J_S$ be continuous almost complex structures on $X$ and $S$
respectively, and let $u\in L^{1,p}(S, X)$. Then
$$
du + J\scirc du\scirc J_S \in L^p(S,\Lambda^{0, 1}S \otimes u^*(TX)).
$$

\state Proof. \rm One can easily see that $du \in L^p(S,\hom_\rr(TS, u^*(TX))$.
On the other hand,
$$
(du + J\scirc du\scirc J_S ) \scirc J_S= -J \scirc
(du + J\scirc du\scirc J_S ).
$$
which means that $du + J\scirc du\scirc J_S$ is $L^p$-integrable
$u^*(TX)$-valued (0,1)-form.
\qed

\smallskip
Consider a Banach bundle ${\cal T}\to {\cal S}\times {\cal J}_S\times {\cal J}$
with a fiber
$$
{\cal T}_{(u, J_S, J)} = L^p(S,\Lambda^{0, 1}S\otimes
u^*(TX)),
$$
where $TS$ and $TX$ are equipped with complex structures $J_S$ and $J$
respectively. ${\cal T}$ has two distinguished sections:

1) $\sigma_0 \equiv 0$, the~zero section of ${\cal T}$;

2) $\sigma_{\dbar}(u, J_S, J) = du + J\scirc du\scirc J_S$;

\noindent
and by definition ${\cal P}$ is the zero-set of $\sigma_{\dbar}$.

\medskip
Let us compute the tangent space to ${\cal P}$ at the point $(u, J_S, J)$. Let
$(u_t, J_S(t),J(t))$ be a curve in ${\cal P}$ such that $(u_0, J_S(0), J(0)) =
(u, J_S, J)$. Let
$$
(v,\dot J_S,\dot J) \deff
\left({du_t\over dt}|_{t=0},
{dJ_S(t)\over dt}|_{t=0},
{dJ(t)\over dt}|_{t=0}\right)
$$
be the tangent vector to this curve and hence to ${\cal P}$ at $t=0$.
The condition $(u_t, J_S(t), J(t)) \in {\cal P}$ means that
$$
du_t + J(u_t, t)\scirc du_t\scirc J_S(t) = 0 \eqno(2.1.4)
$$
in $ L^p(S,\Lambda^{0, 1}(S)\otimes u_t^*(TX))$. Let $\nabla $ be
some symmetric connection on $TX$, {\sl i.e.} $\nabla_YZ - \nabla_ ZY = [Y,
Z]$. The covariant differentiation of (2.1.4) with respect to $t$ gives
$$
\nabla_{\partial \over \partial t}(du_t) + (\nabla_vJ)(du_t\scirc J_S)
+ J(u_t, t)\scirc \nabla_{\partial \over \partial t}(du_t)\scirc J_S +
J\scirc du_t\scirc \dot J_S +
\dot J\scirc du_t\scirc J_S = 0
$$
Let us show that $\nabla_{\partial \over \partial t}(du_t) = \nabla v$.
Indeed, for $\xi \in TS$ one has
$$
(\nabla_{\partial \over \partial t}
du_t)(\xi ) = \nabla_{\partial \over \partial t}[du_t(\xi )] =
\nabla_{\partial \over \partial t}({\partial u_t\over \partial \xi }) =
\nabla_{\xi }({du_t\over dt}) = \nabla_{\xi }v.
$$
So every vector $(v,\dot J_S,\dot J)$ which is tangent to $\cal P$
satisfies the equation

\smallskip
$$
\nabla v + J\scirc \nabla v\scirc J_S + (\nabla_vJ)\scirc (du\scirc
J_S) + J\scirc du\scirc \dot J_S +
\dot J\scirc du\scirc J_S = 0. \eqno(2.1.5)
$$


\bigskip
\noindent
{\bf Definition 2.1.1.} \rm
Let $u$ be a $J$-holomorphic curve in $X$. Define the operator $D_{u, J}$
on $L^{1,p}$-sections $v$ of
$u^*(TX)$ as
$$
D_{u, J}(v) = \msmall{1\over2}\bigl(\nabla v + J\scirc\nabla v\scirc J_S
+ (\nabla_vJ) \scirc (du\scirc J_S) \bigr)
\eqno(2.1.6)
$$

\bigskip
\noindent
{\bf Remark.} This operator plays crucial role in studying  properties of
pseudo-holomorphic curves. In its definition we use the symmetric connections
instead of those compatible with $J$, as is done in [Gro]. The matter is that
one can use the same connection $\nabla$ for changing almost complex
structures $J$. The lemmas below justify our choice.

\bigskip
\noindent
{\bf Lemma 2.1.2. }\it $D_{u, J}$ doesn't depend on the choice of a~symmetric
connection $\nabla $ and is an $\rr$-linear operator from $L^{1,p}(S,
u^*(TX))$ to $L^p(S, \Lambda^{0, 1}S \otimes u^*(TX))$.

\smallskip
\rm
\noindent
{\bf Proof.} Let $\widetilde \nabla $ be another symmetric connection on $TX$.
Consider the bilinear tensor on $TX$, given by formula $Q(Z, Y) {:=} \nabla_ZY
- \widetilde \nabla_ZY$. It is easy to see, that $Q$ is symmetric on $Z$ and
$Y$.
Note further that $\nabla_{\xi }v - \widetilde\nabla_{\xi }v =
\nabla_{du(\xi )}v
- \widetilde \nabla_{du(\xi )}v = Q(du(\xi), v)$ and besides of that
$(\nabla_ZJ)(Y) - (\widetilde \nabla_ZJ)(Y) = Q(Z, JY) - JQ(Z, Y)$. From here
one
gets that
$$
2(D_{u, J}v)(\xi ) - 2(\widetilde D_{u, J}v)(\xi) =
\nabla_{\xi }v - \widetilde \nabla_
{\xi }v + J(\nabla_{J_S\xi }v - \widetilde \nabla_{J_S\xi }v) + (\nabla_vJ -
\widetilde \nabla_vJ)du(J_S\xi ) =
$$
$$
= Q(du(\xi ), v) + JQ(du(J_S\xi ), v) + Q(v, JduJ_S
\xi ) - JQ(v, du(J_S\xi)) =
$$
$$
=Q(du(\xi ), v) + Q(JduJ_S\xi , v) = Q((du + JduJ_S)(\xi ), v) = 0.
$$

Now let us show that $D_{u, J}(v)$ is $J_S$-antilinear:
$$
2D_{u, J}(v)[J_S\xi] = \nabla_{J_S\xi }v + J\scirc \nabla_{J^2
_S\xi }v + (\nabla_vJ)\scirc (du\scirc J^2_S)(\xi ) =
$$
$$
= \nabla_{J_S\xi }v -
J(\nabla_{\xi }v) - (\nabla_vJ)(du(\xi )) = -J[\nabla_{\xi }v + J(\nabla
_{J_S\xi }v) - J\nabla_vJdu(\xi )] =
$$
$$
= -J[\nabla_{\xi }v + J\scirc \nabla_{J_S\xi}v + (\nabla_vJ)
(du\scirc J_S(\xi )] = -2J\, D_{u, J}(v)[\xi].
$$
Here we use the~fact that $J\scirc \nabla_vJ + \nabla_vJ\scirc J = 0$ and that
$du\scirc J_S = J\scirc du $.\qed

\smallskip
This lemma allows us to obtain expression (2.1.5) also by computation in local
coordinates $x_1, x_2$ on $S$ and $u_1,\ldots, u_{2n}$ on $X$.

\bigskip\noindent
\sl 2.2. Operator $\dbar_{u, J}$ and the~holomorphic structure on
the~induced bundle.\rm

\nobreak
Now we need to understand the structure of operator $D_{u, J}$ in more
details. The problem arising here is that $D\deff D_{u,J}$ is only
$\rr$-linear. So we decompose it into the $J$-linear and $J$-antilinear parts.
For a $J$-holomorphic curve $u:S\to X$ denote $u^*(TX)$ by $E$. Note that the
bundle $E$ carries a complex structure, namely $J$ itself or more accurately
$u^*J$. For $\xi \in C^1 (S, TS)$ and $v\in L^{1,p}(S, E)$ write $D_{\xi }v =
{1\over2}[ D_{\xi }v - JD_{\xi }(Jv)] + {1\over2} [D_{\xi }v + JD_{\xi
}(Jv)] = \dbar_{u, J}[v](\xi) + R(v,\xi )$.

\medskip
\noindent
\bf Definition 2.2.1. \rm The operator $\dbar_{u, J}$, introduced
above as $J$-linear part of $D_{u, J}$, we shall call the \sl $\dbar$-operator
for a~$J$-holomorphic curve $u$. \rm

\medskip
\noindent
\bf Lemma 2.2.1. \it $R$ is continuous $J$-antilinear operator from
$E$ to $\Lambda^{0,1}\otimes E$ of the order zero, satisfying
$R\scirc du \equiv 0 $.

\noindent \rm
{\bf Proof.} $J$ antilinearity of $R$ is given by its definition. Compute
$R(v,\xi )$ for $v\in L^{1,p}(S, E)$ and $\xi \in C^1(S, TS)$, setting
$w\deff du(\xi)$ and $D\deff D_{u, J}$ to simplify the notations:

\smallskip
\noindent
$$
4R(v,\xi ) = 2D[v](\xi ) + 2JD[Jv](\xi ) =
$$

\smallskip
$$
=\nabla_{\xi }v + J\nabla_{J_S\xi }v + \nabla_vJ\scirc du\scirc J_S(\xi )
$$
$$
+J\nabla_{\xi }(Jv) + J^2\nabla_{J_S(\xi )}(Jv) +
J \nabla_{Jv}J\scirc du\scirc J_S(\xi ) =
$$

\smallskip
$$
= \nabla_\xi v + J\nabla_{J_S\xi }v + \nabla_vJ \scirc du\scirc J_S(\xi ) +
$$
$$
+ J^2\nabla_{\xi }v + J(\nabla_wJ)v + J^3\nabla_{J_S\xi }v +
J^2(\nabla_{Jw}J)v + J\scirc \nabla_{Jv}J\scirc Jw.
$$

\smallskip
Here we used the relations $\nabla_{du(\xi )}J = \nabla_wJ$, $du(J_S\xi) =Jw$
and $\nabla_{du(J_S\xi )}J = \nabla_{Jw}J$. Contracting terms
we obtain that
$$
4R(v,\xi ) = \nabla_vJ(Jw) + J(\nabla_wJ)v - (\nabla_{Jw}J)v +
J(\nabla_{Jv}J(Jw)) =
$$
$$
= (\nabla_vJ\scirc J)w - (\nabla_wJ\scirc J)v - (\nabla_{Jw}J)v +
(\nabla_{Jv}J)w = 4N(v, w),
$$
where $N(v, w)$ denotes the torsion tensor of the almost complex structure
$J$, see [Li], p.183, or [Ko-No], vol.II., p.123, where another normalization
constant for the almost complex torsion is used. Finally we obtain
$$
R(v,\xi ) = N(v, du(\xi ))\eqno(2.2.1)
$$

\smallskip
$N$ is antisymmetric and $J$-antilinear on both arguments, so
$$
-JR\bigl(du(\eta ),\xi \bigr) = R\bigl(du(\eta ), J_S\xi \bigr) =
N\bigl(du(\eta ), du(J_S\xi ) \bigr)
= N\bigl(du(\eta ), du(\eta )\bigr) =0
$$
if $\xi$ and $\eta$ were chosen in such a way that $J_S(\xi )=\eta $.
The relation $R(du(\xi ),\xi )=0$  obviously follows from (2.2.1).
\qed


\bigskip
The $J$-linear operator $\dbar_{u, J} : L^{1,p}(S, E) \to
L^p(S,\Lambda^{0, 1}S\otimes E)$ defines in a standard way, see
Appendix for the~proof, a holomorphic structure on the bundle $E$.
This follows from the fact that due to one-dimensionality of
$S$ there is no integrability condition on the operator $\dbar$.
The sheaf of holomorphic sections of $E$ we shall denote by ${\cal
O}(E)$. The tangent bundle $TS$ to our Riemannian surface also carries a
natural holomorphic structure. We shall denote by ${\cal O}(TS)$
the~corresponding analytic sheaf.

\medskip
\noindent
\bf Lemma 2.2.2. \it Let $u : (S, J_S) \to (X, J)$ be a nonconstant
pseudo-holomorphic curve in almost complex manifold $X$. Then $du$ defines
an injective analytic morphism of analytic sheaves
$$
0\longrightarrow {\cal O}(TS)
\buildrel du \over \longrightarrow {\cal O}(E) \eqno(2.2.2)
$$

\smallskip
\noindent
where $E = u^*(TX)$ is equipped with holomorphic structure defined as above
by operator $\dbar_{u, J}$.

\smallskip
\noindent
\rm {\bf Proof.} Injectivity of a sheaf homomorphism is equivalent to its
nondegeneracy, which is our case.

To prove holomorphicity of $du$ it is sufficient to show that for $\xi,\eta \in
C^1(S, TS)$ one has
$$
\bigl(\dbar_{u, J}(du(\xi ))\bigr)(\eta ) =
du\bigl((\dbar_S\xi ) (\eta )\bigr) \eqno(2.2.3)
$$

\smallskip
\noindent
where $\dbar_S$ is the~usual $\dbar $-operator on $TS$. We shall use
the relation, which is in fact the~definition for $\dbar_S$:

\smallskip
$$
(\dbar_S\xi )(\eta ) =
\msmall{1\over2}\left( \nabla_{\eta }\xi + J_S\nabla_{J_S\eta}\xi \right)
\eqno(2.2.4)
$$
\noindent Here $\nabla$ is a symmetric connection on $S$, compatible with
$J_S$. One has
$$
2\cdot (\dbar_{u, J}du(\xi ))(\eta ) = \nabla_{\eta }(du(\xi) ) +
J\nabla_{J_S\eta }(du(\xi))  + (\nabla_{du(\xi )}J)(du(J_S\eta )) =
$$
$$
= (\nabla_{\eta }du)(\xi ) + du(\nabla_{\eta }\xi ) +
J(\nabla_{J_S\eta}du)(\xi ) + J(du(\nabla_{J_S\eta }\xi ))
+ (\nabla_{du(\xi )}J)(du(J_S\eta )) =
$$
$$
= du\bigl[ \nabla_{\eta }\xi + J_S \nabla_{J_S\eta }\xi \bigr] +
\bigl[ (\nabla_{\eta }du)(\xi ) + J(\nabla_{J_S\eta}du)(\xi )
+ (\nabla_{du(\xi )}J)(du(J_S\eta )) \bigr]. \eqno(2.2.5)
$$
The first term of $(2.2.5)$ is $2\cdot du(\dbar_S\xi)(\eta)$.
To cancel the second one we use the identities $(\nabla_\xi du)[\eta] =
(\nabla_\eta du)[\xi]$, $\nabla_wJ\circ J = - J\circ \nabla_w J$, and
$(\nabla_\xi du)\circ J_S= J\circ(\nabla_\xi du) + \nabla_{du(\xi)}J\circ du$.
The last identity is obtained via covariant differentiation of $du\circ J_S=
J\circ du$. Consequently we obtain
$$
(\nabla_{\eta }du)(\xi ) + J(\nabla_{J_S\eta}du)(\xi ) +
 (\nabla_{du(\xi )}J)(du(J_S\eta ))=
$$
$$
(\nabla_{\xi }du)(\eta ) + J(\nabla_\xi du)(J_S\eta ) +
 (\nabla_{du(\xi )}J)(du(J_S\eta ))=
$$
$$
(\nabla_{\xi }du)(\eta ) +
J^2(\nabla_\xi du)(\eta ) +
J (\nabla_{du(\xi)}J)(du(\eta))
+ (\nabla_{du(\xi )}J)(du(J_S\eta) )=
$$
$$
=(J \circ\nabla_{du\xi}J)(du\eta)
+ (\nabla_{du(\xi )}J\circ J)(du(\eta ))=0.
$$
\qed

\medskip\noindent
\bf Remark. \rm We can give an alternative proof of both {\sl Lemmas 2.2.1} and
{\sl 2.2.2}, which does not use direct calculation. Fix a~complex structure
$J_S$ on $S$ and let $\phi_t$ be the~one
parameter group of diffeomorphisms of $S$, generated by a~vector field
$\xi$. Then ${d\over dt}|_{t=0}(\dbar\phi_t)=D_{J_S, \id }\xi=\dbar \xi$.
Let a~$J$-holomorphic map $u:S\to X$ be also fixed. Then $
{d\over dt}|_{t=0}(u\circ \phi_t)=du(\xi)$ and consequently
$$
D_{J, u}(du(\xi))={d\over dt}\bigm|_{t=0}\dbar_J(u\circ \phi_t)=
{d\over dt}\bigm|_{t=0}(du \scirc \dbar \phi_t)=
du\bigl({d\over dt}\bigm|_{t=0}(\dbar_{J_S} \phi_t)=
du(\dbar\xi)
$$
or equivalently $D_{J, u}\circ du=du \circ \dbar$. Taking $J$-antilinear
part of last equality we obtain $R\circ du=0$. Nevertheless we shall use the
explicit form of the antilinear part of the Gromov operator given in (2.2.1).

\bigskip
The zeros of analytic morphism $du : {\cal O}(TS) \to {\cal O}
(E)$ are isolated. So we have the following

\medskip
\noindent
\bf Corollary 2.2.3. \rm([Sk]). \it The set of critical points of
pseudo-holomorphic curve in almost complex manifold $(X, J)$ is discrete,
provided $J$ is of class $C^1$.

\smallskip\rm
For $C^\infty$-structures this result is due to McDuff, see [McD-3].

\smallskip
The {\sl Lemma 2.2.2} makes possible to give the following

\smallskip
\noindent\bf
Definition 2.2.2. \sl By the order of zero ${\ss ord}_p du$ of the differential
$du$ at a point $p\in S$ we shall understand the order of vanishing at $p$ of
the holomorphic morphism $du : {\cal O}(TS)\to {\cal O}(E)$.
\rm

\smallskip
From (2.2.2) we obtain the following short exact sequence:

\smallskip
$$
0\longrightarrow {\cal O}(TS) \buildrel du \over\longrightarrow {\cal O}(E)
\longrightarrow {\cal N}\longrightarrow 0. \eqno(2.2.6)
$$

\smallskip
Here ${\cal N}$ is a quotient-sheaf ${\cal O}(E)/du(TS)$. We can decompose
${\cal N} = {\cal O}(N_0)\oplus {\cal N}_1$, where $N_0$ is a holomorphic
vector bundle and ${\cal N}_1 =
\bigoplus_{i=1}^P \cc_{a_i}^{n_i}$. Here $\cc_{a_i}^{n_i}$ denotes the sheaf,
supported at the~critical points $a_i\in S$ of $du$ and having a stalk
$\cc^{n_i}$ with $n_i = {\ss ord}_{a_i}du$,
the~order of zero of $du$ at $a_i$.

Denote by $[A]$ the divisor $\sum_{i=1}^Pn_i[a_i]$, and by ${\cal O}([A])$ a
sheaf of meromorphic functions on $S$ having poles in $a_i$ of order at most
$n_i$. Then (2.2.6) give rise to the exact sequence

\smallskip
$$
0\longrightarrow {\cal O}(TS)\otimes {\cal O}([A])\buildrel{du}\over
{\longrightarrow} {\cal O}(E)
\longrightarrow {\cal O}(N_0)\longrightarrow 0. \eqno(2.2.7)
$$

\smallskip
Denote by $L^p_{(0, 1)}(S, E)$ the spaces of $L^p$-integrable (0, 1)-forms
with coefficients in $E$. Then (2.2.3) together with {\sl Lemma 2.2.1}
implies that the following diagram is commutative

\smallskip
$$

\def\mapright#1{\smash{\mathop{\longrightarrow}\limits^{#1}}}
\def\mapdown#1{\Big\downarrow\rlap{$\vcenter{\hbox{$\scriptstyle#1$}}$}}
\matrix{0&\mapright{}&L^{1,p}(S, TS\otimes [A])&\mapright{du}&
L^{1,p}(S, E)&\mapright{pr}&L^{1,p}(S, N_0)&\mapright{}&0\cr
 & &\mapdown{\dbar_S}& &\mapdown{D_{u, J}}& &\mapdown{}& & \cr
0&\mapright{}&L^p_{(0,1)}(S, TS\otimes [A])&\mapright{du}&
L^p_{(0,1)}(S, E)&\mapright{}&
L^p_{(0,1)}(S, N_0)&\mapright{}&0\cr}
\eqno(2.2.8)
$$

\smallskip
This defines an operator $D_{u, J}^N : L^{1,p}(S, N_0) \longrightarrow
L^p_{(0,1)}(S, N_0)$ which is of the form $D_{u, J}^N =
\dbar_N + R$. Here $\dbar_N$ is a usual $\dbar $-operator on $N_0$
and $R\in C^0(S, \hom_\rr(N_0,\Lambda^{0, 1} \otimes N_0))$. This
follows from the fact that $D_{u, J}$ has the same form, see
{\sl Definition 2.2.1}.

\medskip
\state Definition 2.2.3. Let $E$ be a holomorphic vector bundle over
a~compact Riemannian surface $S$ and let $D:L^{1,p}(S, E)\to L^p(S,
\Lambda^{0,1}S \otimes E)$ be an operator of the~form $D=\dbar + R$ where
$R\in L^p\bigl(S,\,\hom_\rr(E,\,\Lambda^{0, 1}S\otimes E) \bigr)$
with $2<p<\infty$. Define $\ssh^0_D(S, E)\deff \ker D$ and $\ssh^1_D(S, E)
\deff \coker D$.

\smallskip
\state Remark. It is shown in {\sl Lemma 2.3.2} below that for given $S$, $E$
and $R\in L^p$, $2<p<\infty$, one can define $\ssh^i_D(S, E)$
as (co)kernel of the operator $\dbar +R: L^{1,q}(S, E) \to
L^q(S,\,\Lambda^{0, 1}S\otimes E)$ for any $1<q\le p$. So the definition
is independent of the choice of functional space.

By the standard lemma of homological algebra we get from (2.2.8) the following
long exact sequence of $D$-cohomologies.

\smallskip
$$

\def\mapright#1{\smash{\mathop{\longrightarrow}\limits^{#1}}}
\def\mapdown{\Big\downarrow}
\matrix{
0& \mapright{}& \ssh^0(S, TS\otimes [A]) &\mapright{} & \ssh^0_D(S, E)
 & \mapright{}& \ssh^0_D(S, N_0)   &\mapright{\delta} &\cr
\vphantom{\mapdown}&&&&&&&&\cr
 & \mapright{}& \ssh^1(S, TS\otimes [A]) &\mapright{} & \ssh^1_D(S, E)
 & \mapright{}& \ssh^1_D(S, N_0)         &\mapright{} &0.
}
\eqno(2.2.9)
$$

\bigskip\smallskip
\noindent \sl 2.3. Surjectivity of $D_{u, J}$.

\nobreak\smallskip\rm
For the {\sl Step 2} of the~{\sl Theorem 1} we shall need a~result of
Gromov ([Gro]) and Hofer-Lizan-Sikorav ([Hf-L-Sk]) about
surjectivity of  $D^N_{u, J}$, namely a vanishing theorem
for $D$-cohomologies. First we prove some technical statements.

\medskip
\state Lemma 2.3.1. {\it Let $X$ and $Y$ be Banach spaces and $T:X\to Y$
a~closed dense defined unbounded operator with the~graph $\Gamma=\Gamma_T$
endowed with the~graph norm $\Vert x \Vert_\Gamma= \Vert x \Vert_X + \Vert Tx
\Vert_Y$. Suppose that the~natural map $\Gamma\to X$ is compact. Then

 \sli $\ker(T)$ is finite-dimensional;

 \slii $\im(T)$ is closed;

 \sliii the~dual space $\bigl(Y/\im(T)\bigr)^*$ is naturally isomorphic
to $\ker(T^*:Y^*\to X^*)$.
}

\state Proof. Obviously, for $x\in \ker(T)$ one has $\Vert x \Vert_X =
\Vert x \Vert_\Gamma$. Let $\{x_n\}$ be a sequence in $\ker(T)$ which is
bounded in $\Vert \cdot \Vert_X$ norm. Then it is bounded in $\Vert \cdot
\Vert_\Gamma$ norm and hence relatively compact with respect to
$\Vert \cdot \Vert_X$ norm. Thus the~unit ball in $\ker(T)$ is compact
which implies the~statement \sli of the~lemma.

 Due to finite-dimensionality, there exists a~closed complement $X_0$ to
$\ker(T)$ in $X$. Let $\{x_n\}$ be a~sequence in $X$ such that $Tx_n \lrar y
\in Y$. Without losing generality we may assume that $x_n$ belong to $X_0$.
Suppose that $\Vert x_n \Vert_X \lrar \infty$. Denote $\tilde x_n \deff {x_n
\over \Vert x_n \Vert_X }$. Then $\Vert \tilde x_n \Vert_\Gamma$ is bounded
and hence some subsequence of $\{\tilde x_n\}$, still denoted by $\{\tilde
x_n\}$, converges in $X_0$ to some $\tilde x$. Note that $\tilde x\not=0$
because $\Vert \tilde x \Vert_X = \lim \Vert \tilde x_n \Vert_ X =1$. On
the~other hand, one can see that $T\tilde x_n \lrar 0\in Y$. Since $\Gamma$ is
closed, $(\tilde x, 0)\in \Gamma$ and hence $\tilde x\in \ker(T) \cap X_0
=\{0\}$. The~contradiction shows that the~sequence $\{x_n\}$ must be bounded
in $X$. Since $\{Tx_n \}$ is also bounded in $Y$, some subsequence of
$\{x_n\}$, still denoted by $\{x_n\}$, converges in $X_0$ to some $x$. Due to
the~closeness of $\Gamma$, $Tx=y$. Thus $\im (T)$ is closed in $Y$.

 Denote $Z\deff\ker(T^*:Y^*\to X^*)$ and let $h\in \bigl(Y/\im(T)
\bigr)^*$. Then $h$ defines a~linear functional on $Y$,\ie some element
$h'\in Y^*$, which is identically zero on $\im(T)$. Thus for any $x$ from
the~domain of definition of $T$ one has $\langle h', Tx \rangle =0$, which
implies $T^*(h')=0$. Consequently, $h'$ belongs to $Z$. Conversely, every
$h'\in Z$ is a~linear functional on $Y$ with $h'(Tx)= \langle T^*h', x\rangle
=0$ for every $x$ from the~domain of definition of $T$. Thus $h'$ is
identically zero on $\im(T)$ and is defined by some unique $h\in
\bigl(Y/\im(T)\bigr)^*$.\qed

\bigskip
\state Lemma 2.3.2. {\sl (Serre duality for $D$-cohomologies.)}
{\it Let $E$ be a holomorphic vector bundle over a~compact Riemannian
surface $S$
and let $D:L^{1,p}(S, E)\to L^p(S, \Lambda^{0, 1}S\otimes E)$ be
an operator of the~form $D=\dbar + R$ where $R\in L^p\bigl(S,\,
\hom_\rr(E,\,\Lambda^{0, 1}S\otimes E) \bigr)$ with $2<p <\infty$. Let
also $K\deff \Lambda^{1, 0}S$ be the~canonical holomorphic line
bundle of $S$. Then there exists the~naturally defined operator
$$
D^*=\dbar- R^* : L^{1,p}(S, E^* \otimes K)
\to L^p(S,\Lambda^{0, 1} \otimes E^* \otimes K)
$$
with $R^* \in L^p\bigl(S,\,\homr(E^*\otimes K,\,
\Lambda^{0, 1}S \otimes E^*\otimes K) \bigr)$ and the~natural isomorphisms
$$
\ssh^0_D(S,\, E)^*\cong \ssh^1_{D^*}(S,\, E^*\otimes K),
$$
$$
\ssh^1_D(S,\, E)^*\cong \ssh^0_{D^*}(S,\, E^*\otimes K).
$$
If, in addition, $R$ is $\cc$-antilinear, then $R^*$ is also
$\cc$-antilinear.
}

\medskip
\state Proof. For any $1<q\le p$ we associate with $D$ an~unbounded
dense defined operator $T_q$
from $X_q\deff L^q(S, E)$ into $Y_q\deff L^q(S,\Lambda^{0, 1}S\otimes E)$
with the~domain of definition $L^{1, q}(S, E)$.
The elliptic regularity of $D$ (see {\sl Lemma 3.2.1} below ) implies that
$$
\Vert \xi \Vert_{L^{1,q}(S, E)} \le C(q)
\left(
\Vert \dbar\xi +R\xi\Vert_{L^q(S, E)} +
\Vert \xi \Vert_{L^q(S, E)} \right).
$$
Consequently, $T_q$ are closed and satisfy the~hypothesis of
{\sl Lemma 2.3.1}. For $q>q_1$ we have also the~natural imbedding $X_q
\hookrightarrow X_{q_1}$ and $Y_q \hookrightarrow Y_{q_1}$ which commutes with
the~operator $D$. Moreover, due to regularity  of $D$
this imbedding maps $\ker T_q$ {\sl identically} onto $\ker T_{q_1}$.
Thus we can identify $\ssh^0_D(S, E)$ with any $\ker T_q$.

 Now note that for $q'\deff q/(q-1)$ we have the~natural isomorphisms
$$
\eqalign{
X_q^*\equiv &(L^q(S, E))^*
\cong L^{q'}(S,\Lambda^{0, 1}S \otimes E^* \otimes K),\cr
Y_q^*\equiv &(L^q(S,\Lambda^{0, 1}S\otimes E))^*
\cong L^{q'}(S,  E^* \otimes K),\cr
}
$$
induced by pairing of $E$ with $E^*$ and by integration over $S$. One can
easily check that the~dual operator $T_q^*$ is induced by the~differential
operator $-D^*: L^{1,q'}(S, E^* \otimes K)
\to L^{q'}(S,\Lambda^{0, 1} \otimes E^* \otimes K)$
of the~form $D^* =\dbar - R^*$. Really, for $\xi \in L^{1, q}(S, E)$
and $\eta\in L^{1, q'}(S, E^*\otimes K)$ one has
$$
\langle T_q\xi,\,\eta \rangle=
\int_S \langle \dbar\xi + R\xi,\,\eta \rangle =
\int_S \dbar\langle \xi ,\,\eta \rangle +
\int_S \langle \xi,\, -(\dbar - R^*)\eta \rangle =
\int_S \langle \xi ,\, -D^*\eta \rangle,
$$
since the~integral of any $\dbar$-exact $(1, 1)$-form vanishes. From {\sl
Lemma 2.3.1} we obtain the~natural isomorphisms $\ssh^1_D(S, E)^*\equiv (\coker
T_p)^*\cong \ker T_p^*$ and $\ssh^0_D(S, E)^*\equiv (\ker
T_p)^*\cong \coker T_p^*$,  which yields the~statement of the lemma. \qed

\medskip
\state Corollary 2.3.3. ([Gro], [H-L-Sk].) {\sl (Vanishing theorem for
$D$-cohomologies.)} {\it Let $S$ be a~Riemannian surface $S$ of the~genus $g$.
Let also $L$ be a~holomorphic {\sl line} bundle over $S$, equipped with
a~differential operator $D=\dbar + R$ with $R\in L^p\bigr(S,
\homr(L,\,  \Lambda^{0, 1}S\otimes L) \bigl)$, $p>2$. If $c_1(L)<0$,
then $\ssh^0_D (S,\, L)=0$. If $c_1(L)>2g-2$, then $\ssh^1_D (S,\, L)=0$.
}

\medskip
\state Proof. Suppose $\xi$ is a~nontrivial $L^{1,p}$-section of $L$ satisfying
$D\xi=0$. Then due to {\sl Lemma 3.1.1}, $\xi$ has only finitely many zeros
$p_i\in S$ with {\sl positive} multiplicities $\mu_i$. One can easily see
that $c_1(L)=\sum \mu_i \ge 0$. Consequently $\ssh^0_D (S,\, L)$ vanishes
if $c_1(L)<0$. The~vanishing result for $\ssh^1_D$ is obtained via the~Serre
duality of {\sl Lemma 2.3.2}.\qed

\bigskip\bigskip
\noindent \sl 2.4. Deformation of cusp-curves.

\nobreak\medskip \rm
Recall that in (2.2.9) we have obtained the following long exact sequence

\smallskip
$$

\def\mapright#1{\smash{\mathop{\longrightarrow}\limits^{#1}}}
\def\mapdown{\Big\downarrow}
\matrix{
0& \mapright{}& \ssh^0(S, TS\otimes [A]) &\mapright{} & \ssh^0_D(S, E)
 & \mapright{}& \ssh^0_D(S, N_0)   &\mapright{\delta} &\cr
\vphantom{.}&&&&&&&&\cr
 & \mapright{}& \ssh^1(S, TS\otimes [A]) &\mapright{} & \ssh^1_D(S, E)
 & \mapright{}& \ssh^1_D(S, N_0)         &\mapright{} &0.
}
$$
It is more important for us, that we can associated a similar long exact
sequence of $D$-cohomologies also to the short exact sequence (2.2.6). Note,
that due to the {\sl Lemmas 2.2.1, 2.2.2} and {\sl Corollary 2.2.3} we obtain
the short exact sequence of complexes
$$

\setbox1=\hbox{$\to$}
\def\mapright#1{\smash{\mathop{{\hbox to
\wd1{\hss\hbox{$\displaystyle\longrightarrow$}\hss}}}\limits^{#1}}}
\def\mapdown#1{\Big\downarrow\rlap{$\vcenter{\hbox{$\scriptstyle#1$}}$}}
\matrix{0&\mapright{}&L^{1,p}(S, TS)&\mapright{du}&
L^{1,p}(S, E)&\mapright{\overline\pr}&
L^{1,p}(S, E)\bigm/du(L^{1,p}(S, TS))
&\mapright{}&0\cr
 & &\mapdown{\overline \partial_S}& &\mapdown{D}& &\mapdown{\overline D}& & \cr
0&\mapright{}&L^p_{(0,1)}(S, TS)&\mapright{du}&
L^p_{(0,1)}(S, E)&\mapright{}&
L^p_{(0,1)}(S, E)\bigm/
du(L^p_{(0,1)}(S, TS))&\mapright{}&0\cr\vphantom{.}\cr}
\eqno(2.4.1)
$$

\smallskip\noindent
where $\overline D$ is induced by $D\equiv D_{u, J}$.

\medskip\noindent
{\bf Lemma 2.4.1.} {\it For $\overline D$ just defined,
$\ker \overline D=\ssh^0_D(S, N_0) \oplus \ssh^0(S, {\cal N}_1)$ and
$\coker \overline D=\ssh^1_D(S, N_0)$.
}

\smallskip\noindent
\bf Proof. \rm Let $\pr_0: {\cal O}(E) \to {\cal O}(N_0)$ and
$\pr_1: {\cal O}(E) \to {\cal N}_1$ denote the natural projections induced
by $\pr: {\cal O}(E) \to {\cal N}\equiv {\cal O}(N_0) \oplus {\cal N}_1$.
Let also $A$ as in (2.2.7) denote the support of ${\cal N}_1$, {\sl i.e.} the
finite set of the vanishing points of $du$. Then $\pr_0$ defines maps
$$
\matrix{
\pr_0:L^{1,p}(S , E)  &  \longrightarrow
& L^{1,p}(S , N_0)
\cr
\vphantom{.}\cr
\pr_0:L^p_{(0, 1)}(S , E)  &  \longrightarrow
& L^p_{(0, 1)}(S, N_0)
\cr}
\eqno(2.4.2)
$$
Further, in a neighbourhood of every point $p\in A$ the sequence (2.2.6)
can be represented in the form
$$
0\longrightarrow {\cal O}
\buildrel \alpha_p \over\longrightarrow {\cal O}^{n-1}
\oplus {\cal O}
\buildrel \beta_p \over\longrightarrow {\cal O}^{n-1}
\oplus{\cal N}_1\vert_p
\longrightarrow 0 \eqno(2.4.3)
$$
with $\alpha_p(\xi)= (0, z^{\nu_p} \xi)$ and $\beta(\xi,\eta)=
(\xi, j^{(\nu-1)}_p\eta)$. Here $z$ denotes local holomorphic coordinate
on $S$ with $z(p)=0$, $\nu_p$ is the multiplicity of $du$ in $p$,
$j^{(\nu-1)}_p\eta$ is a $(\nu-1)$-jet of $\eta$ in $p$, and ${\cal N}_1
\vert_p$ is a stalk of ${\cal N}_1$ in $p$.

 Now let $\overline\xi\in L^{1,p}(S, E)/du(L^{1,p}(S, TS))$ and
satisfies $\overline D(\overline\xi)=0$. This means that $\overline\xi$ is
represented by some $\xi \in L^{1,p}(S, E)$ with $D\xi= du(\eta)$ for some
$\eta \in L^p_{(0, 1)}(S, TS)$. It is obvious that there exists $\zeta
\in L^{1,p}(S, TS)$ such that $\overline\partial \zeta = \eta$ in a
neighbourhood of $A$. Consequently, $D(\xi- du(\zeta))=0$ in a neighbourhood
of $A$. Denote $\xi_1 := \xi -du(\zeta)$. Due to (2.2.1), in a neighbourhood
of $p\in A$ the equation $D\xi_1=0$ is equivalent to
$$
\overline\partial \xi_1+ N(\xi_1, du)=0.
$$
Due to the {\sl Lemma 3.1.2} below, $\xi_1= P(z)+o(|z|^{\nu_p} )$ with some
(holomorphic) polynomial $P(z)$. This gives possibility to define
$\varphi_0 := \pr_0(\xi) \in L^{1,p}(S, N_0)$ and $\varphi_1 :=
\pr_1(\xi- du(\zeta)) \in \ssh^0(S, {\cal N}_1)$. Due to (2.2.7) and
(2.2.8), $D_N\varphi_0=0$. If $\overline\partial \zeta' = \eta$
in a neighbourhood of $A$ for some other $\zeta'\in L^{1,p}(S, TS)$,
then $\zeta - \zeta'$ is holomorphic in a neighbourhood of $A$, and
consequently $\pr_1(du(\zeta - \zeta'))=0$. Thus the map
$\iota^0: \ker \overline D \to \ssh^0_D(S, N_0) \oplus \ssh^0(S, {\cal N}_1)$,
$\iota^0(\overline\xi)=(\varphi_0,\varphi_1)$ is well-defined.

\smallskip
   Assume that $\iota^0(\overline\xi)=0$ for some $\overline\xi\in
\ker\overline D$ and that $\overline\xi$ is represented by $\xi\in
L^{1,p}(S, E)$ with $D\xi= du(\eta)$ for some $\eta \in
L^p_{(0, 1)}(S, TS)$. Let $\zeta \in L^{1,p}(S, TS)$ satisfies
$\overline\partial \zeta = \eta$ in a neighbourhood of $A$. The assumption
$\iota^0(\overline\xi)=0$ implies that $\pr_1(\xi- du(\zeta))=0$ and that
$\pr_0(\xi- du(\zeta))=0$ in a neighbourhood of $A$. Consequently, $\xi-
du(\zeta)=du(\psi)$ for some $\psi \in L^{1,p}(S, TS)$ and $\xi \in
du(L^{1,p}(S, TS))$. This means that $\overline\xi=0\in \ker\overline D$
and $\iota^0$ is injective.

\smallskip
   Let $\varphi_0 \in \ssh^0_D(S, N_0)$ and $\varphi_1 \in \ssh^0(S, {\cal
N}_1)$. For every $p\in A$ fix a neighbourhood $U_p$ and representation of
(2.2.6) in form (2.4.3) over $U_p$. In every $U_p$ we can find $\xi=(\xi_0,
\xi_1) \in L^{1,p}(U_p,\cc^{n-1} \times\cc)$ satisfying the following
properties:

\smallskip
 {\sl a)} $D\xi=0$;

\smallskip
 {\sl b)} $\xi_0$ coincide with $\varphi_0\vert_{U_p}$ under
the identification ${\cal O}(N_0)\vert_{U_p} \cong {\cal O}^{n-1}$;

\smallskip
 {\sl c)} $j^{(\nu-1)}_p \xi = \varphi_1\vert_p \in
{\cal N}_1 \vert_p$.

\smallskip
 Corresponding $\xi_1 \in L^{1,p}(U_p,\cc)$ can be
constructed as follows. Let $D(\xi_0,\xi_1) = (\eta_0,\eta_1)$. From
$\pr_0(D\xi)= D_N(\pr_0\xi)=0$ one has $\eta_0=0$. In the representation
(2.4.3) the identity $R\scirc du =0$ of {\sl Lemma 2.2.1} means that
$D(0,\xi_1) = (0,\dbar\xi_1)$. So one can find $\xi_1$ with $\dbar\xi_1
=-\eta_1$, which gives $D(\xi_0,\xi_1) =0$. Consequently, $j^{(\nu-1)}_p
\xi$ is well-defined. Adding an appropriate {\sl holomorphic} term to
$\xi_1$ one can satisfy the condition {\sl c)}.
Using an appropriate partition of unity we can construct $\xi'
\in L^{1,p}(S, E)$ such that $\xi'$ coincide with $\xi$ in a
(possibly smaller) neighbourhood of every $p\in A$ and such that $\pr_0 \xi'
= \varphi_0$ in $S$. Thus from $D_N\varphi_0=0$ and (2.2.8) one
get $D\xi' \in du(L^p_{(0, 1)}(S, TS \otimes [A]))$. But $D\xi'
=0$ in a neighbourhood of every $p\in A$, which means that $D\xi' \in du(
L^p_{(0, 1)}(S, TS))$. This shows surjectivity of $\iota^0: \ker
\overline D \to \ssh^0_D(S, N_0) \oplus \ssh^0(S, {\cal N}_1)$.

\medskip
Now consider the case of $\coker \overline D$. Since the operator
$\pr_0$ satisfies the identities $\pr_0\circ D_E=D_N\circ \pr_0$ and
$\pr_0\circ du=0$, the induced map
\smallskip$\displaystyle
\iota^1 : \coker \overline D\equiv L^p_{(0,1)}(S, E)\bigm/
\bigl(D_E(L^{1,p}(S, E) \oplus  du(L^p_{(0,1)}(S, TS)\bigr)
\longrightarrow $

\smallskip \hfill$\displaystyle
\longrightarrow\ssh^1_D(S, N_0)\equiv
L^p_{(0,1)}(S, N_0)/ D_N(L^{1,p}(S, N_0)),
\qquad$

\smallskip\noindent
is well-defined. More over, the surjectivity of $\pr_0:
L^p_{(0,1)}(S, E) \to L^p_{(0,1)}(S, N_0)$
easily yields the surjectivity of $\iota^1$.

Assume now that $\iota^1 \overline\xi =0$ for some $\overline\xi \in
\coker \overline D$, and that $\overline \xi$ is represented by $\xi
\in  L^p_{(0,1)}(S, E)$. Then the condition $\iota^1
\overline\xi=0$ means that $\pr_0 \xi= D_N \eta$ for some $\eta
\in L^{1,p}(S, N_0)$. Find $\zeta \in L^{1,p}(S, E)$, such that
$\pr_0 \zeta =\eta$. Then $\pr_0 (\xi - D\zeta)=0$, and due to
(2.2.8), $\xi- D\zeta=du(\varphi)$ for some $\varphi \in
L^p_{(0,1)}(S, TS\otimes [A])$. Find $\psi \in
L^{1,p}(S, TS\otimes [A])$, such that $\dbar \psi =\varphi$
in some neighbourhood of $A$. Then
$$
\xi- D\zeta - D(du(\psi))= du(\varphi-\dbar\psi) \in
du(L^p_{(0,1)}(S, TS)).
$$
Consequently, $\overline \xi \in \im\overline D$. This shows the injectivity
of $\iota^1$.\qed

\bigskip
\state Corollary 2.4.2. \it The~short exact sequence $(2.2.6)$ induces
the~long exact sequence of $D$-cohomologies
$$

\def\mapright#1{\smash{\mathop{\longrightarrow}\limits^{#1}}}
\def\mapdown{\Big\downarrow}
\matrix{
0& \mapright{}& \ssh^0(S, TS) &\mapright{} & \ssh^0_D(S, E)
 & \mapright{}& \ssh^0_D(S, N_0)\oplus \ssh^0(S, {\cal N}_1)
 &\mapright{\delta} &\cr
\vphantom{.}&&&&&&&&\cr
 & \mapright{}& \ssh^1(S, TS) &\mapright{} & \ssh^1_D(S, E)
 & \mapright{}& \ssh^1_D(S, N_0)         &\mapright{} &0.
}
$$

\medskip\rm
Now suppose that  $(u, J_S, J) \in {\cal P}$, {\sl i.e.} $u$ belongs to
$C^1(S, X)$ and satisfies the equation $du\scirc J_S = J \scirc
du$. Set
$$
T_{(u, J_S, J)}{\cal P}\deff \bigl\{\, (v,\dot J_S,\dot J) \in T_u{\cal S}
\times T_{J_S} {\cal J}_S \times T_J{\cal J} \, :\,
2D_{u, J}v  +  \dot J \scirc du \scirc J_S +
 \dot J_S \scirc du \scirc J =0 \,\bigr\}.
$$
Let $\pr_{\cal J}: {\cal S} \times {\cal J}_S \times {\cal J} \to {\cal J}$
and  $\pr_{(u, J_S, J)} : T_{(u, J_S, J)} {\cal P} \to T_J{\cal J}$ denote
the natural projections.

\smallskip\nobreak
\state Theorem 2. {\it The map $\pr_{(u, J_S, J)} : T_{(u, J_S, J)}{\cal P}
\to T_J{\cal J}$ is surjective \iff\/  $\ssh^1_D(S, N_0)=0$, and then the
following holds:

\sli  the kernel $\ker(\pr_{(u, J_S, J)})$ admits a closed complementing
space;

\slii for $(\tilde u,\tilde J_S,\tilde J)\in {\cal P}$ close enough to
$(u, J_S, J)$ the projection $\pr_{(\tilde u,\tilde J_S,\tilde J)}$ is
also surjective;

\sliii for some neighbourhood $U\subset {\cal S} \times {\cal J}_S
\times {\cal J}$ of $(u, J_S, J)$ the set ${\cal P} \cap U$ is a Banach
submanifold of $U$ with the tangent space $T_{(u, J_S, J)}{\cal P}$
at $(u, J_S, J)$;

\sliv there exists a $C^1$-map $f$ from some neighbourhood $V$ of
$J\in {\cal J}$ into $U$ with $f(V)\subset {\cal P}$, $f(J) = (u,J_S,J)$,
 and $\pr_{\cal J}
\scirc f = \id_V$, such that $\im (df: T_J{\cal J} \to
T_{(u, J_S, J)}{\cal P})$ is complementing to  $\ker(\pr_{(u, J_S, J)})$.
}

\state Proof. Denote by $\tilde A$ the set of singular point of $M=u(S)$, {\sl
i.e.} the set of cuspidal and self-intersection points of $M$.
Let $\xi \in L^p_{(0, 1)} (S, E)$. Using the local
solvability of the $D$-equation (proved essentially in {\sl Lemma 3.2.1} below)
and an appropriate partition of unity we can find $\eta\in L^{1,p}(S, E)$,
such that $\xi -2D\eta \in C^1_{(0, 1)}(S, E)$ and $\xi -2D\eta =0$ in a
neighbourhood of $\tilde A$. Then we find $\dot J\in T_J{\cal J} \equiv
C^1_{(0, 1)}(X, TX)$ such that $\dot J \scirc du \scirc J_S =\xi-
2D\eta$. The surjectivity of $\pr_{(u, J_S, J)}$ and {\sl Lemma 2.4.1} easily
yield the identity $\ssh^1_D(S, N_0)=0$.

\smallskip\sl
From now on and till the~end of the~proof we suppose that $\ssh^1_D(S, N_0)=0$.

\smallskip\rm
Let $ds^2$ be some Hermitian metric on $S$ and let ${\cal H}_{(0, 1)} \subset
C^1_{(0, 1)}(S, TS)$ be a space of $ds^2$-harmonic $TS$-valued
$(0, 1)$-forms on $S$. Then the~natural map ${\cal H}_{(0, 1)} \to
\ssh^1(S, TS)$ is an isomorphism. Further, due to the {\sl Corollary 2.4.2}
the~map
$$
g\deff(du, 2D): {\cal H}_{(0, 1)} \oplus L^{1,p}(S, E) \to
L^p_{(0, 1)}(S, E)
$$
is surjective and has a finite-dimensional  kernel. Consequently there exists
a~closed subspace $Y\subset {\cal H}_{(0, 1)} \oplus L^{1,p}(S, E)$,
such that $g : Y \to L^p_{(0, 1)}(S, E)$ is an isomorphism. Let
$$
h=(h_{TS},h_E): L^p_{(0, 1)}(S, E) \to Y\subset {\cal H}_{(0, 1)}
\oplus L^{1,p}(S, E)
$$
be the inversion of $g\vert_Y$.

Take $\dot J\in T_J{\cal J}= \{\, I\in C^1(X, {\ss End}(TX))\, :\,
JI+IJ=0\,\} \equiv C^1_{(0, 1)}(X, TX)$. Then $\dot J \scirc du
\scirc J_S$ lies in $C^0_{(0, 1)}(S, E)$. Let $h(\dot J \scirc du
\scirc J_S) =(\xi,\eta)$ with $\xi\in {\cal H}_{(0, 1)}
\subset C^1_{(0, 1)}(S, TS)$ and $\eta\in L^{1,p}(S, E)$.
Then we obtain that
$$
2D(-\eta) + \dot J \scirc du \scirc J_S +
J \scirc du ( J_S \xi)=0,\eqno(2.4.4)
$$
where we used the identity $J\scirc du \scirc J_S =-du$. Using the identity
$T_{J_S}{\cal J}_S= C^1_{(0, 1)}(S, TS)$
we conclude that the formula $F(\dot J)= \bigl(J_S h_{TS}(\dot J \scirc du
\scirc J_S), -h_E(\dot J \scirc du \scirc J_S),\dot J \bigr)$ defines a linear
{\sl bounded} operator $F: T_J{\cal J} \to T_{(u, J_S, J)}{\cal P}$, such that
$\pr_{(u, J_S, J)}\circ F = \id_{T_J{\cal J}}$. In particular, $\ssh^1_D(S,
N_0)=0$ implies the surjectivity of $\pr_{(u, J_S, J)}$.

\smallskip
The image of the just defined operator $F$ is closed, because the convergence
of $F(\dot J_n)$, $J_n\in T_J{\cal J}$, obviously yields the convergence of
$J_n=\pr_{(u, J_S, J)}\circ F(\dot J_n)$. One can easily see that $\im(F)$ is
a closed complementing space to $\ker(\pr_{(u, J_S, J)})$.

\smallskip
Further, for $(\tilde u,
\tilde J_S,\tilde J) \in {\cal P}$ close enough to $(u, J_S, J)$ the~map
$\tilde g \deff (d\tilde u, 2D_{\tilde u,\tilde J}) : Y \to
L^p_{(0, 1)}(S,\tilde E)$ is also an~isomorphism. This implies
the~surjectivity $\pr_{(\tilde u,\tilde J_S,\tilde J)}$.

\smallskip
The~statements \sliii and \sliv can be easily obtained from {\sl ii)} and
the implicit function theorem.\qed
\smallskip

\def\diff{{\cal D}\mskip -.6mu i\mskip -5.1mu f\mskip-6mu f}
\def\star{\mathop{\msmall{*}}}

Let now ${\cal J}_S$ denote the space of all
$C^{0,\alpha}$-continuous complex structures on $S$ with some $0<\alpha<1$ and
let $g$ be the genus of $S$. Denote by ${\cal G}$ the space $\diff^{1,\alpha }
(S)$ of all $C^{1, \alpha }$-diffeomorfisms of $S$.
This is a~Banach group with the $C^1$-continuous group operation given by
composition, and with the tangent space $T_\id{\cal G}= C^{1,\alpha}(S, TS)$.
One has also the natural $C^1$-actions of ${\cal G}$ on
${\cal S}$ and on ${\cal J}_S$ given by relations
$$
\eqalign{
h\in {\cal G}, u\in{\cal S}     &\longrightarrow h\star u\deff u\scirc h\inv
\cr
h\in {\cal G}, J_S\in{\cal J}_S  &\longrightarrow
h\star J_S\deff dh\scirc J_S\scirc dh\inv,
\cr}
\eqno(2.4.5)
$$
so that ${\cal P}\subset {\cal S}\times {\cal J}_S \times {\cal J}$ is
invariant with respect to the diagonal action of ${\cal G}$. The {\sl Moduli
space of nonparametrized pseudo-holomorphic curves} is now define as a~quotient
space ${\cal M}_{[\gamma], g}\deff {\cal P}/ {\cal G}$ with the natural
projection $\pr_{\cal M} :{\cal M}_{[\gamma], g} \to {\cal J}$, and the fiber
${\cal M}_{[\gamma], g, J}\deff  \pr_{\cal M}\inv(J)$ is called the {\sl Moduli
space of nonparametrized $J$-holomorphic curves}. Informally speaking, this is
a topologization of a set of $J$-holomorphic curves $M\subset X$, which can be
parameterized by $S$.

The relation $h\star J_S = J_S'$ is equivalent to $(J_S, J_S')$-holomorphicity
of the map $h: S\to S$. So due to the regularity property of $\dbar$-operator,
see the {\sl Corollary 3.2.2}, the chosen smoothnesses of complex structures
and diffeomorphisms of $S$ are more convenient for description of the Moduli
spaces then, for example, $J_S \in C^1$.

The Moduli space
${\cal M}_{[\gamma], g, J}$ appears as an important object in the theory of
{\sl quantum cohomology}, which is also known as {\sl Gromov-Witten invariants
theory}, see [McD-S].
One of the problems, arising by studying Moduli spaces, is to obtain the
regularity property of ${\cal M}_{[\gamma], g, J}$ for given $[\gamma]$, $S$,
and $J$.

One get a~partial result in this directions using {\sl Lemma 2.4.1} and
{\sl Theorem 2}.

\state Corollary 2.4.3. \it  Let $u: (S, J_S) \to (X, J)$ be a nonconstant
irreducible pseudo-holomorphic map, such that $\ssh^1_D(S, N_0)=0$. Than in
a neighbourhood of $M\deff u(S)$ the Moduli space
${\cal M}_{[\gamma], g, J}$ is a~manifold with the tangent space $T_M
{\cal M}_{[\gamma], g, J} = \ssh^0_D(S, N_0) \oplus \ssh^0(S, {\cal N}_1)$.

\state Proof. \rm  The considerations similar to ones from the {\sl Theorem 2}
shows that ${\cal P}_J$ is a~Banach manifold with the tangent space
$$
T_{(u, J_S)}\vph{\cal P}_J =
\{\, (v, \dot J_S,) \in L^{1,p}(S, E) \times C^{0,\alpha}_{(0,1)}(S, TS)
\; : \; 2D_{u, J}\vph v + J\scirc du \scirc \dot J_S =0 \},
$$
which is linearly (not topologically!) isomorphic to
$$
\{\, v \in L^{1,p}(S, E) : D_{u,J}v\in du (C^{0,\alpha}_{(0,1)}(S, TS)) \,\}.
$$

Further, due to {\sl Corollary 2.2.3} the group ${\cal G}$ acts {\sl freely}
on all the pairs $(\tilde u, \tilde J_S)$, which are close enough to the
given $(u, J_S)$. So the tangent space to the orbit
${\cal G}\star (\tilde u, \tilde J_S) \deff \{\,(h\star \tilde u,
h\star \tilde J_S) : h\in{\cal G}\,\}$ in the point $(\tilde u, \tilde J_S)$
is equal to
$$
T_{(\tilde u, \tilde J_S) } {\cal G}\star (\tilde u, \tilde J_S)=
\{\, (\,d\tilde u(\xi), -2\tilde J_S \cdot \xi)
\in L^{1,p}(S,\tilde E) \times C^{0,\alpha}_{(0,1)}(S, TS)
 :  \xi \in C^{0,\alpha}_{(0,1)}(S, TS) \,\}.
$$
Here we use the relation (2.2.4) and the fact, that for $C^1$-curve $h_t$ in
${\cal G}$ with $h_0=\id$ and $\dot h_0 =\xi \in C^{1,\alpha}(S, TS)$ one has
$$
\msmall{d\over dt} (h_t\star J_S) = \nabla\xi \scirc J_S -J_S \scirc \nabla\xi
\equiv -2J_S\dbar_{J_S}\vph\xi.
$$

\medskip
Consider the quotient spaces
$$\eqalign{
&   Z_1\deff \{\, v \in L^{1,p}(S, E) : D_{u,J}v\in
du(C^{0,\alpha}_{(0,1)}(S, TS)) \,\}  \!\!\bigm/\!\!
du(C^{1,\alpha}(S, TS)),
\cr
&   Z_2\deff \{\, v \in L^{1,p}(S, E) : D_{u,J}v\in
du(L^p_{(0,1)}(S, TS)) \,\}  \!\!\bigm/\!\!
du(L^{1,p}(S, TS)),
\cr}
$$
with the natural map $\iota Z_1 \to Z_2$. Using methods, similar to ones from
the proof of {\sl Theorem 2}, one can show that $\iota$ is an isomorphism. So
every complementing space $Z$ to $T_{(u,J_S)}{\cal G} \star (u,J_S)$ in
$T_{(u, J_S)}\vph{\cal P}_J$ is finite-dimensional, closed, and isomorphic
to $\ssh^0_D(S, N_0) \oplus \ssh^0(S, {\cal N}_1)$. One can define a chart on
${\cal M}_{[\gamma], g, J}$, taking an appropriate submanifold $U \subset
{\cal P}_J$ with $(u,J_S) \in U$ and $T_{(u,J_S)}U$ complementing to
$T_{(u,J_S)}{\cal G} \star (u,J_S)$, and cosidering the restriction of the
projection ${\cal P}_J \longrightarrow {\cal M}_{[\gamma], g, J}$.
\qed

\vskip 0pt plus25pt
\bigskip
\bigskip
\noindent{\bf 3. Nodes of pseudo-holomorphic curves.}

\nobreak\smallskip\noindent \sl 3.1. Unique continuation lemma.

\nobreak\smallskip
\rm We start with the following unique continuation-type lemma
(compare with [Ar] and [Hr-W]):

\medskip
\state Lemma 3.1.1. {\it Suppose that the~function $f\in L^2_\loc(\Delta,
\cc^n)$ is not identically $0$, $\dbar f\in L^1_\loc(\Delta,\cc^n)$
and satisfies {\sl a.e.} the~inequality
$$
\vert \dbar f \vert \le h\cdot \vert f \vert\eqno(3.1.1)
$$
for some nonnegative $h\in L^p_\loc(\Delta)$ with $2<p<\infty$. Then

     \sli $f\in L^{1, p}_\loc(\Delta)$, in particular $f\in C^{0,\alpha}_\loc
(\Delta)$ with $\alpha\deff1-{2\over p}$;

     \slii for any $z_0\in\Delta$ such that $f(z_0)=0$ there exists
$\mu\in\nn$---the~multiplicity of zero of $f$ in $z_0$---such that
$f(z)=(z-z_0)^\mu \cdot g(z)$ for some $g\in L^{1, p}_\loc (\Delta)$ with
$g(z_0)\not=0$.
}

\bigskip
\state Proof. \sli is easily obtained by increasing smoothness argument, see
\eg [Mo]. Really, from (3.1.1) and H\"older inequality one gets that $\bar
\partial f\in L^{2p\over p+2}_\loc(\Delta)$. Consequently $f\in L^{1, {2p\over
p+2}}_\loc(\Delta)$ due to ellipticity of $\dbar$ and since ${2p\over p+2} >
1$. Again by Sobolev imbedding we have $f\in L^p_\loc(\Delta)$. Suppose now
that $f\in L^q_\loc(\Delta)$ with some $q>2$. If ${pq\over p +q} <2$, then
$f\in L_\loc^{\tilde q}$ with $\tilde q= 1/({1\over q} + {1\over p} - {1\over
2} ) > q$ by the~Sobolev imbedding theorem. Repeating this argument several
times we obtain that $f\in L^{1,\tilde q} _\loc(\Delta)$ for some $\tilde
q>2$ which implies that $f$ is continuous. Again by (3.1.1) $\dbar f \in
L^p_\loc(\Delta)$ and $f \in L^{1, p}_\loc(\Delta)$. Another application of
the~Sobolev imbedding theorem, this time for $2<p<\infty$, yields H\"older
$\alpha$-continuity.

\smallskip
    \slii  Now suppose that $f(z_0)=0$. Then due to the~H\"older continuity we
have $\vert f(z) \vert \le C\vert z-z_0 \vert^\alpha$ for $z$ close enough to
$z_0$ and consequently $f_1(z) \deff f(z) / (z-z_0)$ is from $L^2_\loc(\Delta)
$. The~theorem of Harvey-Polking [Hv-Po] provides that $\dbar f_1\in
L^1_\loc(\Delta)$ and $f_1$ also satisfies
inequality (3.1.1). In particular, $f_1$ is also continuous. Iteration of this
procedure gives a possibility to define the~multiplicity of zero of $f$ in
$z_0$ provided we show that after  finite number of steps we obtain the
function $f_N$ with $f_N(z_0)\not= 0$. To do this we may assume that $z_0=0$.
Let $\pi_\tau (z)\deff\tau\cdot z$ for $0<\tau<1$. Then $f_\tau\deff
\pi_\tau^*(f)$ satisfies inequality $\vert \dbar f_\tau \vert \le
\tau\pi_\tau^*h \cdot \vert f_\tau \vert$. Since
$$
\Vert \tau\pi_\tau^*h \Vert_{L^p(\Delta)} =
\tau^{1-2/p}\cdot \Vert h \Vert_{L^p(\pi_\tau(\Delta))}
$$
we can also assume that $\Vert h \Vert_{L^p(\Delta)}$ is small enough.
Fix a~cut-off function $\varphi\in \cinf_0(\Delta)$ which is identically $1$
in ${1\over2}\Delta$, the~disk of radius $1\over2$. Then
$$
\Vert \varphi f \Vert_{L^p(\Delta)}\le
C_1\cdot \Vert \varphi f \Vert_{L^{1, {2p\over2+p}}(\Delta)}\le
C_2\cdot \Vert \dbar(\varphi f) \Vert_{L^{2p\over2+p}(\Delta)}\le
$$
$$
\le C_2\cdot (\Vert \varphi \dbar f \Vert_{L^{2p\over2+p}(\Delta)}+
   \Vert \dbar\varphi  f \Vert_{L^{2p\over2+p}(\Delta)})\le
   C_2\cdot ( \Vert \varphi h f  \Vert_{L^{2p\over2+p}(\Delta)}
   +  \Vert \dbar\varphi f \Vert_{L^{2p\over2+p}(\Delta)} \le
$$
$$
    \le  C_3(\Vert \varphi f \Vert_{L^p(\Delta)}\cdot \Vert h \Vert _
    {L^p(\Delta )} + \Vert f\Vert _{L^p(\Delta \setminus {1\over2}\Delta )}
 \cdot \Vert \dbar\varphi \Vert _{L^2(\Delta \setminus {1\over2}\Delta )}).
$$

\smallskip
Here we used the~fact that
the~support of $\dbar\varphi$ lies in $\Delta\bss{1\over2}\Delta$. Since
$\Vert h  \Vert_{L^p(\Delta)}$ is small enough we obtain the~estimate
$$
\Vert f \Vert_{L^p({1\over3}\Delta)}
\le C\cdot \Vert f \Vert_{L^p(\Delta\setminus {1\over2}\Delta)}
$$
with the~constant $C$ independent of $f$. Thus if the~multiplicity of
zero of $f$ in $z_0=0$ is at least $\mu$, then
$$
\Vert z^{-\mu}f \Vert_{L^p({1\over3}\Delta)}
\le C\cdot \Vert z^{-\mu}f \Vert_{L^p(\Delta\bss{1\over2}\Delta))}
$$
which easily gives
$$
\Vert f \Vert_{L^p({1\over3}\Delta)} \le C
\left(\hbox{$2\over3$}\right)^{-\mu}
\cdot \Vert f \Vert_{L^p(\Delta\bss{1\over2}\Delta))}.
$$
Now one can easily see that either $f$ has isolated zeros of finite
multiplicity, or $f$ is identically zero. Yet the~last case is excluded by
hypothesis of the~lemma. \qed

\medskip
\state Lemma 3.1.2. {\it
Under the~hypothesis of {\sl Lemma 3.1.1} suppose
additionally that $f$ satisfies {\sl a.e.} the inequality
$$
\vert \dbar f(z) \vert \le \vert z-z_0 \vert^\nu
h(z)\cdot \vert f(z) \vert,\eqno(3.1.2)
$$
with $z_0\in \Delta$, $\nu\in\nn$, and $h\in L^p_\loc(\Delta)$, $2<p<\infty$.
Then
$$
f(z)=(z-z_0)^\mu\bigl(P^{(\nu)}(z) + (z-z_0)^\nu g(z)\bigr),
\eqno(3.1.3)
$$
where $\mu\in\nn$ is the multiplicity of $f$ in $z_0$, defined above,
$P^{(\nu)}$ is a polynomial in $z$ of degree $\le \nu$ with $P^{(\nu)}(z_0)
\not=0$, $g\in L^{1, p}_\loc(\Delta,\cc^n)\hookrightarrow C^{0,\alpha}$,
$\alpha=1-{2\over p}$, and $g(z)=O(\vert z-z_0\vert^\alpha)$.
}

\medskip\noindent
\bf Proof \rm uses the~same idea as in the~previous lemma. Define $f_0(z)
\deff \msmall{f(z)\over (z-z_0)^\mu}$ and $h_1(z) \deff h(z)\cdot \vert f_0(z)
\vert$. Due to {\sl Lemma 3.1.1}, $f_0\in C^{0,\alpha}$, $f_0(z_0)\not=0$,
$h_1\in L^p_\loc$, and $f_0$ satisfies {\sl a.e.} the~inequality $\vert \dbar
f_0(z) \vert \le \vert z-z_0 \vert^\nu h_1(z)$.

 Set $a_0=f_0(z_0)$. Since $f_0(z)-a_0 = O(\vert z-z_0 \vert^\alpha)$,
we have $f_1\deff\msmall{f_0(z) - a_0 \over z-z_0}\in L^2_\loc$. Applying
the~theorem of Harvey-Polking once more we obtain that $\vert \dbar f_1 \vert
\le \vert z-z_0 \vert^{\nu-1} h_1$, and consequently $f_1\in C^{0,\alpha}$,
$f_1(z)-f_1(z_0)=O(\vert z-z_0 \vert^\alpha)$.
Repeating this procedure $\nu$ times we obtain the~polynomial
$$
P^{(\nu)}(z)=a_0 + (z-z_0)a_1 + \cdots + (z-z_0)^\nu a_\nu
$$
with
$$
a_k\deff\lim_{z\to z_0}
\msmall{ f(z)-\sum_{i=0}^{k-1} (z-z_0)^i a_i \over (z-z_0)^k },
\qquad 0\le k\le\nu,
$$
and the~function
$$
g(z)\deff \msmall{f(z)-P^{(\nu)}(z) \over (z-z_0)^\nu },
$$
which satisfies the~conclusion of the~lemma. \qed

\medskip

\state Corollary 3.1.3. {\it Let $J$ be a Lipschitz-continuous almost complex
structure in a neighbourhood $U$ of $0\in \cc^n$, such that $J(0)=J\st$,
the~standard complex structure in $\cc^n$. Suppose that $u:\Delta\to U$
is a $J$-holomorphic $C^1$-map with $u(0)=0$. Then there exist uniquely
defined $\mu\in\nn$ and a (holomorphic) polynomial $P^{(\mu-1)}$ of degree
$\le\mu-1$, such that $u(z)=z^\mu\cdot P^{(\mu-1)} +z^{2\mu -1}v(z)$
with $v(0)=0$ and $v\in L^{1,p}(\Delta,\cc^n)$ for any $p<\infty$.
}

\medskip
\state Proof. Really, by the hypothesis of the lemma $du + J(u)
\scirc du \scirc J_\Delta =0$ and hence
$$
\vert \dbar u \vert  =
\bigl\vert \msmall{1\over2}( du + J\st \scirc du \scirc J_\Delta)\bigr\vert =
\bigl\vert \msmall{1\over2} (J\st - J(u))\scirc du \scirc J_\Delta \bigr\vert
 \le \Vert J \Vert_{C^{0,1}(B)} \cdot \vert u \vert  \cdot
\vert du \vert.
\eqno(3.1.4)
$$

\smallskip\noindent
Thus by {\sl Lemma 3.1.1} $u(z) = z^\mu w(z)$ with some $\mu \in\nn$
and $w\in L^{1,p}_\loc(\Delta, \cc^n)$ for any $p<\infty$.
If $\mu =1$ we are done. Otherwise $ du(z)/z^{\mu-1} =
\mu\,w\,dz + z\,dw \in L^p_\loc(\Delta, \cc^n)$ for any $p<\infty$.
So the corollary follows now from {\sl Lemma 3.1.2}.
\qed

\bigskip
\state Corollary 3.1.4. {\it Let $u:S\to (X, J)$ be a~$J$-holomorphic map.
Then for any $p\in X$ the~set $u\inv(p)$ is discrete in $S$, provided that
$J$ is Lipschitz-continuous.
}

\smallskip
\state Proof. Take an~{\sl integrable} complex structure $J_1$ in some
neighbourhood $U$ of $p\in X$, such that $J(p)=J_1(p)$. Take also
$J_1$-holomorphic coordinates $w_1,\ldots, w_n$ in $U$, such that $w_i(p)=0$.
Then the~statement follows easily from {\sl Corollary 3.1.3.}\qed

\bigskip
\noindent \sl 3.2. Inversion of $\dbar_{J}+R$.

\smallskip\nobreak
\rm
Now suppose that $J\in C^0(\Delta, {\ss End}_\rr(\rr^{2n}))$ satisfies
the~identity $J^2\equiv -1$,\ie $J$ is a continuous almost complex structure
in
the trivial $\rr^{2n}$-bundle over $\Delta$. Define the~$\rr$-linear
differential
first order operator $\dbarj:L^{1, p}(\Delta,\rr^{2n}) \to L^p(\Delta,\rr^{2n})$
by setting
$$
\dbarj(f)=\msmall{1\over2}\left(
\msmall{\d f \over \d x} + J \msmall{\d f\over \d y} \right).\eqno(3.2.1)
$$
For example, for $J\equiv J\st$, the~standard complex structure in $\rr^{2n}=
\cc^n$, the~operator $\dbarj$ is usual Cauchy-Riemann operator $\dbar$.
The~operator $\dbarj$ is elliptic and possesses nice regularity properties in
Sobolev spaces $L^{k, p}$ with $1<p<\infty$ and H\"older spaces $C^{k,\alpha}$
with $0<\alpha<1$. The~following two statements are typical for (nonlinear)
elliptic PDE and gathers result which we need for purpose of this paper.

\smallskip
\state Lemma 3.2.1. {\it Let $J$ be $C^k$-continuous with $k\ge0$ and $\dbarj$
be as above. Let also $R$ be an~${\ss End}(\rr^{2n})$-valued function in
$\Delta$ of class $L^{k,p}$ with $1<p<\infty$. If $k=0$ we also assume that
$p>2$. Suppose that for $\dbarj f +Rf \in L^{k,p}$ for some $f\in L^{1,
1}(\Delta,\rr^{2n})$. Then $f\in L^{k+1,p}_\loc (\Delta,\rr^{2n})$ and
for $r<1$
$$
\Vert f \Vert_{L^{k+1, p}(\Delta(r))}
\le C_1(J, \Vert R\Vert_{L^{k,p}}, k, p, r)
\bigl( \Vert \dbarj f + Rf \Vert_{L^{k, p}(\Delta)} +
\Vert f \Vert_{L^1(\Delta)}\bigr)
 \leqno{\sl i)}
$$

\smallskip\noindent
If, in addition, $J$ and $R$ are $C^{k,\alpha}$-smooth with $0<\alpha<1$,
then $f\in C^{k+1,\alpha}_\loc(\Delta,\rr^{2n})$ and

$$
\Vert f \Vert_{C^{k+1,\alpha}(\Delta(r))}
\le C_2( J , \Vert R \Vert_{C^{k,\alpha}},
k, \alpha, r)\cdot \bigl(
\Vert \dbarj f + Rf \Vert_{C^{k,\alpha}(\Delta)}+
\Vert f \Vert_{L^1(\Delta)}\bigr)
\leqno{\sl ii)}
$$

\smallskip\noindent
If additionally  $\Vert J-J\st \Vert_{C^k(\Delta)} +
\Vert R \Vert_{L^{k,p}(\Delta)}$ $($resp.\ $\Vert J-J\st \Vert_{C^{k,\alpha}
(\Delta)}+ \Vert R \Vert_{C^{k,\alpha}(\Delta)}\,)$
is small enough, then there exists a~linear bounded operator
$T_{J, R}:L^{k,p}(\Delta,\rr^{2n}) \to L^{k+1, p}(\Delta,\rr^{2n})$
$($resp.\ $T_{J, R}:C^{k, \alpha}(\Delta,\rr^{2n}) \longrightarrow
C^{k+1,\alpha}(\Delta,\rr^{2n})$\thinspace{\rm)},
such that $(\dbarj + R)\scirc T_{J, R}\equiv \id$ and $T_{J,
R}(f)\vert_{z=0}=0$.
}

\bigskip
\state Proof. The~estimates \sli and \slii of the~theorem are obtained from
the~ellipticity of $\dbarj$ by standard methods, see [Mo] for the~case of
general elliptic systems or [Sk] for the~direct proof. We note only that
for $J$ being only continuous the~constant $C_1$ in \sli crucially depends on
modulus of continuity of $J$, whereas in all the~other cases of the~theorem
the~dependence is essentially only on the~corresponding norm of $J$.

\smallskip
     Let $T$ be a composition ${\d\over\d z} \scirc G$ where $G(f)$ is
the~solution of the~Poisson equation $\Delta u=f$ with the~Dirichlet boundary
condition on $\d\Delta$. Then for $J\equiv J\st$ and $R\equiv0$ we can define
$T_{J\st,0}(f) \deff T(f)- T(f)(0)$. In general case we set
$$
T_{J, R}=\sum_{n=0}^\infty (-1)^nT_{J\st, 0}\scirc
\bigl((\dbarj-\dbar_{J\st}+R)\scirc T_{J\st, 0}\bigr)^n.
$$
Due to the Sobolev imbedding theorem, the~series converges in an~appropriated
norm provided $\Vert J-J\st \Vert + \Vert R \Vert$ is small enough.
\par\qed

\bigskip
\state Corollary 3.2.2. {\it Let $k\in \nn$, $q>2$, and $J$
a~$C^k$-smooth almost complex
structure in $X$. Let also $(S, J_S)$ be a~complex curve. Suppose that
$L^{1,q}$-map $u:S\to X$ satisfies the~equation
$$
du + J \scirc du \scirc J_S = 0.
$$
Then $u$ is $L^{k+1,p}$-smooth for any $p<\infty$. If, in addition, $J$ is
$C^{k,\alpha}$-smooth with $0<\alpha<1$, then $u$ is $C^{k+1,\alpha}$-smooth.

Let $J^{(n)}$ (resp.\ $J^{(n)}_S$) be a sequence of almost complex
structures on $X$ (resp.\ on $S$), which $C^k$-converges to $J$
(resp.\ to $J^{(n)}_S$) and let $u_n:S \to X$ be a sequence of
$(J^{(n)}_S, J^{(n)})$-holomorphic maps. Then the $C^0$-convergence
$u_n \longrightarrow u$ implies  the $L^{k+1,p}$-convergence for
any $p<\infty$, and $C^{k+1, \alpha}$-convergence if
$(J^{(n)}_S, J^{(n)})$ converge to $(J_S, J)$ in $C^{k,\alpha}$,
$0<\alpha<1$.
}

\smallskip
\state Proof. The map $u$ is continuous and in local coordinates
$w_1,\ldots, w_{2n}$ on $X$ and $z=x+iy$ on $S$ the~equation has the~form
$$
du(z) + J(u(z)) \scirc du \scirc J_M = 0,
$$
which is equivalent to $\dbar_{J\scirc u} u =0$.  Using {\sl Lemma 3.2.1}
and induction in $k'=0\ldots k$, one can obtain the~regularity of $u$.

Similarly, for $J^{(n)}$ and $u^{(n)}$ satisfying the hypothesis of the
corollary one gets
$$
\dbar_{J\scirc u}(u^{(n)}-u)=
(\dbar_{J\scirc u} - \dbar_{J^{(n)} \scirc u^{(n)}})u^{(n)}
\longrightarrow 0 \quad\hbox{in $L^{k',p}_\loc$
(resp.\ in $C^{k',\alpha}_\loc$),}
$$
by induction in $k'=0\ldots k$.

\par\qed

\smallskip
\state Remark. The corollary implies that for a compact Riemannian surface
$S$ the topology in the space $\cal P$ of $(J_S, J)$-holomorphic maps $u:S\to
X$ is independent of the particular choice of the functional space
$L^{k', p}(S,X)\supset \cal P$ with $1\le k' \le k+1$, $1< p<\infty$, and
$kp>2$, provided $J_S$ and $J$ are changing $C^k$-smoothly. In the same
way, {\sl Lemma 3.2.1} implies that for $J\in C^k$ with $k\ge1$ the
differential structure on $\cal P$ is also independent of the particular
choice of a functional space $L^{k',p}(S,X)$ with $1\le k' \le k$,
$1<p<\infty$ and $k'p>2$.


\medskip
\state Lemma 3.2.3. {\it For any natural numbers $\mu>\nu\ge0$ and real
numbers $p>2$, $\alpha<{2\over p}$ and $\gamma>0$ there exists $C =
C(\mu, \nu, \alpha, p, \gamma)>0$ with the following property.
Let $J$ be an almost complex structure in
$B \subset \cc^n$ with $J(0)=J\st$ and let $u:\Delta \to B$ be a
$J$-holomorphic map of the form $u(z)=z^\mu(P(z) +z^\nu v(z))$, where
$P(z)$ is some (holomorphic) polynomial of the degree $\nu$ and
$v\in L^{1,p}(\Delta, \cc^n)$ with $v(0)=0$. Suppose that $\Vert J-J\st
\Vert_{C^1(B)} \le \gamma$ and $\Vert u\Vert_{L^{1,p}(\Delta)} \le \gamma$.
Then for any $0<r<{1\over2}$
$$
\bigl\Vert z dv \bigr\Vert_{C^0(\Delta(r))} +
\bigl\Vert z dv \bigr\Vert_{L^{1,p}(\Delta(r))}
\le
C \cdot r^\alpha\cdot \Vert u \Vert_{L^{1,p}(\Delta)},
\eqno(3.2.2)
$$

In particular, $du(z) = d(z^\mu P(z)) + o(|z|^{\mu+ \nu -1 +\alpha})$
for any $\alpha<1$.
}

\state Proof. For $0<r<{3\over4}$ we define the map $\pi_r:B\to B$, setting
$\pi_r(w) \deff r^\mu w$. We also set $J^{(r)} \deff \pi_r^*J$, $u^{(r)}(z)
\deff \pi_r^{-1} \scirc u(rz)$, $P^{(r)}(z)\deff P(rz)$ and $v^{(r)}(z)=
r^\nu v(rz)$. Then $J^{(r)}$ is an almost complex
structures in $B$ with $\Vert J^{(r)} -J\st \Vert_{C^1(B)} \le r^\mu
\Vert J-J\st\Vert_{C^1(B)}$, and $u^{(r)} \equiv z^\mu (P^{(r)}(z) +
z^\nu v^{(r)}(z))$ is a $J^{(r)}$-holomorphic.

Without losing generality we can suppose that $\alpha>0$. Set $\beta \deff
1+\alpha - {2\over p}$ and $q\deff {2\over1-\beta}$. Then $\alpha<\beta<1$,
$\beta + {2\over p} -1=\alpha$, and $q>2$. The {\sl Lemma 3.1.3} implies that
$\Vert v \Vert_{C^{0,\beta}(\Delta(2/3))} + \Vert dv \Vert_{L^q(\Delta(2/3))}
\le C_1 \cdot \Vert u \Vert_{L^{1,p}(\Delta)}$. Here the constant $C_1$, as
also constants $C_2,\ldots,C_6$ below, depend only on $\mu, \nu, p, \alpha$,
and $\gamma$, but are independent of $r$. Consequently,
$$
\Vert u^{(r)} - z^\mu P^{(r)}(z) \Vert_{L^{1,q}(\Delta)}
\le C_2\cdot r^{\nu + \beta} \cdot
\Vert u \Vert_{L^{1,p}(\Delta)}.
$$
Further, due to the {\sl Corollary 3.2.2} we have $\Vert u^{(r)}
\Vert_{L^{2,p}(\Delta)} \le C_3\cdot \Vert u \Vert_{L^{1,p}(\Delta)}$. Thus
$$
\Vert \dbar_{J\st}( u^{(r)} - z^\mu P^{(r)}(z) ) \Vert_{L^{1,p}(\Delta)} =
\Vert (\dbar_{J\st} -\dbar_{J^{(r)} \circ u^{(r)} } ) ( u^{(r)}) \Vert_{L^{1,
p}(\Delta)}
\le C_4\cdot r^\mu \cdot
\Vert u \Vert_{L^{1,p}(\Delta)}.
$$
Applying {\sl Lemma 3.2.1} we obtain
$$
\bigl\Vert z dv^{(r)} \bigr\Vert_{L^{1,p} (\Delta(2/3))} \le
C_5 \cdot r^{\nu + \beta} \cdot \Vert u \Vert_{L^{1,p}(\Delta)},
\eqno(3.2.3)
$$
\smallskip\noindent which is equivalent to
$$
\bigl\Vert z dv \bigr\Vert_{L^{1,p}(\Delta(2r/3))} \le
C_5 \cdot r^{\beta + 2/p -1} \cdot \Vert u \Vert_{L^{1,p}(\Delta)}.
$$
\smallskip\noindent On the other hand, $(3.2.3)$ implies that
$$
\bigl\Vert z dv^{(r)} \bigr\Vert_{C^0 (\Delta(2/3))} \le
C_6 \cdot r^{\nu + \beta} \cdot \Vert u \Vert_{L^{1,p}(\Delta)},
$$
\smallskip\noindent and consequently
$$
\bigl\Vert z dv \bigr\Vert_{C^0(\Delta(2r/3))} \le
C_6 \cdot r^\beta \cdot \Vert u \Vert_{L^{1,p}(\Delta)}.
$$

\qed


\smallskip
\state Lemma 3.2.4. \it
For given $p$, $2<p<\infty$, and $\gamma>0$ there exist constants
$\epsi=\epsi(p,\gamma)$ and $C=C(p,\gamma)$ with the following
property. Suppose that $J$ is a $C^1$-smooth almost complex
structure in the unit ball $B \subset \cc^n$
with $\Vert J-J\st \Vert_{C^1(B)}
\le \epsi$ and $u\in  L^{1,p}(\Delta, B(0,{1\over2}))$ is
a $J$-holomorphic map with $\Vert u \Vert_{L^{1,p}(\Delta)}
\le \gamma$. Then for every almost complex structure $\tilde J$ in
$B$ with $\Vert \tilde J-J \Vert_{C^1(B)} \le \epsi$
there exists a $\tilde J$-holomorphic map $\tilde u: \Delta \to
B$ with
$$
\Vert \tilde u-u \Vert_{L^{1,p}(\Delta)} \le
C\cdot \Vert \tilde J-J \Vert_{C^1(B)},
$$
such that $\tilde u(0) = u(0)$.

\smallskip\rm
\state Proof. Let $J_t$ be a curve of $C^1$-smooth almost
complex structures in $B$, starting at $J_0=J$ and depending
$C^1$-smoothly in $t\in [0,1]$. Consider an ordinary differential
equation for $u_t \in L^{1,p} (\Delta, B)$ with initial
condition $u_0=u$ and
$$
\msmall{du_t\over dt} = - T_{J_t \scirc u_t,\, R_t}
\biggl(\msmall{dJ \over dt} \scirc du_t \scirc J_\Delta \biggr),
\eqno(3.2.4)
$$
where $R_t$ is defined by relation $D_{J_t, u_t}=\dbar_{J_t, u_t}
+ R_t$, see {\sl paragraph 2.2}. Since $J_0$ and $R_0$ satisfy
the hypothesis of the {\sl Lemma 3.2.1} and $R_t$ depends
$L^p$-continuously on $J_t\in C^1$ and $u_t \in
L^{1,p}$, the solution exists in some small interval $0\le t \le
t_0$. For such solution using (2.2.1) one has
$$
\left\Vert \msmall{du_t\over dt} \right\Vert_{L^{1,p}} +
\left\Vert \msmall{dR_t\over dt} \right\Vert_{L^p}
\le C\cdot (\Vert u_t\Vert_{L^{1,p}} +
\Vert J_t  - J\st \Vert_{C^1} )
\cdot \left\Vert \msmall{dJ_t\over dt} \right\Vert_{C^1}.
$$
This implies the existence of the solution of (3.2.2) for all
$t\in [0,1]$, provided
$$
\int_{t=0}^1
\left\Vert \msmall{dJ_t\over dt} \right\Vert_{C^1}
dt \le \epsi.
$$\qed

\medskip
\noindent \sl 3.3. Perturbing cusps of pseudo-holomorphic curves. \rm

\smallskip
In this paragraph we suppose that some $p$ with $2<p<\infty$ is fixed.

\smallskip\nobreak
\state Lemma 3.3.1. {\it
For every  $\gamma >0$ and every pair of integers $\nu\ge1$, $\mu\ge1$
there exists an $\epsi=\epsi(\mu,\nu,\gamma)>0$
such that for an almost complex structure $J_0$ in $B$ with
$\Vert J_0-J\st \Vert_{C^1(B)}\le \epsi$, $J_0(0)=J\st$, and for
a $J_0$-holomorphic map $u_0: \Delta \to B(0, {1\over2})$ with $\Vert u_0
\Vert_{L^{1,p}(\Delta)} \le\gamma $ and with the~multiplicity $\mu$ at
$0\in \Delta$ the~following holds:

\smallskip
    \sli If $\nu\le2\mu-1$ then for every almost complex structure $J$
in $B$ with $J(0)=J\st$ and for every $v\in L^{1,p}(\Delta,\cc^n)$ with
$\Vert v \Vert_{L^{1,p}(\Delta)} + \Vert J-J\st \Vert_{C^1(\Delta)}
\le\epsi$ there exists $w\in L^{1,p}
(\Delta,\cc^n)$ with $w(0)=0$, satisfying
$$
\dbarj(u_0+z^\nu(v+w))=0\eqno(3.3.1)
$$
and
$$
\Vert w \Vert_{L^{1,p} (\Delta)}\le
C\cdot \left(
\Vert v \Vert_{L^{1,p}(\Delta)} +
\Vert J-J_0 \Vert_{C^1(\Delta)}
\right).\eqno(3.3.2)
$$

\smallskip
 \slii If $\nu\ge2\mu$ then for every $v\in L^{1,p}(\Delta,\cc^n)$
with $\Vert v \Vert_{L^{1,p}(\Delta)}\le\epsi$ there exists $w\in
L^{1,p}(\Delta,\cc^n)$ with $w(0)=0$, satisfying
$$
\bar\partial_{J_0}(u_0+z^\nu(v+w))=0\eqno(3.3.3)
$$
and
$$
\Vert w \Vert_{L^{1,p} (\Delta)}\le
C\cdot \Vert v \Vert_{L^{1,p}(\Delta)}.
\eqno(3.3.4)
$$
}

\smallskip
\state Proof. If $\nu\le 2\mu -1$, we fix a curve $J_t$, $t\in [0,1]$, of
$C^1$-smooth almost complex structures in $B$, which starts in $J_0$,
finishes in $J$ and satisfies the conditions $J_t(0)=J\st$
and $\int_{t=0}^1\Vert \dot J_t \Vert_{C^1(B)} dt \le 2\epsi$,
where $\dot J_t \deff dJ_t/dt$. If $\nu\ge2\mu$, we simply set
$J_t\equiv J_0$.

As in the previous lemma, we want to find a needed $w$ by solving
for $t\in [0,1]$ and $w_t\in L^{1,p}(\Delta,\cc^n)$ the equation
$$
z^{-\nu}\dbar_{J_t}(u_0 + z^\nu(t\cdot v + w_t)) =0.
$$
However, this time we need to take into consideration the fact that
now we deal with different (almost) complex structures on $B$, namely
$J\st$ and $J_t$ for any fixed $t\in [0,1]$. So we write the
last equation in the more correct form:
$$
(x + y J\st)^{-\nu}\dbar_{J_t}
(u_0 + (x + y J\st)^\nu(t\cdot v + w_t)) =0.
\eqno(3.3.5)
$$

\smallskip
The differentiation of (3.3.5) with respect to $t$ gives
$$
(x + y J\st)^{-\nu}D_{u_t, J_t}
((x + y J\st)^\nu(v + \dot w_t))
+(x + y J\st)^{-\nu}\dot J_t
\scirc du_t \scirc J_\Delta
=0,
\eqno(3.3.6)
$$
where $u_t\deff u_0+ (x + y J\st)^\nu(t\cdot v + w_t)$
and $J\st$ denotes also the pull-back of the standard
complex structure on $E\deff u_t^*TB\cong \cc^n$.

\smallskip
Firsts we show that the operator $(x + y J\st)^{-\nu} \scirc
D_{u_t, J_t} \scirc (x + y J\st)^\nu$ has the form
$\dbar_{J^{(\nu)}_t} + R^{(\nu)}_t$ for appropriate (almost)
complex structure $\dbar_{J^{(\nu)}_t}$ in $E\cong \cc^n$ and
$\rr$-linear operator $R^{(\nu)}_t$. Really, the explicit formula
(2.1.6) for $D_{u, J}$ shows that taking the standard connection
$\nabla$ in $TB\cong \cc^n$ and identifying $\Lambda^{(0,1)}\Delta
\cong \cc$ one gets
$$
D_{u_t, J_t}(v)=\msmall{1\over2}\left( \msmall{\d \over \d x}v +
J_t \msmall{\d \over \d x}v \right) + R^{(0)}_t(v)
$$
with $R^{(0)}_t\in C^0(\Delta, {\ss End}_\rr(\cc^n))$. More over,
the formula (2.2.1) implies that
$$
\Vert R^{(0)}_t \Vert_{L^p(\Delta)} \le  \Vert J_t \Vert_{C^1(B)}
\cdot \Vert du_t \Vert_{L^p(\Delta)}.
$$
Hence
$$
2(x + y J\st)^{-\nu} \scirc D_{u_t, J_t}
\scirc (x + y J\st)^\nu   =  \Bigl[\msmall{ \d \over \d x} +
(x + y J\st)^{-\nu} \scirc J_t \scirc (x + y J\st)^\nu
\msmall{\d \over \d y}  \Bigr]+
$$
$$
+ \Bigl[ \nu \cdot(x + y J\st)^{-\nu} \scirc
(1 + J_t \scirc J\st) \scirc (x + y J\st)^{\nu -1}
+ 2(x + y J\st)^{-\nu} \scirc R^{(0)}_t \scirc (x + y J\st)^\nu
\Bigr]=
$$
$$
=2\dbar_{J^{(\nu)}_t} + 2\,R^{(\nu)}_t.
$$
One has the obvious identities $(J^{(\nu)}_t)^2 \equiv -1$,
$$
\Vert J^{(\nu)}_t -J\st \Vert_{C^0(\Delta)}=
\Vert (x + y J\st)^{-\nu} \scirc (J_t -J\st)
\scirc (x + y J\st)^\nu \Vert_{C^0(\Delta)} =
\Vert J_t -J\st \Vert_{C^0(B)},
$$
and
$$
\Vert (x + y J\st)^{-\nu} \scirc R^{(0)}_t \scirc (x + y J\st)^\nu
\Vert_{L^p(\Delta)} =\Vert R^{(0)}_t \Vert_{L^p(\Delta)} .
$$
Further, $1 + J_t(0) J\st=0$ and hence $\bigl\Vert \> |z|^{-1}(1 +J_t J\st)
\bigr\Vert_{C^0(\Delta)} \le \Vert J_t -J\st \Vert_{C^1(B)}$. This gives
the estimate
$$
\Vert R^{(\nu)}_t \Vert_{L^p(\Delta)} \le \bigl(C_1 \cdot \nu +
C_2 \cdot \Vert du_t \Vert_{L^p(\Delta)} \bigr) \cdot
\Vert J_t -J\st \Vert_{C^1(B)}.
$$

\smallskip
To show the existence of the solution of (3.3.5) for all
$t\in [0,1]$ we need the estimate
$$
\left\Vert z^{-\nu}\cdot \dot J_t (u_0 + z^\nu w)
\scirc d(u_0 + z^\nu w) \right\Vert_{L^p(\Delta)}
\le C_1(\mu, \nu) \cdot \Vert u_0  \Vert_{L^{1,p}(\Delta)}^2
\cdot \Vert \dot J_t \Vert_{C^1(B)}
\eqno(3.3.7)
$$
for all sufficiently small $w\in L^{1,p}(\Delta, \cc^n)$.
The estimate trivially holds for $\nu\ge 2\mu$, since
in this case $\dot J_t\equiv 0$. Otherwise for
$\lambda \deff \min\{\mu, \nu\}$ from
{\sl Corollary 3.1.3.} one gets
$$
\left\Vert z^{-\lambda}(u_0 + z^\nu w)
\right\Vert_{L^\infty (\Delta)} +
\left\Vert z^{-\lambda +1} d(u_0 + z^\nu w)
 \right\Vert_{L^p(\Delta)}
\le C_2(\mu, \nu)\cdot
\left\Vert u_0 \right\Vert_{L^{1,p}(\Delta)}
\eqno(3.3.8)
$$
which gives the estimate (3.3.7).

\smallskip
Now consider the ordinary differential equation for $w_t \in L^{1,p}
(\Delta, B)$ with initial condition $w_0\equiv 0$ and
$$
\msmall{dw_t\over dt} = - T_{J^{(\nu)}_t,\, R^{(\nu)}_t}
\Bigl(z^{-\nu}\cdot \dot J_t \scirc du_t \scirc J_\Delta  + v \Bigr).
\eqno(3.3.9)
$$
As in {\sl Lemma 3.2.4} one has the estimate
$$
\biggl\Vert \msmall{dw_t\over dt} \biggr\Vert_{L^{1,p}(\Delta)} +
\biggl\Vert \msmall{dR_t^{(\nu)}\over dt} \biggr\Vert_{L^p(\Delta)}
\le
$$
$$
\le C\cdot (\Vert u_t\Vert_{L^{1,p}(\Delta)} +
\Vert J_t - J\st \Vert_{C^1(B)} )
\cdot \left( \bigl\Vert \msmall{dJ_t\over dt} \bigr\Vert_{C^1(B)}
+ \bigl\Vert v \bigr\Vert_{L^{1,p}(\Delta)} \right),
$$
which implies the existence of the solution of (3.3.9) for all
$t\in [0,1]$.
\par\nobreak\qed

\bigskip\noindent
\sl 3.4. Positivity of intersections.

\rm
Let us recall the notion of intersection number of two surfaces in $\rr^4$.
So, let $M_1$ and $M_2$ be two-dimensional oriented surfaces in $\rr^4$
passing through the origin, {\sl i.e.} $M_i$ is $C^1$-image of some
two-dimensional oriented surface $S_i$, $i=1, 2$. Further we suppose that both
$M_1, M_2$ intersect the unit sphere $S^3$ transversally by curves $\gamma_1$
and $\gamma_2$ respectively and that $\gamma_1$ and $\gamma_2$ does not meet
each other. Let $\widetilde M_i$ are small perturbations of $M_i$ making them
intersect transversally.

\smallskip\noindent
\bf Definition 3.4.1. \rm The {\sl intersection number of $M_1$ and $M_2$} is
defined to be the algebraic intersection number of $\widetilde M_1$ with
$\widetilde M_2$. If $M_1, M_2$ intersect only at zero and both intersect
every sphere $S^3_r$, $0<r<0 $ transversally, we shall also say that number
just defined is the \sl intersection index of $M_1$ and $M_2$ \rm at zero.
It will be denoted by $\delta_0(M_1, M_2)$ or $\delta_0$.

\smallskip\noindent\bf
Remark. \rm This number is independent of the particular choice of
perturbations $\widetilde M_1$ and $\widetilde M_2$. Later we shall use the~fact
that the~intersection number of $M_1$ and $M_2$ is equal to
the~{\sl linking number} $l(\gamma_1,\gamma_2)$ of the reducible
curves $\gamma_1$ and $\gamma_2$ on $S^3$, see [Rf].

\smallskip
\bigskip\noindent
\bf Theorem 3.4.1. \it Let $u_i:\Delta\to(\rr^4, J)$, $i=1, 2$ be two
primitive nonconstant $J$-holomorphic disks  such that $u_1(\Delta)=M_1
\not= M_2=u_2(\Delta )$ and $u_1(0)=u_2(0)$. Let $Q = M_1\cap M_2$ be
the set of intersection of those disks. If  $J$ is $C^1$-smooth, then:

\smallskip
\sli $u_i^{-1}(Q)$ is a discrete subset of $\Delta $ for $i=1, 2$.

\smallskip
\slii  The~intersection index in any such point of $Q$ is strictly positive.
More over, if $\mu_1$ and $\mu_2$ are the~multiplicities of $u_1$ and $u_2$ in
 $z_1$ and $z_2$ respectively, with $u_1(z_1)=u_2(z_2)=p$, then the~
intersection number $\delta_p$  of $M_1$ and $M_2$ at  $p$ is at
least $\mu_1\cdot \mu_2$.

\smallskip
\sliii $\delta_p=1$ iff $M_1$ and $M_2$ intersect at $p$ transversally.

\smallskip
\noindent\bf Proof. \rm
First we consider the case when $u_1:\Delta \to \rr^4$ is an imbedding. Let $N$
be the~normal bundle to $u_1(\Delta )$. Fix a~Hermitian metric $h$ on $X$ and
consider the~map
$$
f: \xi\in N_z \mapsto \exp_{u_1(z)}(\xi) \in X,
$$
where $\exp_p$ denote the~usual $h$-exponential map. Then $f$ is
a~diffeomorphism of some neighbourhood of zero-section of $N$ onto
a~neighbourhood $V$ of $u_1(\Delta )$ in $X$. Since $N$ admits a~holomorphic
structure, in the~set $V$ there exist an~{\sl integrable} complex structure
$J_1$ and $J_1$-holomorphic coordinates $w_1$ and $w_2$, such that $J=J_1$ on
$u_1(\Delta )$ and the~map $u_1$ has the~form $(z, 0)$ in coordinates
$\{w_i\}$. In particular,
$$
\vert J(w_1, w_2) - J_1(w_1, w_2)  \vert \le C\vert w_2 \vert. \eqno(3.4.1)
$$

     Consider now  $ u_2:\Delta \to V$ our second $J$-holomorphic map. We have
that $ u_2(0)=u_1(0)=(0, 0)\ddef p$ and $u_2(\Delta) \not\subset u_1(\Delta )$.
Then due to {\sl Corollary 3.1.3} in the~coordinates $\{w_1, w_2\}$ the~map
$u_2$ has the~form $(w_1(z), w_2(z))=(z^\mu a(z), z^\mu b(z))$ with $a, b\in
L^{1,p}(\Delta,\cc)$, $\bigl(a(0),\, b(0) \bigr)\not= 0$ and $\mu>0$. On
the~other hand, the~condition $(3.4.1)$ implies that the~second coordinate of
$u_2$ satisfies the~inequality $\vert \dbar w_2(z) \vert \le C\cdot \vert w_2
(z)\vert$, and consequently $w_2(z)=z^\nu b'(z)$ with $b' \in
L^{1,p}(\Delta,\cc)$, $b'(0) \not= 0$ and $\nu>0$. One can see that the
only two following cases are possible.

     {\sl Case 1.} $b(0)\not=0$, so the~vector $\bigl(a(0),\, b(0)\bigr)$ is
transversal to $u_1(\Delta )\subset \{\, w_2=0\,\}$, and $\mu=\nu$. In this
case the~intersection number of $u_1(\Delta )$ and $u_2(\Delta)$ in
the~point $p$ equals exactly $\mu$.

     {\sl Case 2.} $b(0)=0$ and $\nu > \mu$. In this case the~map $u_2$
has the~form $\bigl(z^\mu\cdot a(z),\, z^\nu\cdot b'(z)\bigr)$, so
the~intersection number is equal to $\nu$.

\smallskip
     Now one can conclude that $u_2\inv\bigl(u(\Delta )\bigr)$ is
discrete in $\Delta$ and the~intersection index of $u_1(\Delta)$ and $
u_2(\Delta)$ in any common point $p\in u_1(\Delta )\cap u_2(\Delta)$ is
positive and not less than the~multiplicity of $u_2$ in $p$. One also
concludes that, under the~hypothesis of the~theorem, there exists at most
countable set of intersection points of $M_1$ and $M_2$, and that
the~limit points of $M_1\cap M_2$ are cuspidal for both $M_1$
and $M_2$.

\medskip
     Let us consider the~situation when both curves $M_1$ and $M_2$
have a~cusp in a~point $p\in Q$ with multiplicities $\mu_1$ and $\mu_2$
respectively. Since the~problem is essentially local we can assume that $u_1$
and $u_2$ are proper maps from $\Delta$ into the~unit ball $B \subset \cc^2$,
which are holomorphic with respect to the~almost complex structure $J$. We can
also assume that $0\in\cc^2$ is a~cuspidal point for both $u_1$ and $u_2$, and
that $J(0)=J\st$.

 Let the~map $\pi_t:B\to B$ be defined by formula $\pi_t(w)\deff t\cdot
w$. The~direct calculation shows that $\Vert \pi_t^*J - J\st \Vert_{C^1(B)}
\le C\cdot t$. Thus considering $\pi_t^*J$ instead of $J$ and
$\pi_t\inv \scirc u_j(r^{\mu_j}z)$ instead $u_j(z)$, with $t$
small enough, we can assume that
$$
\Vert J\st-J \Vert_{C^1(B)} \le \epsi
\qquad \hbox{and} \qquad
\Vert u_j\Vert _{L^{1,p}(\Delta )}\le \gamma ,
$$
for any given $\epsi>0$. Further, we can also assume that $B$ contains no
other cuspidal points of $M_1$ or $M_2$, and that $\d B$ contains no
intersection points. Note, that we can also fix $r^-$ with $r > r^- > 0$, such
that there are no intersection points in the~spherical shell $B\bss B_{r^-}$.

  By {\sl Corollary 3.1.3} the~maps $u_i$ can be represented in the~form
$u_i(z)=z^{\mu_i}\cdot v_i(z)$ with $v_i\in L^{1,p}(\Delta,\cc^2)$ such
that $v_1(0) \not= 0 \not= v_2(0)$. As above, we consider two cases.

{\sl Case 1.} The~vectors $v_1(0)$ and $v_2(0)$ are not collinear. It is easy
to see that in this case $0\in\cc^2$ is an~isolated intersection point of
$u_1(\Delta)$ and $u_2(\Delta)$ with multiplicity exactly $\mu_1\cdot \mu_2$.

{\sl Case 2.} The~second case is when the~vectors $v_1(0)$ and $v_2(0)$ are
collinear. The idea is to ``turn'' the~map $u_2$ a little and to reduce
the~case to the~previous one. So let $T\in {\ss End}_\cc(\cc^2)$ be a~linear
unitary map which is close enough to identity, such that $T(v_2(0))$ is not
collinear to $v_1(0)$. Define $J^T\deff T\inv \scirc J \scirc T$ so that
$\Vert J^T -J \Vert_{C^1(B)} \le \Vert T - \id \Vert$. Applying {\sl
Lemma 3.3.1} with $u_0=z^{\mu _2}v_2$ and $v=0$ we can find $w\in L^{1,p}
(\Delta,\cc^2)$ with $w(0)=0$, such that $\dbar_{J^T} (z^{\mu_2}(v_2+w)) =0$.
The~map $\tilde u_2\deff T(z^{\mu_2}(v_2+w))$ is needed ``turned''
$J$-holomorphic map. Since such a ``turn'' can be done small enough,
the~intersection number of $u_i(\Delta)$ in the~ball $B_{r^-}$ does not
change.

 This implies, that the~intersection number of $u_1(\Delta)$ and
$u_2(\Delta)$ in any ball $B_r$ is not less than $\mu_1\cdot \mu_2$. Another
conclusion is that $0\in\cc^2$ is an isolated intersection point of
$u_j(\Delta)$. Otherwise we could find a sequence $r_i\msmall{\searrow}0$ with
at least one intersection point of $u_j(\Delta)$ in every spherical shell
$B_{r_i} \!\backslash B_{r_{i+1}}$, and therefore the~intersection number of
$u_j(\Delta)$ in the~balls $B_{r_i}$ would be strictly decreasing in $i$.

Thus the statements (i) and (ii) are proved. The proof of (iii) is now obvious
 and follows from the observation that $\mu _1 = \mu _2 =1$ in this case.\qed

\smallskip
\state Corollary 3.4.2. {\it Let $u_i:S_i\to(X, J)$, $i=1, 2$ be compact
irreducible $J$-holomorphic curves such that $u_1(S_1)=M_1\not=u_2(S_2)=M_2$.
Then
they have finitely many intersection points and the~intersection index in any
such point is strictly positive. More over, if $\mu_1$ and $\mu_2$ are the~
multiplicities of $u_1$ and $u_2$ in such a~point $p$, then the~intersection
number of $M_1$ and $M_2$ in $p$ is at least $\mu_1\cdot \mu_2$.
}

\bigskip
\noindent
\bf 4. The Bennequin index and the Genus formula for pseudo-holomorphic
curves.

\smallskip\noindent
\sl 4.1. Transversality and the Bennequin Index.

\smallskip
\rm Let $u : (\Delta , 0)\to (\cc^2, J, 0) $ be a germ of nonconstant
pseudo-holomorphic curve in zero. Without loss of generality we always suppose
that $J(0)=J\st$. Taking into account that zeros of $du$ are isolated, we
can suppose that $du$ vanishes only at zero. Further, let $w_1, w_2$ be
the standard complex coordinates in $(\cc^2, J\st)$. According to the
{\sl Lemma 3.1.3} we can write our curve in the form

\smallskip
$$
u(z) = z^\mu\cdot a + o(\vert z\vert^{\mu +\alpha}),
\qquad a\in \cc^2 \setminus \{ 0\},\> 0<\alpha<1. \eqno(4.1.1)
$$

\smallskip

For $r>0$ define $F_r\deff TS^3_r \cap J(TS^3_r)$ to be the distribution
of $J$-complex planes in the tangent bundle $TS^3_r$ to the sphere of radius
$r$. $F_r$ is trivial because $J$ is homotopic to $J\st = J(0)$. By $F$
we denote the distribution $\cup_{r>0}F_r\subset \cup_{r>0}TS^3_r \subset TB^*$,
where  $TB^*$ is the tangent bundle to the punctured ball in $\cc^2$.

\smallskip
\noindent
\bf Lemma 4.1.1. \it The (possibly reducible) curve $\gamma_r = M\cap S^3_r$ is
transversal  to $F_r$ for all sufficiently small $r>0$.

\smallskip\noindent
\bf Proof. \rm Really, due to the {\sl Lemma 3.2.3} one has $du(z) = \mu
z^{\mu-1} \cdot a\,dz + o(\vert z\vert^{\mu -1+\alpha})$. Since $J\approx J\st$
for $r$ sufficiently small, $T\gamma_r$ is close to $J\st {\bf n}_r$,
where ${\bf n}_r$ is the field of normal vectors to $S^3_r$.

On the other hand, for sufficiently small $r$ the distribution $F_r$ is
close to the one of $J\st$ complex planes in $TS^3_r$, which is orthogonal
to $J\st{\bf n}_r$.
\qed

\smallskip
This fact permits us to define a Bennequin index of $\gamma_r$. Namely, take
any nonvanishing section $\vec v$ of $F_r$ and move $\gamma_r$ along vector
field $v$ to obtain a curve $\gamma'_r$. We can make this move for a small
enough time, so that $\gamma'_r $ does not intersect $\gamma_r$. Following
Bennequin [Bn], we give

\state Definition 4.1.1. Define the \sl Bennequin index $b(\gamma_r)$ \rm to
be a linking number of $\gamma_r$ and $\gamma'_r$.

\smallskip
\rm
This number does not depend on $r>0$, taken sufficiently small, because
$\gamma_r$ is homotopic to $\gamma_{r_1}$ for $r_1<r$ within the curves
transversal to $F$, see [Bn]. It is also independent of the particular choice
of the field $\vec v$. Thus for the standard complex structure $J\st$ in
$B\subset \cc^2$ we use $\vec v\st(w_1,w_2)\deff  (-\bar w_2, \bar w_1)$
for calculating the Bennequin index of the curves
on sufficiently small spheres. For an arbitrary almost complex structure $J$
with $J(0)=J\st$ we can find the vector field $\vec v_J$, which is defined
in a small punctured neighbourhood of $0\in B$, is a small perturbation
of $\vec v\st$, and lies in the~distribution $F$ defined by $J$.

Denote by $B_{r_1, r_2}$ the spherical shell $B_{r_2}\setminus B_{r_1}$ for
$r_1<r_2$.

\smallskip
\noindent
\bf Lemma 4.1.2. \it Let $\Gamma $ be an immersed $J$-holomorphic curve in the
neighbourhood of $\overline B_{r_1, r_2}$ such that all self intersection points

of $\Gamma $  are contained in $B_{r_1, r_2}$ and for every $r_1\le r\le r_2$
all components of the curve $\gamma_r = \Gamma \cap S^3_r$ are transversal to
$F_r$. Then

\smallskip
$$
b(\gamma_{r_2}) = b(\gamma_{r_1}) + 2\cdot \sum_{x\in Sing(\Gamma) }\delta_x,
\eqno(4.1.5)
$$

\smallskip
\noindent
where the sum is taken over self-intersection points of $\Gamma $.

\smallskip
\rm
\noindent
\bf Proof. \rm Move a bit $\Gamma $ along $v_J$ to obtain
$\Gamma^\varepsilon$. By $\gamma^\varepsilon_{r_1},
\gamma^\varepsilon_{r_2}$ denote the intersections $\Gamma^\varepsilon \cap
S^3_{r_1},\Gamma^\epsi\cap S^3_{r_2}$, which are of course the moves
of $\gamma_{r_j}$ along $v_J$. We have that $l(\gamma_{r_2},\gamma
_{r_1}^\epsi) - l(\gamma_{r_1},\gamma_{r_1}^ \epsi) =
{\ss int}(\Gamma ,\Gamma^\epsi)$, here $l(\cdot ,\cdot )$ is the~linking
number and ${\ss int}(\cdot ,\cdot )$ is the~intersection number.

Now let us calculate ${\ss int}(\Gamma ,\Gamma^\epsi)$. From {\sl Theorem
3.4.1} we now that there are only finite number $\{ p_1,\ldots, p_N\} $ of
self-intersection points of $\Gamma $. Take one of them, say $p_1$. Let $M_1,
\ldots, M_d$ be the disks on $\Gamma $ with common point $p_1$ and otherwise
mutually disjoint. More precisely we take $M_j$ to be an irreducible
components of $\Gamma \cap B_{p_1}(\rho)$ for $\rho >0$ small enough.
Remark that $M_j$ are transversal to $v_J$, so their moves $M_j^\varepsilon$
don't intersect them, {\sl i.e.} $M_j\cap M_j^\epsi = \emptyset $.
Note also that ${\ss int}(M_k, M_j) = {\ss int}(M_k, M^\epsi_j)$
for $\varepsilon >0$ sufficiently small. So ${\ss int}(M_k, M_j) =
{\ss int}(M_k, M^\epsi_j) + {\ss int}(M^\epsi_k, M_j)$.
Finally $\delta_{p_1} = \sum_{1\le k<j\le d} {\ss int}(M_k, M_j) =
{\ss int}(\Gamma \cap B_{p_1}(\rho ),\Gamma^\epsi\cap
B_{p_1}(\rho ))$. This means that ${\ss int}(\Gamma ,\Gamma^\epsi)
= 2\cdot \sum_{j=1}^N\delta_{p_j}$.
\qed

\medskip\noindent
\sl 4.2. Proof of the Genus formula and estimate of the Bennequin index from
below.

\nobreak\rm
In $\S3$ we have proved that compact $J$-holomorphic curve with finite number
of irreducible components $M = \bigcup_{i=1}^d M_i$ has only finite number of
nodes ({\sl i.e.} self-intersection points) points, provided $J$ is of class
$C^1$.

For each such point $p$ we can introduce following the {\sl Definition 3.4.1}
\rm the~self-in\-ter\-sec\-tion
number $\delta_p(M)$ of $M$ at $p$. Namely, let
$S_j$ be a parameter curve for $M_j$, {\sl i.e.} $M_j$ is given as an image of
the $J$-holomorphic map $u_j : S_j\to M_j$. We always suppose that the
parametrization $u_j$ is primitive, {\sl i.e.} they cannot be decomposed like
$u_j=v_j\circ r$ where $r$ is nontrivial covering of $S_j$ by another
Riemannian surface. Denote by $\{ x_1,\ldots, x_N\}$ the set of all preimages
of $p$ under $u : \bigsqcup_{i=1}^d S_j\to X$, and take mutually disjoint
disks $\{ D_1,\ldots, D_N\}$ with centers $x_1,\ldots, x_N$ such that their
images have no other common points different from $p$. For each pair $D_i,
D_j$, $i\not= j$, define an intersection number as in {\sl Definition 3.4.1}
and take the sum over all different pairs to obtain $\delta_p(M)$.

Put now $\delta = \sum_{p\in D(M)}\delta_p(M)$, here the sum is taken over the
set $D(M)$ of all nodes of $M$, {\sl i.e.} points which have at least two
preimages.

Consider now the set $\{ p_1,\ldots, p_L\} \subset \bigcup_{j=1}^d S_j = S$
of all cusps of $M$, {\sl i.e.} points where the differential of the
appropriate parametrization vanishes. Take a small ball $B_r$ around say
$u(p_i)$ and small disk $D_{p_i}$ around $p_i$.
Let $\gamma^i_r\deff u(\Delta_{p_i}) \cap \partial B_r$ and $b_i$ be
the~Bennequin index of $\gamma^i_r$, defined in {\sl Definition 4.1.1}.

\state Definition 4.2.1. The~number $\varkappa_i\deff (b_i +1)/2$
is called the~{\sl conductor of the~cusp $a_i=u(p_i)$}.

\smallskip Put also $\varkappa\deff \sum_{i=1}^L \varkappa_i$. We are able
to state now the Genus formula.

\medskip
\noindent
\bf Theorem 3. \it Let $M =\cup_{j=1}^d M_j$ be a compact $J$-holomorphic
curve in almost complex surface $(X, J)$ with the~distinct irreducible
components $\{M_j\}$, where $J$ is of class $L^{1,p}$. Then

\smallskip
$$
\sum_{j=1}^d g_j = {[M]^2 - c_1(X, J)[M]\over2} + d -\delta - \varkappa ,
\eqno(4.2.1)
$$

\smallskip
\noindent
where $g_j$ are the genera of parameter curves $S_j$.

\smallskip\noindent\bf Proof.
\rm The main line of the proof is the reduction of a general case to
the case when $M$ is immersed, which is already proved in $\S 1$.

Let $u: \bigsqcup_{j=1}^d S_j \longrightarrow X$ be a $J$-holomorphic map,
which parameterizes the curve $M$. Let also $\{x_1,\ldots,x_n\}$ be the
set of cusp-points of $M$, {\sl i.e.} the images of critical points of $u$.
Our reduction procedure is local, we make our constructions in a
neighbourhood of every point $x_j$ separately. So from now on we fix such
a point $x$. Due to the {\sl Corollaries 2.2.3} and {\sl 3.4.2} there exists
a neighbourhood $U$ of $x$ which contains no other cusp-points and no other
self-intersection points. The {\sl Theorem 3.4.1 i)}
implies that taking the neighbourhood $U$ small
enough we may assume that any component of $M\cap U$ goes through $p\in U$.
Without losing generality we may also assume that $U$ is the unit ball $B
\subset \cc^2$, that $x$ corresponds to the center $0$ of $B$, and $J(0)=J\st$.

Denote by $\Gamma_j$ the irreducible components of $M \cap B$ and let
$u_j :\Delta \to B$ be parameterization of $\Gamma_j$, such that $u_j(0)=0
\in B$. Denote $\mu_j \deff \ord_0 du_j +1$, so $0\in \Delta$ is the cusp-point
for $u_j$ \iff $\mu_j\ge2$.

\smallskip\noindent
{\sl Step 1. Rescaling procedure}.

\nobreak
Take some cusp-component $\Gamma_j$. Due to {\sl Corollary  3.1.3} the map
$u_j$ has the form $u_j(z)= z^{\mu_j} \cdot a_j + O\bigl(|z|^{\mu_j+\alpha}
\bigr)$ with the constant $a_j\not=0\in\cc^2$. More over, {\sl Lemma 3.2.3}
provides that
$$
\bigl\vert du_j(z) - \mu_j z^{\mu_j-1} a_j \cdot dz \bigr\vert
\le C(\Vert u_j \Vert_{L^{1,p}(\Delta)}, \alpha, p) \cdot
|z|^{\mu_j-1+\alpha},
$$
for any $0<\alpha<1$ and $2<p<\infty$.

For $0<r\le1$ we consider the maps $\pi_r:B\to B$, $\pi_r(w)\deff r^{\mu_j}
\cdot w$, the rescaled maps $u^{(r)}_j:\Delta \to B$, $\pi_r\scirc
u^{(r)}_j(z) = u_j (rz)$, and the rescaled almost complex structures $J^{(r)}
\deff \pi_r^*J$ in $B$. The map $u^{(r)}_j$ is a parametrization of
$J^{(r)}$-holomorphic curves $\pi_r^{-1}\Gamma_j$. One can see that
$\Vert u^{(r)}_j(z) - z^{\mu_j} a_j \Vert_{C^0(\Delta)} \le C\cdot
r^{\mu_j+\alpha}$ and $\Vert du^{(r)}_j(z) - \mu_j\cdot z^{\mu_j-1} a_j \cdot
dz \Vert_{C^0(\Delta)} \le C\cdot r^{\mu_j-1+\alpha}$. In particular, there
exists $r_j>0$ such that for all $0<r<r_j$ the maps $u^{(r)}$ are transversal
to all spheres $S^3_1$. On the other hand, $\Vert J^{(r)} - J\st \Vert_{C^1(B)}
= O(r^{\mu_j})$, and for any $\epsi>0$ we can choose sufficiently small
$r_j>0$ such that $\Vert J^{(r)} - J\st \Vert_{C^1(B)} \le \epsi$ for all
$0<r<r_j$. So replacing $B$, $J$ and all the maps $u_j$ by their rescaled
counterpart we may assume that the almost complex structure $J$ satisfies
the estimate
$$
\Vert J - J\st \Vert_{C^1(B)} \le \epsi.
\eqno(4.2.2)
$$

\medskip\noindent
{\sl Step 2. Reduction to the case of holomorphic cusp-points}.

\nobreak
Recall, that in the {\sl Lemma 1.1.1} we establish the natural diffeomorphism
between the space $\jj$ of orthogonal complex structures in $\rr^4$ and
the unit sphere $S^2\deff \{\,(c_1, c_2, s)\in \rr^3 : c_1^2 + c_2^2
+ s^2 =1\, \}$. Define the function $\Phi :\jj \to \rr^2$ by setting $\Phi(J)
=(c_1, c_2)$. The map $\Phi$ defines a diffeomorphism between the upper
half-sphere in $S^2$ and the unit disk $\rr^2$, such that the north pole
$(0,0,1)\in S^2$ corresponds to the center of the disk.

Define the function $f:B \to \rr^2$ setting $f(w)\deff \Phi(J(w))$. Due to
$(4.2.3)$, $f$ parameterizes the given almost complex structure $J$,
$\Vert f \Vert_{C^1(B)} \le C\epsi$, and $f(0)=0$. Fix a~cut-off function
$\chi$ in $B$, such that $0\le\chi\le1$, $\chi\bigm|_{B(1/2)}\equiv 1$,
$\supp \chi \comp B$, and $\Vert d\chi \Vert_{C^1(B)} \le 3$. For $0\le t \le1$
and $0<\sigma<1$ set $f_{\sigma, t}(w) \deff (1-t\chi(w/\sigma))\cdot f(w)$
and define $J_{\sigma, t}\deff \Phi^{-1}(f_{\sigma, t})$. One can easily see
that $\Vert f_{\sigma, t} \Vert_{C^1(B)}  \le 4\cdot\Vert f \Vert_{C^1(B)}$,
and consequently $\Vert J_{\sigma, t} - J\st \Vert_{C^1(B)} \le C_1\cdot
\Vert J - J\st \Vert_{C^1(B)}\le C_1\cdot \epsi$. Here the constant $C_1$,
as also the constants $C_2,\ldots, C_5$ below in the proof, are
independent of $\epsi$, $t$ and $\sigma$. For fixed $\sigma$, the curve
$J_{\sigma, t}$, $0\le t\le1$ is a homotopy between $J\equiv J_{\sigma, 0}$
and an almost complex structure $J_\sigma \deff J_{\sigma, 1}$, such that
$J_{\sigma, t} \equiv J$ in $B\backslash B(\sigma)$ and $J_\sigma \equiv J\st$
in $B({\sigma\over2})$. More over, we have
$$
\Vert J_{\sigma, t} - J\st \Vert_{C^1(B)} +
\Vert \msmall{ dJ_{\sigma, t} \over dt} \Vert_{C^1(B)} \le C_2\cdot \epsi.
\eqno(4.2.3)
$$

\medskip
Now we fix some $p>p'>2$ and set $\alpha\deff {2\over p'} -{2\over p}$.
Due to (4.2.3) we can apply {\sl Lemma 3.3.1} \sli to the map $u_j$ with
$\nu_j=\mu_j$, $v_j\equiv0$, and with the curve of $C^1$-smooth almost
complex structures $J_{\sigma, t}$. As result we obtain a curve of maps
$u_{j, \sigma, t} = u_j + z^{\mu_j}w_{j, \sigma, t}$ with $w_{j, \sigma, t}(0)
=0$ and $\Vert w_{j, \sigma, t} \Vert_{L^{1,p}(\Delta)} \le C_3\cdot \Vert u_j
\Vert_{L^{1,p}(B)}$. The condition $\supp \dot J_{\sigma, t} \subset B(\sigma)$
together with {\sl Corollary 3.1.3} provide that $\supp (\dot J_{\sigma, t}
\scirc du_{j, \sigma, t}) \subset \Delta(\rho)$ with $\rho\sim\sigma^{1/\mu_j}$.
 Due to (3.3.8) and the H\"older inequality we obtain
$$
\left\Vert z^{-\mu_j}\cdot \dot J_{\sigma, t} (u_{j, \sigma, t})
\scirc du_{j, \sigma, t} \right\Vert_{L^{p'}_{\vphantom{1}}(\Delta)}
\le C_4 \cdot \sigma^{\alpha/\mu_j}\cdot \Vert u_j  \Vert_{L^{1,p}(\Delta)}^2
\cdot \Vert \dot J_{\sigma, t} \Vert_{C^1(B)},
\eqno(4.2.4)
$$
and consequently
\smallskip\noindent
$$
\Vert w_{j, \sigma, 1} \Vert_{L^{1,p'}_{\vphantom{1}}(\Delta)}
\le C_5 \cdot \sigma^{\alpha/\mu_j}\cdot \Vert u_j  \Vert_{L^{1,p}(\Delta)}^2
\eqno(4.2.5)
$$
The {\sl Lemma 3.2.3} provides that for $\sigma$ small enough the
$J_\sigma$-holomorphic maps $u_{j, \sigma} \deff u_{j, \sigma, 1}$ are
transversal to all spheres $S^3_r$, $0<r<{1\over2}$, have no self-intersection
points in $B({1\over2}) \backslash B({1\over4})$, and the Bennequin index
of $u_{j, \sigma}(\Delta)\cap S^3_r$, ${1\over4}<r<{1\over2}$, coincides with
the one of $\Gamma_j\cap S^3_r$, $0<r<1$. More over, we may match
$u_{j, \sigma}$ to the rest of $M$, changing in an appropriate way the almost
complex structure $J$, see {\sl Lemma 6.2.1}.

\smallskip
As result we conclude that for appropriate small enough neighbourhoods $U_1
\comp U$ of any cusp-point $p\in M$ there exist a perturbed almost complex
structure $\tilde J$
and $\tilde J$-holomorphic curve $\widetilde M$, which is parameterized by
$\tilde u:
\bigsqcup_{j=1}^d S_j \longrightarrow X$ and has the following properties:

\sli $\tilde J$ and $\widetilde M$ coincide on $X\backslash U$ with $J$ and $M$
correspondingly;

\slii $U_1$ is a ball centered in $p$ and $\widetilde M \cap U_1$ is obtained
from
$M\cap U_1$ by perturbing of its components in the above described way;

\sliii the (possibly reducible) curve $\tilde \gamma \deff \d U_1\cap \widetilde
 M$
is isotopic to $\gamma \deff \d U_1\cap M$, in particular all the corresponding
components $\tilde \gamma_j$ of $\tilde \gamma$ and $\gamma_j$ of $\gamma$ have
the same Bennequin index, and the linking number $l(\tilde \gamma_i,
\tilde \gamma_j)$ equals to $l(\gamma_i,\gamma_j)$;

\sliv $\tilde u$ is homotopic to $u$;

\slv the cusp-points of $\widetilde M$ coincide with the ones of $M$ and
$\tilde\delta + \tilde\varkappa= \delta +\varkappa$,

\slvi $\tilde J$ is integrable in a neighbourhood.

The last equality of  \slv follows from \sliii du to the {\sl Lemma 4.1.2}.
Thus the formulas (4.2.1) for $M$ and $\widetilde M$ are equivalent.

\medskip\noindent
{\sl Step 3. Final reduction to the case of an~immersed curve}.

\nobreak
This step is rather obvious and uses the following fact, shown in [Bn]. Let
$B$
be the unit ball in $\cc^2$ and let $\Gamma_0$ be an irreducible holomorphic
curve in $B$, which is transversal to $\d B$ and is defined as zero-divisor
of a holomorphic function $f$. Then for any sufficiently small nonzero
$\epsi\in\cc$ the curve $\Gamma_\epsi$, defined as the zero-divisor
of the function $f+\epsi$, are smooth and of the same
genus $g$. More over, all $\Gamma_\epsi$ are transversal to $\d B$ and the
Bennequin index of $\gamma_\epsi \deff \Gamma_\epsi \cap \d B$ equals to
$2g-1$. In particular, the conductor of a single cusp-point of a holomorphic
curve in $B$ can be defined as a genus of general small perturbation  to a
smooth curve.

In general, let $M$ be a $J$-holomorphic curve in $X$, such that $J$ is
integrable in a neighbourhood of every cusp-point of $M$. One can now
see that we can perturb $J$ and $M$ to an almost complex structure
$\tilde J$ and a $\tilde J$-holomorphic curve $\widetilde M$, satisfying
the conditions {\sl i)}--{\sl iv)} from above, and the desired condition

{\sl v$'$)} $\widetilde M$ has no cusp-points, $\sum\tilde g_j =\sum g_j +
\varkappa$, $\tilde \delta =\delta$.
\qed

\medskip
An important corollary of the proof of {\sl Theorem 3} is the estimate
from below of the conductor number of a cusp point.

\state Corollary 4.2.2. \it Let $J$ be an almost complex structure in
the unit ball
$B\subset \cc^2$ with $J(0)=J\st$, and let $u:\Delta \to B$ be a $J$-holomorphic
 curve with the cusp-point $u(0)=0$. Then $\varkappa_0$ is integer,
$\varkappa_0 \ge \ord_0du$, or equivalently, for all  sufficiently small $r>0$
the Bennequin index of $\gamma_r\deff u(\Delta) \cap S^3_r$ is odd and
satisfies the inequality
$$
b(\gamma_r)\ge 2\cdot \ord_0du - 1. \eqno(4.2.6)
$$

\medskip\noindent
\bf Proof. \rm Rescaling $u$ as in the {\sl Step 1} of the proof of the {\sl
Theorem 3} we may assume that $\Gamma\deff u(\Delta)$  has no nodes and cusps,
excepting $0\in \Delta$ , and is transversal to all spheres $S^3_r$, $0<r<1$,
so that the Bennequin index $b(\gamma_r)$ is the same and equals to
$2\varkappa_0 -1$. We may also assume that for the rescaled almost complex
structure $J$ the estimate $\Vert J - J\st \Vert_{C^1(B)} \le \epsi$ with
the appropriate $\epsi$ is fulfilled. Applying the {\sl Step 2} with
a sufficiently small $\sigma$ we can deform $u$ into a $\tilde J$ holomorphic
map $\tilde u$, which has the following properties:

\smallskip
\sli $\widetilde \Gamma \deff \tilde u(\Delta)$ is transversal to all spheres
$S^3_r$, $0<r<1$;

\smallskip
\slii $\widetilde \Gamma$ has no self-intersection points in $B\backslash
B({1\over2})$, and the Bennequin index of $\tilde\gamma_r \deff \widetilde
\Gamma\cap S^3_r$, ${1\over2}<r<1$, coincides with the one of $\gamma_r$;

\smallskip
\sliii $\tilde J$ coincides with $J$ in $B\backslash B({1\over2})$,
is integrable in the neighbourhood of $0\in B$, and $\ord_0 d\tilde u =
\ord_0 du $.

\smallskip
Now of the corollary follows from {\sl Lemmas 4.1.} and the fact that for
{\sl integrable} complex structures the same statement is true, see [Bn].
\qed

\medskip
Another modification of the proof of {\sl Theorem 3} leads to the following
corollary.

\state Lemma 4.2.3 \it Let $J$ be a $C^1(X)$-almost complex
structure on $X$ and $M\subset X$ a compact $J$ holomorphic curve
parameterized by $u: S \to X$. Then

\sli $u$ can be $L^{2,p}$-approximate by $J_n$-holomorphic immersions
$u_n: S \to X$ with $J_n \longrightarrow J$ in $C^1(X)$.

\slii there exists a $C^1(X)$-approximation $J_n$ of $J$ and a sequence of
$J_n$-holomorphic imbedded curves $M_n$ which converges to $M$ in the Gromov
topology.
\rm

\state Proof. We do not need this result for the purpose of this paper, so we
give only the sketch of the proof.

On the first step, one applies the rescaling procedures in order to find
appropriate small neighbourhoods of the cusp-points of $M$.

On the second step, one applies the {\sl Lemma 3.1.3} to the chosen
neighbourhoods, taking $\nu=1$, $v$ sufficiently small,
and $J$ unchanged. After matching procedure, this gives the statement {\sl i)}.

To obtain the statement \slii, one needs firsts to deform all the nodes
of $M$ into simple transversal ones, and to find an appropriate small
neighbourhood $U$ of every node. In some complex coordinates $(w_1, w_2)$
in $U$, the curve $M$ is defined by equation $w_1\cdot w_2=0$. It remains
to replace a node $M\cap U$ by a ``small handle'' $M_\epsi\deff \bigl\{\,
(w_1, w_2) \in U \> : \> w_1\cdot w_2=\epsi\,\}$ with $\epsi$ sufficiently
small and to use the matching procedure once more.
\qed

\bigskip
\noindent
\bf 5. Continuity principles and envelopes.

\smallskip\noindent
\sl 5.1. Continuity principles relative to the K\"ahler spaces.

\rm
Our aim in this paragraph is to prove {\sl Theorem 1} and give some
corollaries from it and from other results of this paper. First of all we
shall need an appropriate form of so called "continuity principle" for the
extension of meromorphic mappings.

We shall work in Gromov topology on the space of stable curves over $X$, see
Appendix 6.3 for definitions. Continuity of the family of curves will
be understand allwayse with respect to this topology. The reason for this
refinement of the usual cycle topology is the following Proposition, which
will be used in the poof of Continuity principle.

\def\calu{{\cal U}}

Consider a sequence $(C_n,u_n)$ of stable curves over a complex manifold $X$,
which converges in Gromov topology to $(C_0,u_0)$. We suppose $C_n$ to be
smooth, exept $C_0$.
\smallskip\noindent\bf
Proposition 5.1.1. \it There exist a natural $N$ and smooth complexe surface
$\calc $ together with a surjective holomorphic mapping $\pi :\calc \to
\Delta $ and holomorphic mapping $\calu :\calc \to X$ such that the  family
$$
\matrix{\calc & {\buildrel \calu  \over \longrightarrow } & X\cr
\llap{$\pi$}\downarrow & &\cr
\Delta \cr
}
$$
\smallskip\noindent
is a holomorphic family of stable curves over $X$ joining $(C_N,u_N)$ with
$(C_0,u_0)$. More precisely:
\smallskip\noindent
1) For every $\lambda\in \Delta $ $C_{\lambda }=\pi^{-1}(\lambda )$ is a connected
nodal curve with boundary $\d C_{\lambda }$ and the pair $(C_{\lambda }, u_
{\lambda }:=\calu\mid_{C_{\lambda }} )$ is a stable curve over $X$.

\noindent
2) For $\lambda $ outside of zero $C_{\lambda }$ is connected and smooth.

\noindent
3) $(C_0,u_0)$ is our limit and there is $\lambda_N\in \Delta $ such that
$(C_{\lambda_N},u_{\lambda_N})=(C_N,u_N)$.

\noindent
4)  There are an open subsets $U_1,...,U_m$ of $\calc $ such that
$U_j$ is biholomorphic to  $\Delta \times A_j$, where $A_j$ are an annulai
in $\cc $. In particular the following diagram is commutative

$$
\matrix{U_j&\sim &\Delta \times A_j\cr
\llap{$\pi$}\downarrow & &\downarrow\rlap{$\pi_{\Delta }$}\cr
\Delta & = & \Delta \cr}
$$
\smallskip\noindent\rm
The proof will be postponed till Appendix 6.3.

For the notion of meromorphic mapping from a domain $U$ in complex manifold
into a complex manifold (or space) $Y$ we refer to [Rm]. We only point
out here that meromorphic mappings into $Y = \cc\pp^1$ are exactly meromorphic
functions on $U$, see [Rm].

\smallskip\rm
\smallskip\noindent\bf
Definition 5.1.1. \it A Hermitian complex space $Y$ is called disk-convex if
for any sequence
$(C_n,u_n)$ of stable  curves over $Y$ such that

1) ${\sl area}(u_n(C_n))$ are uniformly bounded;

2) $u_n$ $C^1$-converges in the neighborhoods of $\d C_n $;

\noindent \it there is a compact $K\subset Y$ which contain all $u_n(C_n)$.
\smallskip\rm
This definition obviously transforms to the case when $Y$ is a symplectic
manifold.
In this case one should concider $(C_n,u_n)$ as $J_n$-holomorphic curves,
with $J_n$
converging to some $J$ (everything in $C^1$-tology), and all structures
being tamed
by a given symplectic form.

Let $U$ is now a domain in complex manifold $X$ and $Y$ is a complex space .
\smallskip\noindent\bf
Definition 5.1.2. \it An envelope of meromorphy of $U$ relative to $Y$ is a
largest
domain $(\hat U_Y,\pi )$ over $X$, which contains $U$ (i.e. there exists an
imbedding
$i:U\to \hat U_Y$ with $\pi \circ i={\sl Id}$), such that every meromorphic
map
$f:U\to Y$ extends to a meromorphic map $\hat f:\hat U_Y\to Y$.
\smallskip\rm
\smallskip
Using Cartan-Thullen construction for the germs of meromorphic mappings of
open subsets of $X$ into $Y$ instead of germs of holomorphic functions, one
can prove the existence and uniqueness of the envelope:

\smallskip\noindent
\bf Proposition 5.1.2. \it For any domain $U$ in complex space $X$ and for any
complex space $Y$ there exists a maximal domain $(\hat U_Y,\pi )$ over $X$,
containing $U$ such that every meromorphic mapping $f:U\longrightarrow Y$
extends to a meromorphic mapping $\hat f:\hat U_Y\longrightarrow Y$. Such
domain is unique up to a natural isomorphism.

\smallskip
\noindent \rm See [Iv-1] for details.

\smallskip
\smallskip\noindent\bf
Theorem 5.1.3. \tensl (Continuity principle-I). \it Let $X$ be a
disk-convex complex surface and $Y$ a disk-convex K\"ahler space. Then the
envelope of meromorophy $(\hat U_Y,\pi )$ of $U$
relative to $Y$ is also disk-convex.
\smallskip\rm
This result can be reformulated in more usual terms as follows.

Let $\{ (C_t,u_t)\}_{t\in [0, 1]}$
be a family of complex curves over $X$ with boundaries, parameterized by
unit interval. More precisely, for each $t\in [0, 1[$ a smooth Riemannian
surface with
 boundary $(C_t, \d C_t)$ is given together with holomorphic
mapping $u_t: C_t\longrightarrow X$, which is $C^1$-smooth up to the boundary.
Note that $C_1$ is not supposed to be smooth, i.e. it can be a nodal curve!
(see Appendix 6.3).

Suppose that
 in the neighbourhood $V$ of $u_0(C_0)$ a meromorphic map
$f$ into complex space $Y$ is given.

\smallskip
\noindent
\bf Definition 5.1.3. \rm We shall say that the map $f$ meromorphically
extends along the family $(C_t,u_t)$ if for every $t\in [0, 1]$ a
neighbourhood
$V_t$ of $u_t(C_t)$ is given, and given a meromorphic map $f_t:V_t
\longrightarrow
Y$ such that

a) $V_0 = V$ and $f_0= f$;

b) if $V_{t_1}\cap V_{t_2}\not=\emptyset $ then $f_{t_1}\ogran_{V_{t_1}\cap
V_{t_2}} = f_{t_2}\ogran _{V_{t_1}\cap V_{t_2}}$.

\smallskip\noindent

\smallskip
\noindent

\smallskip
\noindent
\bf Theorem 5.1.4. (\sl Continuity principle-II) \it Let $U$ be a domain in
complex surface $X$.
Let $\{ (C_t,u_t)\} _{t\in [0, 1]} $  be a continuous family of complex
curves over $X$ with boundaries in $U_1$ - relatively compact subdomain in $U$.
Suppose moreover that $u_0(C_0)\subset U$ and that $C_t$ for $t\in [0,1[$
are smooth. Then every meromorphic mapping
$f$ from $U$ to the disk-convex K\"ahler space $Y$ extends
meromorphically along the family $(C_t,u_t)$.

\smallskip\noindent\rm

Taking as image manifold $Y$ complex line $\cc $ or $\cc \pp ^1$ we get the
continuity principles for holomorphic or meromorphic functions. When it is
necessary to underline that we consider the mappings into the certain
manifold $Y$ we shall refer to the statement above as to the \sl continuity
principle relative to $Y$ \rm or \sl c.p.\ for the meromorphic mappings into
$Y$.

\rm

The discussion made above leads to the following

\smallskip
\noindent
\bf Corollary 5.1.5.
\it If we have the domain $U$ in complex surface $X$, K\"ahler space $Y$ and
family $\{ (C_t,u_t)\} $ satisfying the conditions of the "continuity
principle", then the family $\{ (C_t,u_t)\} $ can be lifted onto the $\hat
U_Y$, {\sl i.e.} there exists a continuous family $\{ (C_t,\hat u_t) \} $ of
analytic curves in $\hat U$ such that $\pi \circ \hat u_t = u_t$ for
each $t$.
\smallskip\rm Of course the point here is that the mapping can be extended
to the neighborhood of $u_1(C_1)$, which is a reducible curve having in
general compact components.

\smallskip\noindent\sl
Proof of the Theorem 5.1.3. \rm Suppose that there is a subsequnce of our sequence, we still denote
it as $(C_n,u_n)$, which is not contained in any compact subset of the envelope
$(\hat U_Y, \pi )$. Put $v_n=\pi \circ u_n$ and consider our sequence as a sequence
$(C_n,v_n)$ of stable curves over $X$. While the areas are bounded, we can suppose
that this sequence converges in Gromov topology by the disk-convexity of $X$. Denote
by $(C_0,v_0)$ its limit. For $N$ sufficiently large take a holomorphic
family $(\calc ,\pi ,\calv )$ joining $(C_N,v_N)$ with $(C_0,v_0)$ as in
Proposition 5.1.1.

\smallskip
 Let $f:\hat U_Y \to Y$ be some meromorphic mapping of our envelope into a disk-convex
K\"ahler space $Y$. Composing $f$ with $\calv $ we get a meromorphic map $h:=
f\circ \calv :
\bigcup U_j\cup \pi^{-1}(V(\lambda_N))\to Y$. Here $V(\lambda_N)$ is suffuciently
small neighborhood of $\lambda_N \in \Delta $.

Denote by $W$ the maximal connected open subset of $\Delta $ containing
$V(\lambda_N)$
such that $h$ meromorphically
extends onto $\bigcup U_j\cup \pi^{-1}(W)$. We whant to prove that $W=\Delta
 $. Denote
by $\Gamma_{f_{\lambda }}$ the graph of the restriction $f\mid_{C_{\lambda }}$. Fix
some Hermitian mertic on $\calc $.

\smallskip\noindent\bf
Lemma 5.1.6. \it For any compact $K\subset \Delta $ there is a constant
$M_K$, such that ${\sl area}(\Gamma_{f_{\lambda }})\le M_K$ for all $\lambda
\in W\cap K$.
\smallskip\noindent\sl
Proof. \rm Shrinking $\calc $ if nessessary, we can suppose that $f$ has only finite
number of points of indeterminancy in $\bigcup_jU_j$. Denote by $S_W$ the discrete
subset in $W$, which consists of points $s\in W$ such that either the fiber
$C_s$ is singular or  $s$ is the projections of the indeterminancy points of
$f\mid_{\calc \mid_W}$.
Fix a point $\lambda_1\in W\setminus S_W$. Take a path $\gamma :[0,1]\to
W\setminus
S_W$ connecting
$\lambda_1$ with some $\lambda_2\in W\setminus S_W$. Take some relatively compact
in $W\setminus S_W$ neighborhood $V$ of $\gamma ([0,1])$.

Recall that a K\"ahler metric on complex space $Y$ consists of locally finite
covering $\{ V_{\alpha }\} $ of $Y$ and strictly plurisubharmonic functions
$\phi_{\alpha }$ on $V_{\alpha }$, such that $\phi_{\alpha }-\phi_{\beta }$ are
pluriharmonic on $V_{\alpha }\cap V_{\beta }$. $\{ f^{-1}(V_{\alpha })\} $ is a
covering of $\pi^{-1}(V)$ and $dd^cf^*\phi_{\alpha }=w$ is correctly defined
semi-positive closed forme on $\pi^{-1}(V)$.
 We have that
$$
{\sl area}(\Gamma_{f_{\lambda }})={\sl area}(C_{\lambda })+
\int_{C_{\lambda }}w.
$$
\noindent
By the Stokes formula
$$
\int_{C_{\gamma (0)}}w-
\int_{C_{\gamma (1)}}w =
\int_{\cup_tC_{\gamma (t)}}dw -
\int_{\cup_t{\d C_{\gamma (t)}}}w=
$$
$$
=
\int_{\cup_t{\d C_{\gamma (t)}}}w.
$$
\noindent
This is obviously uniformly bounded on $\lambda_2\in K\setminus S_W$. Which
implies the uniform bound in $K$ by Gromov (or Bishop in this case) compactness
theorem.
\smallskip
\hfill{q.e.d.}

\smallskip\noindent\bf
Lemma 5.1.7. \it Mapping $h$ meromorphically extends onto
$\calc $.
\smallskip\noindent\sl
Proof. \rm Let us prove first that $h$ meromorhpically extends onto $\calc $
minus the singular fiber $C_0$, i.e. that $W\supset \Delta \setminus \{ 0\} $.

The same arguments as in [Iv-2] (the only difference is that in [Iv-2] curves $C_
{\lambda }$ where disks) show that $(\Delta \cap \d W)\setminus (\{ 0\})
$ is an
analytic variety. Now, using Thullen-type extension theorem of Siu, see [Si-1],
we can extend $h$ onto $\calc \setminus C_0$.

Denote by $\hat C_0$ the union of compact components of $C_0$.
$\hat C_0$ contracts to a finite number of normal points. To prove this we
can suppose that $\hat C_0$ is connected. Othervice we can apply the same
arguments as below to the connected components. We shall prove now that
$\hat C_0$ contracts to a normal point. Denote by $E_1,...,E_n$ the
irreducible components of $\hat C_0$. All we must to prove is that the
matrix $(E_i,E_j)$ is negatively defined, see [Gra]. Denote by $l_i$ the
multiplicity of $E_i$. Thus $\hat C_0=\Sigma_{i=1}^nl_i\cdot E_i$. Denote
by $M$ the $\zz $-module generated by $E_1,...,E_n$ with scalar product -
intersection of divisors. Put $D=\hat C_0$. Then we have

\noindent
1) $E_i\cdot E_j\ge 0$ for $i\not= j$;

\noindent
2) $D\cdot E_i\le 0$ for all $i$, because this number is not more then the
intersection of $E_i$ with nonsingular fiber $C_{\lambda }$.
\smallskip
By Proposition 3 from [Sh], v.1, Appendix 1, we have that $A\cdot A\le 0$
for all $A\in M$ and $A\cdot A=0$ iff $A$ is proportional to $D$. But
$D\cdot D <\hat C_0\cdot C_{\lambda }=0$, where $C_{\lambda }$ is a smooth
fiber. This proves that $(E_i\cdot E_j)$ is negatively defined.

So all that left to prove,
 is that normal point is a removable singularity for the meromorphic mappings
into a disk-convex K\"ahler space.

Let $(\calc ,s)$ is a germ of two-dimensional variety with isolated normal
singularity $s$. Let a meromorphic map $h:\calc \setminus \{ s\} \to Y$ is given.
Realise $\calc $ as a finite proper analytic cover over a bidisk $\Delta^2$ with
$s$ being the only point over zero. Denote by $p:\calc \to \Delta^2$ this
covering.

Composition $h\circ p^{-1}$ is a multivalued meromorphic map from $\Delta^2
\setminus \{ 0\} $
to $Y$. This can be extended to origin, see [Si-1]. Consider the following
analytic
set in $\calc\times Y$:
$$
\Gamma = \{ (x,y)\in \calc \times Y: (p(x),y)\in \Gamma_{h\circ p^{-1}} \} .
$$
The irreducible component of $\Gamma $ which containes $\Gamma_h$ will be the
graph of the extension of $h$ onto $\calc $.
\smallskip
\hfill{q.e.d.}

Let us turn to the proof of the Theorem. We have a holomorphic
family $\calc \to \Delta $ of  stable curves $(C_{\lambda },
v_{\lambda })$
over $X$ such that:

1) $(C_{\lambda_0},v_{\lambda_0})=(C_0,v_0)$;

2) $h=f\circ \calv $ extends meromorphically onto $\calc $, where
$\calv :\calc \to X$ is an evaluation map.

\smallskip
We whant to lift this family into the envelope $(\hat U_Y, \pi )$. Denote
by $E_1,...,E_l$ the set of all ireducible curves in $\calc $
which are contracted
by $\calv :\calc \to X$ to the points. $\calc \setminus
\bigcup_{j=1}^lE_j$ clearly
lifts to the envelope by the Cartan-Thullen construction. Note further that
$E_j$
do not intersect $\bigcup_{j=1}^mU_j$. Thats why either $E_j\subset \hat C_0$
 - is a compact component of singular fiber, or projects
surjectively
on $\Delta $ under the projection $\pi :\calc \to \Delta $. In
the second case it
intersects $C_0$ and th point $u(C_0 \cap E_j)$ lifts to the
envelope. This proves that $\calc \setminus \hat C_0$ lifts to the
envelope. More precisely, we had shown that $\calc \setminus
\hat C_0^{'}$ lift into the envelope, where $\hat C_0^{'}$
is a union of all compact
components of the singular fiber maped by $\calv $ into the points.

Denote by $\hat C$ some connected subset of $\hat C_0^{'}$
and take its
neighborhood $V$, which do not intersect other components of $\hat C_0^{'}$.
Then $\calv $ maps $V$ onto the neighborhood of the point $\calv (\hat C$ in
$X$. And now is clear that this point lifts to the envelope.

Thus we had proved that $(C_0,v_0)$ lifts to the envelope $(\hat U_Y,\pi )$
of $U$. This implies that $(C_n,u_n)$ should lie in a compct subset of
$\hat U_Y$. Contradiction.
\smallskip
\hfill{q.e.d.}
\smallskip\rm
Let us give one corollary of the Continuity principle just proved.
\smallskip\noindent\bf
Corollary 5.1.8. \it Let $X$ is a complex surface with one singular
normal point $p$. Let $D$ be a domain in $X$, $\d D\ni p$. Suppose there
is a sequence $(C_n,u_n)$ of stable curves over $X$ converging to
$(C_0,u_0)$ in Gromov topology and such that

a) there is a compact $K\subset D$ with $u_n(\d C_n)\subset K$ for all
$n$;

b) $p\in u_0(C_0)$.

\noindent
Then every meromorphic function from $D$ extends to the neighborhood of
$p$.
\smallskip\rm
To our nowledge this statement is new also for holomorphic functions.

\smallskip\noindent
\sl 5.2. Construction of the family of curves.

\rm
Now we are ready to construct the family of curves needed for the step 2
in the proof of the {\sl Theorem 1}.

Let $\{ J_t\}_{t\in [0,1]}$ be a~family of $C^1$-smooth almost
complex structures on the~$4$-manifold $X$, depending $C^1$-continuously
on $t$, such that all $J_t$ are tamed by a fixed symplectic form $\omega$.
We lift the structures $J_t$ onto the envelope $(\hat U_Y,\pi )$ in such
a way that thouse liftings are standart on $\hat U_Y\setminus i(U_1)$,
where $i:U\to \hat U_Y$ is a canonical imbedding.

Let for $t=1$ a $J_1$-holomorphic sphere is given. We suppose that
$c_1(X)[M]-\delta_1 -\varkappa_1=p \ge1$. Here $\delta_1$ is the~geometric
self-intersection of $M_1$ and $\varkappa_1$ is the~sum of conductors of
all cusp-points of $M_1$. From {\sl Theorem 2} we get a family $M_t=f(J_t)$
of $J_t$-holomorphic
curves for $t$ sufficiently close to $1$, where $f:V\subset {\cal J} \to
{\cal S}\times {\cal J}_S \times {\cal J}$ is
a~local section of $\pr_{\cal J}$ constructed in {\sl Theorem 2}.

\smallskip\noindent\bf
Definition 5.2.1 \sl We say that the~family $\{M_t\}_{t\in(\hat t,1]}$ of
(possibly reducible) $J_t$-holomorphic curves is semi-continuous if
there exist a (may be infinite) partition of the~ interval $(\hat t,1]$
of the form $1=t_0 >t_1>\ldots>t_n>\ldots$ with $t_n\searrow \hat t$,
natural numbers $1=N_0 \le N_1 \le\ldots \le N_n \le\ldots$,
$J_t$-holomorphic maps $u_t: \bigsqcup_{j=1}^{N_i} S^j_i\to X$
for $(t_{i+1}, t_i]$ with $M_t\deff u_t(\bigsqcup_{j=1}^{N_i} S^j_i)$,
such that

 \sli $u_t: (t_{i+1}, t_i] \times \bigsqcup_{j=1}^{N_i} S^j_i\to X$
is a continuous map;

 \slii ${\ss area}(M_t)$ are uniformly bounded from above;

 \sliii ${\cal H}$-$\lim_{t\searrow t_{i+1}} M_t \deff \overline
M_{t_{i+1}} \supset M_{t_{i+1}}$.

 The~inclusion $\overline M_{t_{i+1}} \supset M_{t_{i+1}}$ means that
$M_{t_{i+1}}$ has no other components that those of $\overline M_{t_{i+1}}$.

 We say that $\{M_t\}$ is a family of spheres if all $S^j_i$ are spheres.

\smallskip\rm
Let $T$ be the~infimum of such $\hat t$, for which there is
a semi-continuous family $\{M_t\}_{t\in (\hat t, 1]}$ of spheres, such that
for all irreducible components $M^1_t,\ldots M^{N_i}_t$ of $M_t$,
$t\in (t_{i+1}, t_i]$, one has

 {\sl a)} $c_1(X)[M^j_t] -\delta_t^j -\varkappa_t^j=p_t^j \ge1$;

 {\sl b)} $\sum_{j=1}^{N_i} p_t^j \ge p$.

We allow existence of multiple components, {\sl i.e.} that some of $M_t^j$
can coincide.

\smallskip
The set $T$ is obviously open. To prove the~closeness of $T$ we note that
since all $J_t$ are tamed by the~same form $\omega$, the~areas of
$J_t$-holomorphic curve are uniformly in $t$ estimated from above and below by
$\int_{M_t} \omega$. More over, since $\overline M_{t_{i+1}} \supset
M_{t_{i+1}}$ we obtain that $\int_{M_{t_i}} \omega \le \int_{M_{t_{i+1}}}
\omega$. In particular, this implies that the~sequence $\{N_i\}$ is bounded
from above, and hence stabilizes for $i$ big enough. From
Gromov's compactness theorem and disk-convexity of the envelope
(Theorem 5.1.3)we obtain that for every $j={1,\ldots,N_i}$ the
sequence $\{M^j_i\}_{i=1}^\infty$ has a subsequence, still be denoted in
the~same way, which converges to a $J_{\hat t}$-holomorphic curve $\overline
M^j\subset \hat U_Y$.

To simplify the~notations, we drop the~upper index $j$ from now on,
and write $\overline M$ instead $\overline M^j$, having in mind that we can
do the~same constructions for all $j=1,\ldots,N_i$. Note that all irreducible
components of $\overline M$ are $J_{\hat t}$-holomorphic spheres. We write
$\overline M = \bigcup_{k=1}^d m_k\cdot M_k$, where $M_k$ are distinct
irreducible components of $\overline M$ with the~multiplicities $m_j\ge1$.
The genus formula for $\overline M$ have now (because of multiplicities!)
the~following form:

\smallskip
$$
0=\msmall{[\overline M]^2 - c_1(X)[\overline M] \over2}
+ \sum_{k=1}^d m_k - \sum_{k=1}^d m_k(\delta_k +\varkappa_k)
-\sum_{k < l}m_k\cdot m_l [M_k]\cdot [M_l]
$$
$$
-\sum_{k=1}^d \msmall{m_k(m_k -1) \over2 } [M_k]^2.
\eqno(5.2.1)
$$
\smallskip
The~formula can be obtained by taking the~sum of genus formulas for each $M_k$
and then completing the sum $\sum_{k=1}^d m_k[M_k]^2$ to $[\overline M]^2$.
Here $\delta_k$ and $\varkappa_k$ denote self-intersection and conductor
numbers
of $M_k$. We have also
$$
\varkappa^j_i + \delta^j_i = \msmall{ [M^j_i]^2 - c_1(X)[M^j_i] \over2} +1.
\eqno(5.2.2)
$$
Note that $[\overline M]^2 = [M^j_i]^2$ and $c_1(X)[\overline M]=
c_1(X)[M^j_i]$ for $i$ sufficiently big, and that $\sum_{k<l} m_k m_l\,
[M_k][M_l] \ge \sum_{k=1}^d m_k -1$ because the~system of curves
$\bigcup_{k=1}^d M_k$ is connected. From (5.2.1) and (5.2.2) we get
$$
\sum_{k=1}^d m_k(\delta_k + \varkappa_k) +
\sum_{k=1}^d \msmall{m_k(m_k-1) \over2 } [M_k]^2
\le \delta^j_i + \varkappa^j_i
$$
and
$$
\sum_{k=1}^d m_k c_1(X) [M_k] = c_1(X)[M^j_i].
$$
Using positivity of $(X,\omega)$, {\sl i.e.} the~property $[M_k]^2 \ge0$,
we obtain
$$
\sum_{k=1}^d m_k \biggl( c_1(X) [M_k] -\delta_k -\varkappa_k\biggr)
\ge c_1(X)[M^j_i] -\delta^j_i -\varkappa^j_i =p_j.
$$
Thus we can choose among $M_1,\ldots,M_d$ the~subset with properties {\sl a)}
and {\sl b)}.

Returning to the~old notations, we can ``correct'' all $\overline M^j$,
choosing an appropriate subset $M^j \subset \overline M^j$. Now applying {\sl
Theorem 2} to all $M^j$, we can extend our family $\{M_t\}$ in a~neighbourhood
of $\hat t$, and hence define it for all $t\in [0,1]$.

Using continuity principle (see {\sl Corollary 5.1.3}) along each continuous
piece of our semi-continuous family, we get the~next statement, which is
more general then  {\sl Theorem 1}.

\smallskip\noindent
\bf Theorem 5.2.1. \it Let $M$ be a $\omega $-positive immersed  two-sphere
in a
positive disk-convex K\"ahler surface $(X,\omega )$, having only positive
self-intersections. Suppose that $c_1(X)\cdot [M] - \delta =p\ge 1$. Then the
envelope of meromorphy $(\hat U_Y,\pi )$ of any neighbourhood $U$ of $M$
 relative to any disk-convex K\"ahler space $Y$ contains one-dimensional
 compact analytic set $C$ such
that:

\smallskip
$(1)$ $\pi (C) = \bigcup _{k=1}^NC_k$ is a union of rational curves.

\smallskip
$(2)$ $\displaystyle \sum _{k=1}^N \left(c_1(X)\cdot [C_k]
-\delta_k -\varkappa_k \right) \ge p$.
\smallskip
\rm The next rigidity property of symplectic imbeddings is a straitforward
corollary from this theorem.
\smallskip\noindent\bf
Corollary 5.2.2. \it Let $M_1$ and $M_2$ two symplectically imbedded spheres
in $\cc \pp^2$. Then any biholomorphism of the neighbourhood of $M_1$ onto the
neighbourhood of $M_2$ is fractional linear.
\rm

\bigskip\noindent
\sl 5.3. Examples.

\nobreak\smallskip\rm
Here we discuss a few more examples concerning the~envelopes of meromorphy.

\state Example 1. Let $(X,\omega) = (\cc\pp^1 \times \cc\pp^1,\omega_{FS}
\oplus \omega_{FS} )$ where $\omega_{FS}$ denote the~Fubini-Studi metric on
$\cc\pp^1$. Note that $c_1(X)= 2[\omega]$. Let $J$ be an~$\omega$-tame
almost complex structure on $X$ and $C$ be a $J$-holomorphic curve on $X$.
Denote by $e_1$ and $e_2$ the~standard generators  of $\ssh_2(X,\zz)= \zz^2$
and write $[C]=a\cdot e_1 + b\cdot e_2$. Then we get $a+b=\int_C\omega \ge 1$
and $c_1(X)[C]= 2(a+b)$. Further, by Genus formula $0\le g(C) + \delta
+\varkappa = (2ab -2(a+b))/2 +1 =(a-1)(b-1)$. Thus we conclude that both $a$
and $b$ are nonnegative and $[C]^2=2ab \ge0$. So $\cc\pp^1 \times \cc\pp^1$
is nonnegative in our sense.

Let $M$ be imbedded symplectic sphere in $X$. Then $(a-1)(b-1)=0$ by Genus
formula. So we can assume that $a=1$ and $b\ge0$. Now one concludes that
the~following holds:

\state Corollary 5.3.1. \it Let $M$ be a~symplectic sphere in $(\cc\pp^1
\times \cc\pp^1,\omega_{FS} \oplus \omega_{FS})$. Then the~envelope of
meromorphy of any neighbourhood of $M$ contains a graph of a rational
map of degree $0\le d \le b$ from $\cc\pp^1$ to $\cc\pp^1$.

\rm
\smallskip
\state Example 2. $(X,\omega) = (\cc\pp^2,\omega_{FS})$. Let $M$ be
a~symplectic surface in $X$ of degree $m\deff \int_M \omega$ with {\sl
positive} self-intersections. Then obviously $c_1(X)[M]=3m$.
Remark, that we proceed the~construction of
a~family $\{M_t\}$ if the~condition $c_1(X)[M^j_t] > \varkappa(M^j_t)$ is
satisfied for all irreducible components $M^j_t$ of $M_t$. By the~Genus
formula one has $\varkappa(C) \le (d-1)(d-2)/2$ for every pseudo-holomorphic
curve $C$ of degree $d$. So we can proceed if $3m > (m-1)(m-2)/2$, which is
equivalent to $1\le m \le8$. Thus we have the~following

\state Corollary 5.3.2. \it Let $M$ be a~symplectic surface in $\cc\pp^2$ of
degree
$m\le 8$ with positive self-intersections. Then the~envelope of
meromorphy of any neighbourhood of $M$ coincide with $\cc\pp^2$ itself.

\state Remark. \rm Note that the~examples includes all imbedded symplectic
surfaces in $\cc\pp^2$ of genus $g\le 21$.
\smallskip\noindent\bf
Example 3. \rm Let $X$ is a ball in $\cc^2$, $w$ is a standart euclidean
form. Blow up the origin in $\cc^2$ and denote by $E$ the exeptional curve.
 By $\hat X$ denote the blown-up ball $X$.
Blown-up $\cc^2$ is also K\"ahler, denote by $w_0$ some K\"ahler form there.
Consider a sufficiently small $C^1$-perturbation of $E$. This will be a
$w_0$-symplectic sphere in $\hat X$, denote it by $M$. Chern class of the
normal bundle to $M$ is equal to that of for $E$ and thus is $-1$. So
$c_1(\hat X)[M]=1$ and our Theorem 1 applies. $E$ is the only rational
curve in $\hat X$, so we get that

\it the envelope of meromorphy of any neighborhood of $M$ contains $E$.

\rm One can then blow down the picture to obtain downstears a sphere
$M_1$-image of $M$ under the blown-down map. This $M_1$ is homologous
to zero, so cannot be symplectic, and for this $M_1$ our Theorem 1
cannot be apllied.
\medskip\noindent
\bf 6. Appendix.

\smallskip\noindent
\sl 6.1. A complex structure lemma.

\rm For the convenience of a~reader we give here the following (well known for
the smooth case) lemma, which the line bundles case can be found in
[Hf-L-Sk].

\smallskip
{\bf Lemma 6.1.1.} {\it Let $S$ be a Riemannian surface with a complex
structure $J_S$ and $E$ a $L^{1,p}$-smooth complex vector bundle of rank
$r$ over $S$. Let also
$$
\dbar_E : L^{1,p}(S, E) \to
L^p(S,\Lambda^{(0, 1)}S \otimes E)
$$
be a differential operator, satisfying the condition
$$
\dbar_E(f\xi) = \dbar_S f \otimes \xi
+ f\cdot \dbar_E\xi,  \eqno(6.1.1)
$$
where $\dbar_S$ is the Cauchy-Riemann operator, associated to $J_S$.
Then the sheaf
$$
U\subset S \mapsto {\cal O}(E)(U) := \{\,\xi \in  L^{1,p}(U, E) \, :\,
\dbar_E\xi=0 \,\}   \eqno(6.1.2)
$$
is analytic and locally free of rank $r$. This defines the holomorphic structure
on $E$, for which $\dbar_E$ is an associated  Cauchy-Riemann operator.
}

\smallskip\noindent
{\bf Remark.} The condition $(6.1.1)$ means that $\dbar_E$
of order 1 and has the Cauchy-Riemann symbol.

\smallskip\noindent
{\bf Proof.} It is easy to see that the problem is essentially local. So we
may assume that $S$ is a unit disk $\Delta$ with the standard complex
structure. Let $\xi=(\xi_1,\ldots,\xi_r)$ be a $L^{1,p}$-frame of $E$
over $\Delta$. Let also $\Gamma \in L^p(\Delta,\Lambda^{0, 1}\Delta
\otimes {\bf Mat}(r,\cc))$ be defined by relation
$$
\dbar_E \xi_i = \sum_j \Gamma_i^j \xi_j,
\qquad \hbox{or in matrix form, } \qquad
\dbar_E \xi =  \xi \cdot \Gamma.
\eqno(6.1.3)
$$
Then for any section $\eta=\sum g^i \xi_i$ the equation $\dbar_E
\eta =0$ is equivalent to
$$
\dbar g^i + \sum_j \Gamma^i_j g^j=0,
\qquad \hbox{or in matrix form, } \qquad
\dbar g + \Gamma \cdot g =0.
\eqno(6.1.4)
$$
In particular, a frame $\eta=(\eta_1,\ldots,\eta_r)$ with $\eta_i=
\sum_j g^j_i \xi_j$ consists of sections of ${\cal O}(E)$ {\sl iff}
$$
\dbar g^i_k + \sum_j \Gamma^i_j g^j_k=0,
\qquad \hbox{or in matrix form, } \qquad
\dbar g + \Gamma \cdot g =0.
\eqno(6.1.5)
$$

\smallskip
 Let the map $\tau_t : \Delta \to \Delta$ be defined by
formula $\tau_t(z) = t\cdot z$, $0<t<1$. One can easily check that
$$
\Vert \tau_t^*\Gamma \Vert_{L^p(\Delta)}
\le t^{1-2/p} \Vert \Gamma \Vert_{L^p(\Delta)}.
$$
So taking pull-backs $\tau_t^*E$ with $t$ sufficiently small we
may assume that $\Vert \Gamma \Vert_{L^p(\Delta)}$ is small
enough. Now consider the mapping $F$ from $L^{1,p}(\Delta,
{\bf Gl}(r,\cc)) \subset L^{1,p}(\Delta,
{\bf Mat}(r,\cc))$ to $L^p(\Delta,\Lambda^{0, 1}\Delta
\otimes {\bf Mat}(r,\cc))$ defined by formula $F(g):=
\dbar g \cdot g^{-1}$. It is easy to see that
the derivation of $F$ in $g\equiv \id$ equals to $\dbar$.
Due to the {\sl Lemma 3.2.1}  and the implicit function theorem
for any $\Gamma$ with $\Vert \Gamma \Vert_{L^p(\Delta)}$ small
enough there exists $g$ in $L^{1,p}(\Delta, {\bf Gl}(r,\cc))$ with
$g(0)=\id$,
satisfying the equality $F(g)=-\Gamma$, which is equivalent to
(6.1.5). Consequently, in the neighbourhood of every point $p\in S$ there
exists a frame $\xi=(\xi_1,\ldots,\xi_r)$ of $E$ consisting of section
of ${\cal O}(E)$. Further, due to (6.1.2) a local section $\eta=\sum g^i
\xi_i$ of $E$ is a section of ${\cal O}(E)$ {\sl iff} $g^i$ are holomorphic.
This implies that ${\cal O}(E)$ is analytic  and locally free of rank $r$.
\qed

\bigskip\noindent
\sl 6.2. A matching structures lemma.

\smallskip\nobreak
\rm Let $B(r)$ be a ball of radius $r$ in $\rr^4$ centered at zero, and $J$
a $C^1$-smooth almost-complex structure on $B(2)$, $J(0)=J\st$.
Further, let $M = u(\Delta )$ be a closed primitive $J$-holomorphic disk
in $B(2)$ such, that $u(0)=0$ and $M$ transversally meets $S^3_r$ for
$r\ge 1/2$. Here $S^3_r = \d B(r) $ and transversality is understood with
respect to both $TS^3_r$ and $F_r$, see {\sl paragraph 4.1.}

By $B(r_1,r_2)$ we shall denote the spherical shell $\{ x\in \rr^4:r_1<\Vert
x\Vert <r_2 \}$. In the lemma below denote by $D_{\delta }$ the preimage of
$B(1+\delta )$ by $u$.

\smallskip\noindent\bf
Lemma 6.2.1. \it For any positive $\delta >0$ there exists an $\epsi >0$ such
that if an almost complex structure $\tilde J$ in $B(1+\delta )$ and a closed
$\tilde J$-holomorphic curve $\tilde u : D_{\delta }\to B(1+\delta )$ satisfy
$\Vert \tilde J - J\Vert_{C^1(\bar B(1+\delta ))} < \epsi $ and
$\Vert \tilde u - u\Vert_{L^{1,p}(\bar D_{\delta })} < \epsi $, then
there exists an almost-complex structure $J_1$ in $B(2)$ and $J_1$-holomorphic
 disk $M_1$ in $B(2)$ such that:

a) $J_1\ogran_{B(1-\delta )} = \tilde J\ogran_{B(1-\delta )}$ and
$J_1\ogran_{B(1+
\delta ,2)} = J\ogran_{B(1+\delta ,2)}$.

b) $M_1\ogran_{B(1-\delta )} = \tilde u(D_{\delta })\cap B(1-\delta )$ and
$M_1\cap B(1+\delta ,2) = M\cap B(1+\delta ,2)$.
\smallskip
\noindent\bf
Proof. \rm We have chosen the parametrization of $M$ to be primitive. So $u$
is an imbedding on $D_{-\delta ,\delta } = u^{-1}(B_{1-\delta ,1+\delta })$.
Let us identify a neighbourhood $V$ of $u(D_{-\delta ,\delta })$ in $B_{1-
\delta ,1+\delta }$ with the neighbourhood of zero-section in the normal bundle
$N$ to $u(D_{-\delta ,\delta })$. Now the $\tilde u\ogran_{D_{-\delta ,\delta
}}$
can be viewed as a section of $N$ over $u(D_{-\delta ,\delta })$ which is small
{\sl i.e.}
contained in $V$. Using appropriate smooth function $\phi $ on
$D_{-\delta ,\delta }$ (or equivalently on $u(D_{-\delta ,\delta} )$), $\phi
\ogran_{B(1-\delta )\cap D_{\delta }}\equiv 1$, $\phi\ogran_{\partial D_{\delta
}
}\equiv 0$, $0\le \phi \le 1$ we can glue $\tilde u$ and $u$ to get \sl
symplectic \rm surface $M_1$ which satisfies (b).

Patching $J$ and $\tilde J$ and making simultaneously $M_1$ pseudo-holomorphic
can be done in the way we did in the proof of {\sl Lemma 1.1.2}.
\qed





\bigskip\noindent\sl
6.3. Gromov topology and versal deformations of noncompact curves.
\smallskip\rm

\smallskip\nobreak
\rm

\def\cala{{\cal A}}
\def\barr#1{{\overline{#1}}}

We start with giving several definitions and notations needed.
Let $\Delta^2\deff\{\,(z_1,z_2)\in\cc^2\;:\; |z_1|<1, |z_2| <1\;\}$ be
a bidisk. For $\lambda \in \Delta$ we set $\cala_\lambda \deff
\{\,(z_1,z_2)\in\Delta^2\;:\; z_1\cdot z_2 =\lambda\;\}$.
$\cala_0$ is a union of two disks which we denote by
$\cala_0^1$ and $\cala_0^2$ respectively; the function $z_i$ is a holomorphic
coordinate on $\cala_0^i$; the only common point of $\cala_0^i$ is zero.
For $\lambda\not=0$ $\cala_{\lambda }$ is an annulus with a boundary
consisting
of a pair of circles $\gamma_1$ and
$\gamma_2$, equipped with a pair of holomorpfic coordinates $z_1$ and
$z_2$, such that $|z_i|$ restricted to $\gamma_i$ is identically 1.

If $\lambda \lrar 0$, these annuli degenerate into a node, and every function
$z_i$ becomes a coordinate on the correspondind irreducible component of
$\cala_0$. The family $\cala\deff \sqcup_{\lambda\in \Delta} \cala_\lambda$
is a holomorphic deformation of $\cala_0$. Denote $\cala^*\deff
\sqcup_{\lambda\in \Delta^*} \cala_\lambda$.

We call $\cala_0$ the \sl standard node, \rm and the point $(0,0) \in
\cala_0$ the \sl nodal point. \rm It is the only singular point of
$\cala_0$, which is an ordinary double point.

\medskip\noindent
\bf Definition 6.3.1. \it A \sl nodal curve \it is a connected purely
1-dimensional complex analytic space $C$, possibly with boundary,
satisfying the following conditions:

\sli the only singularities of $C$ are finitely many nodal points;

\slii the boundary $\d C$ consists of finitely many smooth circles $\gamma_i
\cong S^1 = \{\,z\in\cc\,:\, |z|=1\,\}$;

\sliii the set $\barr C \deff C \cup \d C$ is compact. \rm

\smallskip
In particular, $C$ is smooth in the neibourhood of $\d C$. Note that
$C$ can consists of several irreducible components. Nodal curves appear
in a natural way as a limit of immersed smooth curves with respect to
the Gromov topology, which we describe below.
Informally speaking, every node in the limit
arise as a result of degeneration  of certain annuli on prelimit curves.
Topologically, it corresponds to contraction of certain circles on prelimit
curves to a nodal point on the limit curve. This explains the next

\smallskip\noindent
\bf Definition 6.3.2. \it We shall say that connected, oriented real
surface
with boundary $(\Sigma ,\d \Sigma )$ parametrizes a complex nodal curve $C$
if there is a continuous map $\sigma :\Sigma \to C$ such that:

\it $1)$ if $a\in C$ is a nodal point then $\gamma_a = \sigma^{-1}(a)$ is
a smooth imbedded circle in $\Sigma \setminus \d \Sigma $, and if $a\not= b$,
then $\gamma_a\cap \gamma_b=\emptyset $;

$2)$ $\sigma :\Sigma \setminus \bigcup_{i=1}^N\gamma_{a_i}\to C\setminus
\{ a_1,\ldots ,a_N\} $ is a diffeomorphism. Here $a_1,\ldots,a_N$ are
nodes of $C$.

\smallskip \rm Roughly speaking, $\Sigma$ is obtained from $C$ by replacing
all nodes of $C$ by corresponding cylinders, the procedure inverse to
contracting circles into nodal points. Note that $\Sigma$ is defined by $C$
uniquely up to diffeomorphism, but $\sigma $ is not determined by $C$.
However,
two parametrisations $\sigma :\Sigma \to C$ and $\sigma ':\Sigma' \to C$
differ by  a diffeomorphism $f:\Sigma \to :\Sigma'$. Denote by $\tau :\hat C
\to C$ a normalisation of C.

\smallskip\noindent\bf
Definition 6.3.3. \it A stable curve over an almost complex manifold
$(X,J)$ is a pair $(C,u)$, where $C$ is a connected nodal curve, possibly
with boundary $\d C=\bigcup_{i=1}^d\gamma_i$, and $u:C\to X$ is a holomorphic
map which satisfies the following condition:

1) if $C_j$ is an irreducible component of $C$ biholomorphic to $\pp^1$ and $u(C_j)
=\{ point \} $, then $\tau ^{-1}(C_j)$ contains at least three points, which
are the preimages of the nodes of $C$;

2) if $C_j$ is a torus and again $u(C_j)=\{ point \} $ then $\tau ^{-1}(C_j)$
containes at least one praimage of the node of $C$.

\smallskip\rm
Now we are going to describe the Gromov topology on the space of stable curves
over $X$. Let a sequence $J_n$ of $C^1$-smooth almost complex structures on $X$ is
given, and suppose that $J_n\to J$ in $C^1$-topology. Further, let a sequence
$(C_n,u_n)$ of stable curves over $(X,J_n)$ is given, such that all $C_n$ are
parametrized by the same real surface $\Sigma $.

\smallskip\noindent\bf
Definition 6.3.4. \it  We say that $(C_n,u_n)$ converge to a stable curve $(C_\infty, u_\infty)$
in \sl Gromov topology, \it if all $C_n$ and $C_\infty$ are parametrised by
the same real surface $\Sigma$ and there exist parametrisations
$\sigma_n:\Sigma \to C_n$ and $\sigma_\infty:\Sigma \to C_\infty$
such that

\sli $u_n \scirc \sigma_n$ converges to $u_\infty \scirc \sigma_\infty$ in
$C^0(\barr\Sigma, X)$-topology;

\slii if $\{a_i\}$ is the set of nodes of $C_\infty$ and $\gamma_i \deff
g_\infty\inv(a_i)$ are the corresponding circles in $\Sigma$, then on any
compact $K\Subset \Sigma \bss \cup_i \gamma_i$ the convergence $u_n \scirc
\sigma_n
\lrar u_\infty \scirc \sigma_\infty$ is in fact smooth;

\sliii on any compact $K\Subset \Sigma \bss \cup_i \gamma_i$
the complex structures $\sigma_n^* J_n$ converge to the complex structure
$\sigma_\infty^* J_\infty$.
\smallskip\rm
The following result is a Gromov compactness theorem, see [Gro].

\medskip\noindent
\bf Theorem. \it Suppose that a sequence $(C_n,u_n)$ of stable over $X$
$J_n$-holomorphic curves is given such that:

a) $J_n$ are of class $C^1$ and are $C^1$-converging to $J$;

b) ${\sl area}[u_n(C_n)]\le M$ for all $n$;

c) $u_n$ converge near the boundary, i.e. for any boundary component $\gamma$
of $\Sigma $ and any $n$ there is a holomorphic imbedding $\phi :A_r
\deff \{ z\in \cc : r<\vert z
\vert <1\} \to C_n$ with $\sigma_n^{-1}\circ \phi_n(\{ |z|=1\} ) =\gamma$,
$r$ not depending on $n$, such that $u_n\phi_n $ converges on $A_r$ to some
$J$-holomorphic map $h_{\infty }$.

Then there is a subsequence $(C_{n_k},u_{n_k})$ which converges to a stable
over
$X$ $J$ - holomorphic curve $(C_{\infty },u_{\infty })$. Moreover for each
boundary component $\gamma $ there is an imbedding $\phi_{\infty }:A_r\to
C_\infty$ such that $u_n\circ \phi_n \to u_\infty\circ \phi_\infty =
h_\infty$ on $A_r$.

\smallskip\rm
Our aim is obtain a nice description of the neighborhood in Gromov topology
of a stable curve over a complex manifold. This description will allow
us to draw analytic families through two sufficiently close curves
({\sl Proposition 5.1.1}). This will be possible, due to the Theorem of J.-P.
Ramis, if one can prove that the moduli space of stable curves over a
complex manifold forms a Banach analytic space of finite codimension.

Let $B$ is an open set in some complex Banach space .
\smallskip\noindent\bf
Definition 6.3.5. \it We say that a closed subset $\calm \subset B$
is a banach analytic set of finite codimension if there is an analytic
map $F:B\to \cc^N$, for some $N$, such that $\calm =\{ x\in B: F(x)=0\} $.
\smallskip\tenrm
In general nothing good can be said aboud Banach analytic sets. Every compact,
ex. interval $[0,1]$, or a converging sequence of points, can be realised as
a Banach analytic set in appropriate Banach space, see [Ra] p. 33. Hovewer,
the structure of
Banach analytic sets of finite codimension is pretty nice according to the
following theorem, which is due to Ramis.

\smallskip\noindent\bf
Theorem. \it Let $\calm \subset B$ be a Banach analytic set
of finite codimension, $0\in \calm $. Then there is a neighborhood $U\ni x_0$
such that $\calm \cap B$ is a finite union of irreducible components $\calm_j$,
each of them being a finite ramified cover of a neighborhood of zero in the
subspace $L_j$ of finite codimension.
\smallskip\tenrm
Fix now a complex manifold $X$ and a stable curve $(C_0,u_0)$ over $X$ parametrised
by $\Sigma $. Our aim in this section is the following
\smallskip\noindent
\bf Theorem 6.3.1. \it There exist a Banach analytic spaces of finite
codimension $\calm $ and $\calc $ and holomorphic maps
$$
\matrix{
\calc &  {\buildrel \calu \over \longrightarrow } & X
\cr
\llap{$\pi$}\downarrow
\cr
\calm
}
$$
\noindent
such that:

a) for any $\lambda \in \calm $ fiber $C_{\lambda }=\pi^{-1}(\lambda )$ is a
nodal curve
parametrised by $\Sigma$ and for some $\lambda_0$ $C_{\lambda_0}\sim C_0$;

b) $(C_\lambda, u_\lambda=\calu\mid_{C_\lambda },)$ is stable over $X$
and $u_{\lambda_0}=u_0$;

c) if $(C',u')$ is sufficiently close to $(C,f)$ in Gromov topology then
there exists $\lambda '\in \calm $ such that $(C',u')=(C_{\lambda '},F\mid_
{\lambda '})$.

\smallskip\noindent\sl
Proof. \rm Let $(C_0,u_0)$ is our stable curve over $X$. Denote by
$E$ the pull-back $u_0^*TX$ of the tangent budle to $X$ together with
natural holomorphic structure on it. We suppose that $u_0$ extends
smoothly onto the boundary of $C_0$ to be able to consider the
$L^{1,p}$-sections of $E$.


Cover $C_0$ by a finite family
of open sets $U_i$ with boundaries in such a way that intersections
$U_{ij}:=U_i\cap U_j$ have piecewisely smooth boundaries. In the case
when $U_i$ contains a nodal point it should be a union of two disks.
We take $U_i$ sufficiently small to find a coordinate chart $V_i$ containing
$u_0(U_i)$.

If $U_i$ is smooth consider a complex Banach manifold of holomorphic
mappings $u_i\in L^{1,p}(U_i,V_i)$ with a tangent space at $u_i^0=u_0\mid_
{U_i}$ equal to $\calh^{1,p}(U_i,E)$-space of holomorphic
$L^{1,p}$-sections of $E$ over $U_i$.

If $U_i$ is a neighborhood of a node, we  consider a complex Banach
manifold  of
holomorphic $L^{1,p}$-maps from $\cala $ to $V_i$ with  tangent space
at $u_i$ equal to $\calh^{1,p}(U_i,E)$, the space of pairs of
holomorphic $L^{1,p}$-sections over the components
of a standart node, which coincide at the origin.

Denote by $B_i$ open neighborhoods of $u_i$ in these Banach manifolds.
Repeat the same construction for $U_{ij}, i<j$, using as a coordinate
chart $V_i$, and get the Banach manifolds $L^{1,p}(U_{ij},V_i)$ with
tangent spaces $L^{1,p}(U_{ij},E)$.

Denote by $\phi_{ij}:V_j\to V_i$ the coordinate change. We can consider
the following analytic map between complex Banach manifolds
$$
\Phi:\Pi_iB_i\to \Pi_{i<j}L^{1,p}(U_{i,j},V_i),
$$
$$
(\Phi(\{ h_i\} ))_{ij}:=\phi_{ij}(h_j) - h_i.
$$
Zero level set of this mapping is a Banach analytic set
and is by construction some neighborhood in Gromov topology of $(C_0,u_0)$
in the space of stable complex curves over $X$. Denote this neighborhood
by $\calm $. Differential $d\Phi_{u_0}$ of this mapping at $u_0$
coincides with the differential of \v{C}ech complex
\def\calh{{\cal H}}
\def\1{{1\mkern-5mu{\rom l}}}
$$
\matrix{
\delta :& \sum_{i=1}^l \calh ^{1,p}(U_i, E)& \lrar&
   \sum_{i<j} \calh ^{1,p}(U_{ij}, E)
\cr
\delta :& (v_i)_{i=1}^l        & \longmapsto & (v_i -v_j)
}
$$
We shall make use now from the following lemma, which states the possibility
of solving a Cousin-type
problem with estimates.

\smallskip
\noindent
\bf Lemma 6.3.2. \it
Let $C$ be a nodal curve with piecewisely smooth boundary
and $E$ a holomorphic vector bundle over $C$.
Let $\{ U_i\}_{i=1}^l$ be a finite covering of $C$ by Stein domains
with piecewisely smooth boundaries. Let $2\le p<\infty$ and let
$\calh ^{1,p}(U_i, E)$ denote the space of sections of $E$,
which are holomorphic in $U_i$ and $L^{1,p}$-smooth till boundary $\d U_i$.
Set $U_{ij} \deff U_i \cap U_j$ and suppose that any triple intersection
$U_i \cap U_j \cap U_k$ with $i \not= j \not= k \not= i$ is empty.

Then the \v{C}ech complex
$$
\matrix{
\delta :& \sum_{i=1}^l \calh ^{1,p}(U_i, E)& \lrar&
   \sum_{i<j} \calh ^{1,p}(U_{ij}, E)
\cr
\delta :& (v_i)_{i=1}^l        & \longmapsto & (v_i -v_j)
}
\eqno(6.3.1)
$$
has the following properties:

\sli the image $\im(\delta)$ is of finite codimension and closed; more over,
$\coker(\delta) = \ssh^1(C, E)
\allowbreak
 = \ssh^1(C_{comp}, E)$, where $C_{comp}$
denotes the union of compact irreducible components of $C$;

\slii the kernel $\ker(\delta)$ admits a closed complementing.

\smallskip
\noindent\rm
Proof is left to the reader.

We can finish now the proof of {\sl Theorem 6.3.1}. Denote by $T=\im d\Phi_
{u_0}$ and by $S$ its finite dimensional complement. Let $\pi_T$ be a
projection onto $T$ parallel to $S$. Implicit function theorem applied to
$\pi_T\circ \Phi $ tells us that $\calm $ is contained in complex Banach
manifold $\caln $ with tangent space at $u_0$ equal to $\ker d\Phi =\ker
\delta =
\calh^{1,p}(C_0,E)$. Now our $\calm $ is a Banach analytic subset of
$\caln $ defined by the equation $\pi_S\circ \Phi =0$, where $\pi_S$ is
a projection onto finite dimensional vector space $S$ parallel to $T$.
Thus $\calm $ is of finite codimension.
\smallskip
\qed
\smallskip\noindent\sl
Proof of Proposition 5.1.1. \rm Let a sequence of stable curves over $X$
converges to $(C_0,u_0)$. Take a neighborhood $\calm $ of $(C_0,u_0)$ in
Gromov topology, which is realised as a Banach analytic set of finite
codimension in Banach ball ${\bf B}$. Denote by $\lambda_0$ the point on
$\calm $ which corresponds to $(C_0,u_0)$. By the theorem of Ramis $\calm $
in the neighborhood of $\lambda_0$ has finite number of irreducible
components. One of them, denote it as $\calm_1$ should contain a point
$\lambda_N$ which
coresponds to $(C_N,u_N)$. Further, there is a subspace $L$ of ${\bf B}$ of
finite codimension such that $\calm_1$ is a finite covering of $L\cap {\bf B}
$. Now one can easily find an analytic disk passing through $\lambda_0$
and $\lambda_N$.
\smallskip
\qed

\magnification=\magstep1
\spaceskip=4pt plus3.5pt minus 1.5pt
\xspaceskip=5pt plus4pt minus 2pt
\font\csc=cmcsc10
\font\tenmsb=msbm10
\def\rr{\hbox{\tenmsb R}}
\def\cc{\hbox{\tenmsb C}}
\newdimen\length
\newdimen\lleftskip
\lleftskip=2.5\parindent
\length=\hsize \advance\length-\lleftskip
\def\entry#1#2#3#4\par{\parshape=2  0pt  \hsize%
\lleftskip \length%
\noindent\hbox to \lleftskip%
{\bf[#1]\hfill}{\csc{#2 }}{\sl{#3}}#4%
\medskip
}
\ifx \twelvebf\undefined \font\twelvebf=cmbx12\fi

\bigskip\bigskip
\bigskip\bigskip
\centerline{\twelvebf References.}
\bigskip


\entry{Ar}{Aronszajn, N.:}{A unique continuation theorem for elliptic
differential equation or inequalities of the second order.} J.\ Math.\ Pures
Appl., {\bf36}, 235--339, (1957).

\entry{B-K}{Bedford, E., Klingenberg W.:}{On the~envelope of meromorphy of
a 2-sphere in $\cc^2$.}  J.\ Amer.\ Math.\ Soc., {\bf4}, 623-646, (1981).

\entry{B-P-V}{Barth W., Peters, C., Van de Ven, A.:}{Compact complex
surfaces.} Springer Verlag, (1984).

\entry{Bn}{Bennequin, D.:}{Entrelacement et \'equation de Pfaff.}
Ast\'erisque {\bf107--108}, 87--161 (1982).

\entry{Cn-Sp}{Chern, S.-S., Spanier, E.:}{A theorem on orientable surfaces
in four-dimensional space.} Comment. Math. Helv. {\bf25}, 205--209 (1951).

\entry{Ch-St}{Chirka, E., Stout L.:}{A Kontinuit\"atssatz.} Preprint.

\entry{E}{Eliashberg, Y.:}{Filling by holomorphic discs and its applications.}
London Math.\ Soc.\ Lecture Notes, {\bf151}, Geometry of low dimensional
manifolds, (1991).

\entry{Gra}{Grauert, H.:}{\"Uber Modifikationen und exzeptionelle analytische
Mengen.} Math.Ann. {\bf146}, 331--368 (1962).

\entry{Gr-Hr}{Griffiths, P., Harris, J.:}{Principle of algebraic geometry.}
John Wiley \& Sons, N.-Y., (1978).

\entry{Gro}{Gromov, M.:}{Pseudo holomorphic curves in symplectic
manifolds.} Invent.math. {\bf82}, 307--347 (1985).

\entry{Hv-Po}{Harvey, R., Polking, J.:}{Removable singularities of
solutions of partial differential equations.} Acta Math. {\bf125}, 209--226
(1970).

\entry{Hr-W}{Hartman, P., Winter, A.:}{On the~local behavior of solutions
of non-parabolic partial differential equations.} Amer.\ J.\ Math.,
{\bf75}, 449--476, (1953).

\entry{Hf}{Hofer, H.:}{Pseudoholomorphic curves in symplectizations with
applications to the Weinstein conjecture in dimension three.} Invent. Math.
{\bf114}, 515--563 (1993).

\entry{Hf-L-Sk}{Hofer, H., Lizan, V., Sikorav, J.-C.:}{On genericity for
holomorphic curves in 4-dimensional almost-complex manifolds.} Preprint,
(1994).

\entry{Iv-1}{Ivashkovich, S.:}{Extension of analytic objects by the method
of Cartan-Thullen.}In Proceedings of Conference on Complex Analysis and Math.
Phys., 53--61, Krasnojarsk, (1988).

\entry{Iv-2}{Ivashkovich, S.:}{The Hartogs-type extension theorem for
meromorphic maps into compact K\"ahler manifolds.} Invent. Math. {\bf109},
47--54 (1992).

\entry{Ko-No}{Kobayashi, S., Nomizu, K.:}{Foundations of differential
geometry.} Vol.II, Interscience Publishers, (1969).

\entry{Li}{Lichnerovicz, A:}{Global theory of connections and holonomy
groups.} Noordhoff Int. Publishing, Leiden, (1976).

\entry{Mi-Wh}{Micallef, M., White, B.:}{The structure of branch points in
minimal surfaces and in pseudoholomorphic curves.} Preprint, (1994).

\entry{McD-1}{McDuff, D.:}{The local behavior of holomorphic curves in
almost complex 4-manifolds.} J.Diff.Geom. {\bf 34}, 143-164 (1991).

\entry{McD-2}{McDuff, D.:}{Singularities of $J$-holomorphic curves in
almost complex  4-manifolds.} J. Geom. Anal. {\bf2}, 249--265 (1992).

\entry{McD-3}{McDuff, D.:}{Singularities and positivity of intersections of
$J$-holomorphic curves.} In \it ``Holomorphic curves in symplectic geometry.''
 \rm Ed. by M.~Audin and J.~Lafontaine, Birkh\"auser, (1994).

\entry{Mo}{Morrey, C.:}{Multiple integrals in the calculus of variations.}
Springer Verlag, (1966).

\entry{Ra}{Ramis J.-P.:}{Sous-ensembles analytiques d'une vari\'et\'e
banachique complexe.} Springer, Berlin (1970).

\entry{Rf}{Rolfsen, D.:}{Knots and links.} Publish or perish. {\bf N7} (1976).

\entry{Sk}{Sikorav, J.-C.:}{Some properties of holomorphic curves in almost
complex manifolds.} In \it ``Holomorphic curves in symplectic geometry .''
\rm Ed. by M.~Audin and J.~Lafontaine, Birkh\"auser, (1994).

\entry{Si-1}{Siu, Y.-T.:}{Extension of meromorphic maps into K\"ahler
manifolds.} Annals of Math. {\bf102}, 421--462 (1975).

\entry{Si-2}{Siu, Y.-T.:}{Every Stein subvariety admits a Stein
neighbourhood.} Invent.math. {\bf38}, 89--100 (1976).

\entry{Sh}{Shafarevich I.:}{Basic algebraic geometry. Second edition.}
 Springer-Verlag.

\entry{Rm}{Remmert, R.:}{Holomorphe und meromorphe Abbildungen komplexer
R\"aume.} Math. Ann. {\bf133}, 328--370 (1957).

\entry{Sh}{Shcherbina, N.:}{On the~polynomial hull of a two-dimensional
sphere in $\cc^2$.} Soviet Math.\ Doklady, {\bf49}, 629--632, (1991).

\entry{Wn}{Weinstein, A.:}{Lectures on symplectic manifolds.} CBMS. Reg.
Conf. Series, {\bf N29}, AMS(1977).

\entry{Wa}{Walker, R.:}{Algebraic curves.} Springer Verlag, (1978).

\bigskip
\bigskip
\settabs 2 \columns
\+ Institute of Appl. Prob. Mech. \& Math.
       & Institute of Appl. Prob. Mech. \& Math.\cr
\+ Ukrainian Acad. Sci.,
       & Ukrainian Acad. Sci.,\cr
\+ vul. Naukova 3b, 290053 L'viv
       &  vul. Naukova 3b, 290053 L'viv \cr
\+ Ukraine    & Ukraine \cr
\medskip
\+ U.F.R. de Maht\'ematiques
                          & Ruhr-Universit\"at Bochum \cr
\+ Universit\'e de Lille-I
       & Mathematisches Institut \cr

\+ Villeneuve d'Ascq Cedex
              & Universit\"atsstrasse 150 \cr
\+ 59655 France
                 & NA 4/67  44780 Bochum  Germany \cr
\+ ivachkov@gat.univ-lille1.fr
                    & sewa@cplx.ruhr-uni-bochum.de   \cr

\end